\documentclass[a4paper,12pt]{amsart}

\usepackage{setspace} 
\usepackage{array} 
\usepackage{enumerate} 
\usepackage[pdftex]{graphicx} 
\usepackage[svgnames,table,hyperref, dvipsnames]{xcolor}
\usepackage{adjustbox}

\usepackage{empheq} 
\usepackage{amsthm}
\usepackage{mathrsfs} 
\usepackage{amssymb}
\usepackage{csquotes}
\usepackage{comment}

\usepackage{graphicx} 
\usepackage{tikz}
\usetikzlibrary{shadows,arrows,patterns,cd,babel}
\usepackage{ytableau}
\usepackage[raiselinks, pageanchor, hyperindex, citebordercolor=green, linkbordercolor=cyan, urlbordercolor=cyan, filebordercolor=cyan, menubordercolor=cyan]{hyperref}

\DeclareMathOperator{\dime}{dim}
\DeclareMathOperator{\db}{D^b}
\DeclareMathOperator{\derived}{D}
\DeclareMathOperator{\dbstwo}{D^b_{S_2}}
\DeclareMathOperator{\dbs3}{D^b_{S_3}}
\DeclareMathOperator{\dbsn}{D^b_{S_n}}

\DeclareMathOperator{\coh}{\textbf{Coh}}
\DeclareMathOperator{\qcoh}{\textbf{Qcoh}}
\DeclareMathOperator{\category}{\textbf{C}}

\DeclareMathOperator{\chara}{char}

\DeclareMathOperator{\spec}{Spec}

\DeclareMathOperator{\mor}{Mor}
\DeclareMathOperator{\homo}{Hom}
\DeclareMathOperator{\enndo}{End}
\DeclareMathOperator{\id}{id}

\DeclareMathOperator{\hilb}{Hilb}
\DeclareMathOperator{\quot}{Quot}
\DeclareMathOperator{\degr}{deg}
\DeclareMathOperator{\kernel}{ker}
\DeclareMathOperator{\coker}{coker}
\DeclareMathOperator{\image}{im}

\DeclareMathOperator{\red}{red}

\DeclareMathOperator{\blow}{Bl}
\DeclareMathOperator{\ext}{Ext}
\DeclareMathOperator{\codim}{codim}
\DeclareMathOperator{\tor}{Tor}

\DeclareMathOperator{\GL}{GL}

\DeclareMathOperator{\rank}{rk}
\DeclareMathOperator{\supp}{supp}
\DeclareMathOperator{\length}{length}

\DeclareMathOperator{\picard}{Pic}
\DeclareMathOperator{\modd}{mod}

\theoremstyle{plain} 
\newtheorem{prop}{Proposition}[section] 
\newtheorem{lemma}[prop]{Lemma} 
\newtheorem{kor}[prop]{Corollary}
\newtheorem{thrm}[prop]{Theorem}

\theoremstyle{definition} 
\newtheorem{defi}[prop]{Definition}
\newtheorem{bsp}[prop]{Example}

\newtheorem{bem}[prop]{Remark}

\title[A Semi-orthogonal Sequence for the Hilbert Scheme of Three Points]{A Semi-orthogonal Sequence in the Derived Category of the Hilbert Scheme of Three Points}
\author{Erik Nikolov}
\address{Institute of Algebraic Geometry, Leibniz University Hannover, Welfengarten 1, 30167 Hannover, Germany.}
\email{nikolov@math.uni-hannover.de}
\date{April 19, 2024.}

\begin{document}

\begin{abstract}
For a smooth projective variety $X$ of dimension $d \geq 5$ over an algebraically closed field $k$ of characteristic zero, it is shown in this paper that the bounded derived category $\db(X^{[3]})$ of the Hilbert scheme of three points admits a semi-orthogonal sequence of length $\binom{d-3}{2}$.\ Each subcategory in this sequence is equivalent to $\db(X)$ and realized as the image of a Fourier--Mukai transform along a Grassmannian bundle $\mathbb{G} \rightarrow X$ parametrizing planar subschemes in $X^{[3]}$.\ The main ingredient in the proof is the computation of the normal bundle of $\mathbb{G}$ in $X^{[3]}$.\ An analogous result for generalized Kummer varieties is deduced at the end.
\end{abstract}

\maketitle

\textbf{Disclaimer.} This preprint has not undergone peer review or any post-submission improvements or corrections. The Version of Record of this article is published in \textit{Selecta Mathematica} and is available online at \href{https://doi.org/10.1007/s00029-026-01140-2}{doi:10.1007/s00029-026-01140-2}.

\section{Introduction}

The Hilbert scheme of $n$ points, denoted by $X^{[n]}$, parametrizes zero-dimensional subschemes of length $n$ of a given (quasi-)projective variety $X$ of dimension $d = \dime X$ over a field.\ It is smooth and connected whenever $d \in \lbrace 1,2 \rbrace$ or $n \in \lbrace 1,2,3\rbrace$ and $X$ is smooth and connected \cite[Sect.\ 7.2]{fgahilbert}.\ It is also projective whenever $X$ is (for $d \geq 1$ arbitrary).\par 
\noindent If $d = 1$ and thus $X = C$ is a smooth curve, there is not much choice for non-reduced zero-dimensional subschemes on $C$, wherefore $C^{[n]} \cong C^{(n)}$ is isomorphic to the $n$-th symmetric product of the curve.\ For every $d \geq 1$, the Hilbert scheme $X^{[2]}$ can be realized as the quotient of the blow-up $\blow_{\Delta}X^2$ by an action of the symmetric group $S_2$ \cite[Sect.\ 7.3]{fgahilbert}.\ Already $X^{[3]}$ does not admit such a handy description anymore and its geometry becomes more complicated starting from $d \geq 2$. \par
\noindent On the level of derived categories, there is a description of the equivariant derived category of $C^n$ in terms of symmetric powers of $C$ \cite[Thm.\ B]{polishchukvdb} which for $n = 3$ reads 
\begin{align} \label{curvecase}
    \dbs3(C^3) = \left \langle  \db(C^{(3)}), \db(C^2), \db(C) \right \rangle.
\end{align}
Considering the derived category of the Hilbert scheme of a smooth surface $S$, the derived McKay correspondence in combination with a result of Haiman (\cite{bkr}, \cite{haiman}) yields the equivalence of derived categories
\begin{align} \label{surfacecase}
    \dbsn(S^n) \cong \db(S^{[n]}).
\end{align}
For $n=2$ and in higher dimensions $d$, there is a fully faithful embedding of derived categories $\dbstwo(X^2) \hookrightarrow \db(X^{[2]})$ as well as a semi-orthogonal decomposition
\begin{align*}
    \db(X^{[2]}) = \left \langle  (d-2)\cdot \db(X), \dbstwo(X^2) \right\rangle
\end{align*}
due to Krug, Ploog and Sosna \cite[Thm.\ 4.1(ii)]{krugploogsosna}.\ The notation indicates that there are $d-2$ fully faithful functors $\db(X) \hookrightarrow \db(X^{[2]})$ with semi-orthogonal essential images.\ In the spirit of Fantechi and G\"ottsche's description of the cohomology ring $H^{\ast}(X^{[3]},\mathbb{C})$ in \cite{cohomologyx3}, a description of $\db(X^{[3]})$ in terms of $\db(X), \db(X^2)$ and an equivariant part is desirable.\ A possible conjecture (generalizing (\ref{surfacecase})) is that
\begin{align} \label{conjecture}
    \db(X^{[3]}) \overset{?}{=} \left \langle \dbs3(X^3), (d-2)\cdot \db(X^2), \left(\binom{d}{2} - 1\right) \cdot \db(X) \right\rangle, \quad d \geq 2.
\end{align}
\noindent The expected multiplicities in (\ref{conjecture}) are consistent with \cite[Thm.\ 2.5.18]{goettsche1}.\ If $d = 1$, moving factors with negative sign to the other side of (\ref{conjecture}) also leads back to (\ref{curvecase}).\ This work studies the presumably easiest $\binom{d-3}{2}$ components of a potential semi-orthogonal decomposition of $\db(X^{[3]})$ like the one in (\ref{conjecture}), arising here from a correspondence of the form $X \leftarrow \mathbb{G} \rightarrow X^{[3]}$.\ To be more precise, the focus lies on the locus of zero-dimensional subschemes having support at single points $P \in X$.\ This locus comes along with a projection to $X$ whose fibres are (independently of $P$) precisely the \emph{Local} Punctual Hilbert schemes $\hilb_d^3$ (studied e.g.\ in  \cite{bertin}, \cite{goettsche1}, \cite{iarrobino1}, \cite{iarrobino3}), with $\hilb_d^3$ being singular for $d \geq 2$.\ Not all zero-dimensional subschemes of length three supported at one point are important in what follows:\ Only \emph{planar} subschemes (i.e.\ those with two-dimensional tangent spaces) need to be considered, whereas the so-called \textit{curvilinear} subschemes are neglected.  \par
\noindent A reason for concentrating on the planar locus is that the curvilinear locus is only \emph{locally} closed inside $X^{[3]}$ and has the structure of an affine fibre bundle over a projective bundle over $X$.\ In contrast, the planar locus has the structure of a locally trivial Grassmannian bundle $\mathbb{G}$ over $X$, being easier to handle especially from the derived category point of view.\ More work by the author on derived categories of Hilbert schemes of points is in progress, involving also the curvilinear locus and addressing missing components in (\ref{conjecture}).

\subsection*{Structure of the paper}

\noindent Let $\mathbb{G}$ be the locus of planar subschemes concentrated at single points of $X$.\ Once a closed embedding $\iota: \mathbb{G} \rightarrow X^{[3]}$ has been constructed \emph{functorially} following \cite{krugrennemo}, the Fourier--Mukai transform along the diagram\footnote{Diagrams of this form will also be called \textit{roofs}.}
\begin{center}
    \begin{tikzcd}
\mathbb{G}  \arrow[d, "p"'] \arrow[r, "\iota", hook] & {X^{[3]}} \\
X                                                    &          
\end{tikzcd}
\end{center}
can be considered.\ Tensoring with Schur functors $\Sigma^{\alpha}$ applied to tautological bundles $\mathcal{Q}^{\vee}$ on the relative Grassmannian $\mathbb{G}$ yields many more functors $\Phi_{\alpha}: \db(X) \rightarrow \db(X^{[3]})$ of Fourier--Mukai type.\ In order to examine $\Phi_{\alpha}$ for fully-faithfulness, the cohomology of the normal bundle $\mathcal{N}_{\mathbb{G}/X^{[3]}}$, tensored with different tautological bundles and restricted to the Grassmannian fibres of $p$, has to be computed.\ The first important result is the following.

\begin{thrm}[Thm.\ \ref{fullyfaithfultheorem1} below]\label{fullyfaithfulthm}
    Let $X$ be a smooth projective variety of dimension $d \geq 5$ over a field $k = \overline{k}$ of characteristic zero. Suppose $\alpha = (\alpha_1,\alpha_2)$ is a Young diagram inscribed into a rectangle with $2$ rows and $d-2$ columns.\ If $\lambda_{\alpha} = \alpha_1 - \alpha_2 \leq d-5$, then the Fourier--Mukai transform  $\Phi_{\alpha}: \db(X) \rightarrow \db(X^{[3]})$ with kernel $\Sigma^{\alpha}\mathcal{Q}^{\vee}$ is fully faithful.
\end{thrm}

\noindent This theorem yields $\binom{d-3}{2} + 3(d-4)$ fully faithful functors and consequently admissible subcategories of $\db(X^{[3]})$.\ For two different Young diagrams $\alpha \prec \beta$, it turns out that the necessary condition for semi-orthogonality of $\Phi_{\alpha}$ and $\Phi_{\beta}$ is $\alpha_1 - \beta_2 \leq d-5$.\ Excluding some Young diagrams in Theorem \ref{fullyfaithfulthm} yields the wanted semi-orthogonal sequence:

\begin{thrm}[Thm.\ \ref{semiorthogonalsequencethrm} below] \label{shortversioncollection}
    Restricting the attention to the fully faithful functors $\Phi_{\alpha}: \db(X) \rightarrow \db(X^{[3]})$ in Theorem \ref{fullyfaithfulthm} such that $\alpha_2 \geq 3$ yields a collection of $\binom{d-3}{2}$ semi-orthogonal subcategories of $\db(X^{[3]})$, all equivalent to $\db(X)$.
\end{thrm}

\noindent The proofs of the above theorems use cohomological computations on Grassmannians and the Borel--Weil--Bott Theorem like in Kapranov's classical work \cite{kapranov1}, as well as the following ingredient which is also interesting in its own right from the point of view of \textit{deforming planar zero-dimensional subschemes} inside a smooth variety.

\begin{thrm}[Thm.\ \ref{maintheoremnormalbundle} \& Cor.\ \ref{moreexplicitsplitting} below] \label{maintheoremnormalbundleshort}
    Let $X$ be a smooth projective variety over an algebraically closed field $k$ of any characteristic.\ Let $d = \dime(X) \geq 2$, write $\mathbb{G} = \mathbb{G}(2,\Omega_X) \hookrightarrow X^{[3]}$ and let $\mathcal{N} = \mathcal{N}_{\mathbb{G}/X^{[3]}}$.\ Then there is a short exact sequence
    \begin{align} \label{descriptionN}
        0 \longrightarrow \mathcal{Q}^{\vee} \longrightarrow \mathcal{Q} \otimes S^2\mathcal{Q}^{\vee} \longrightarrow \mathcal{N} \longrightarrow 0,
    \end{align}
    where $\mathcal{Q}$ denotes the tautological quotient bundle of rank $2$ on $\mathbb{G}$.\ If $d \geq 3$ and $\chara k = 0$, the sequence (\ref{descriptionN}) splits and $\mathcal{N} \cong S^3\mathcal{Q}^{\vee} \otimes \det \mathcal{Q}$.
\end{thrm}

\noindent As a corollary to Theorem \ref{shortversioncollection}, the following is obtained in the end:

\begin{thrm}[Prop.\ \ref{longversionkummer} below] \label{shortversionkummer}
    Let $A$ be an abelian variety of dimension $d \geq 5$ over an algebraically closed field of characteristic zero.\ Let $K_n(A)$ be the $n$-th generalized Kummer variety of $A$.\ Then the fully faithful functors from Theorem \ref{shortversioncollection} yield an exceptional sequence of length $\binom{d-3}{2}\cdot 3^{2d}$ in $\db(K_3(A))$.
\end{thrm}

\subsection*{Notations and conventions} Unless mentioned otherwise, $X$ denotes a smooth and projective variety over an algebraically closed field $k$, with $k$ of characteristic zero from Section \ref{sectionphialpha} on.\footnote{The characteristic of $k$ is important e.g.\ for representation-theoretic reasons, while the assumption $k = \overline{k}$ allows to check many properties of schemes and morphisms between them on $k$-points.}\ Smoothness and projectivity make life easier on the level of derived categories, but are also important for other reasons like choices of local parameters. \par
\noindent Categories are written in bold letters and abbreviated in a self-explaining way, like e.g. $\mathcal{O}_X$-$\mathbf{Mod}$ and $\qcoh(X)$.\ Derived categories $\db(X)$ are bounded derived categories of coherent sheaves on $X$, cf.\ Section \ref{subsectionderivedcats}.\ Complexes of sheaves are written as $\mathcal{F}^{\bullet}$ to distinguish them from ordinary sheaves $\mathcal{F} \in \coh(X)$. Stalks are written as $\mathcal{F}_{P}$, whereas fibres are written as $\mathcal{F}(P)$ or $\iota_P^{\ast}\mathcal{F}$.\ Here $\iota_P$ is the inclusion of $P \in X$.\ The skyscraper sheaf at a closed point $P \in X$ is denoted by $k(P) \in \coh(X)$.\par
\noindent For $S$-schemes $X$ and $T$, write $X(T) = \mor_S(T,X) = h_X(T)$.\ Grassmannians parametrize quotients, not subspaces, following Grothendieck's convention.\ The $k$-points of moduli spaces are usually denoted in square brackets $[-]$.\ Vector subspaces are denoted as $U \leq V$, while ideals of a commutative ring are denoted as $I \trianglelefteq R$.

\subsection*{Acknowledgements}

This work will be part of the author's dissertation under supervision of Andreas Krug, whom I would like to thank for his constant support, for helpful comments and many inspirations.\par
\noindent Furthermore, I am thankful for enlightening discussions with (former) members of the Institute of Algebraic Geometry at Leibniz University Hannover, especially Nebojsa Pavic, and fruitful conversations with Leonie Kayser from MPI MiS Leipzig.

\section{Preliminaries} \label{preliminaries}

In this section, the assumptions on $X$ and the base field (or scheme) are relaxed a little.

\subsection{Semi-orthogonal decompositions in Derived Categories} \label{subsectionderivedcats}

\noindent Let $X$ be a noetherian scheme over a field $k$ and let $\db(X) = \db(\coh(X))$ denote the bounded derived category of complexes of coherent sheaves on $X$.\par 
\noindent Many basic properties of $\db(X)$ can be found in \cite{huybrechts}.\ General derived categories and derived functors are definded e.g.\ in \cite[Ch.\ 3]{gmanin}.\ Geometric derived functors are discussed in \cite[Ch.\ 21]{goertzwedhorn2} or in \cite[Ch.\ 3.3]{huybrechts}.\ They are defined for quasi-coherent complexes first but restrict to bounded complexes of coherent sheaves in many cases.\par 
\noindent Since unbounded, non-coherent or non-perfect complexes will never occur in the main text, the derived functors primarily considered are derived pushforward $Rf_{\ast}:\db(X) \rightarrow \db(Y)$ for \textit{proper} morphisms $f:X \rightarrow Y$, derived pullback $Lf^{\ast}: \db(Y) \rightarrow \db(X)$ for any morphism with the convention $L^if^{\ast}(-) = \mathcal{H}^{-i}(Lf^{\ast}(-))$, derived tensor product $\otimes^L$ as well as derived sheaf homomorphisms  $R\mathcal{H}om_{\mathcal{O}_Y}(-,-)$, where $\mathcal{E}xt^i(-,-) := \mathcal{H}^i(R\mathcal{H}om(-,-))$.\footnote{For the latter three derived functors, it is assumed that the ambient scheme $Y$ has the resolution property and that it is regular.\ The resolution property is implied for example by quasi-projectivity over $k$ \cite[Prop.\ 22.57]{goertzwedhorn2}.\ These assumptions make it possible to compute derived functors using finite locally free resolutions instead of the more general flat resolutions.} The geometric derived functors can be concatenated to obtain an important class of functors, namely the Fourier--Mukai transforms $\Phi_{\mathcal{K}^{\bullet}}$ (FM transforms for short):

\begin{defi}
    If $p_X$ and $p_Y$ are the projections from $X \times_k Y$ to $X$ and $Y$ and if $\mathcal{K}^{\bullet}$ is an object in $\db(X \times_k Y)$, then $\Phi_{\mathcal{K}^{\bullet}} := Rp_{Y\ast}(\mathcal{K}^{\bullet} \otimes^L p_X^{\ast}(-))$ with kernel object $\mathcal{K}^{\bullet}$.\
\end{defi}

\noindent Any \textit{roof} of $S$-schemes $X \overset{p}{\longleftarrow} V \overset{q}{\longrightarrow} Y$ (regular and (quasi-)projective over $k$) for a $k$-scheme $S$ with $\mathcal{K}^{\bullet} \in \db(V)$ yields an FM transform $\Phi_{Ri_{\ast}\mathcal{K}^{\bullet}} = Rq_{\ast}(\mathcal{K}^{\bullet} \otimes^L Lp^{\ast}(-))$ via an induced morphism $i: V \rightarrow X \times_S Y \hookrightarrow X \times_k Y$ as long as $i$ is proper.\ Properness of $i$ is unimportant if unbounded complexes are allowed.\ These functors may be called \textit{relative} FM transforms \cite[Ch.\ 6]{fouriernahm} or \emph{kernel functors} \cite[2.5]{homprojdual}.\par 
\noindent The criterion of Bondal and Orlov (\cite{bondalorlov}, taken here from \cite[Prop.\ 7.1]{huybrechts}) allows to examine Fourier--Mukai transforms for fully-faithfulness by looking at closed points:

\begin{thrm}[{\cite[Thm.\ 1.1]{bondalorlov}}] \label{bondalorlovff}
    Let $\Phi: \db(X) \rightarrow \db(Y)$ be an FM transform between smooth projective varieties.\ Then $\Phi$ is fully faithful if and only if for all closed $x,y \in X$, 
    \begin{align*}
        \homo_{\db(Y)}\bigl(\Phi(k(x)),\Phi(k(y))[i]\bigr) = \begin{cases}
            k \cdot \id,~ x=y \text{ and } i = 0, \\ 0,~ x\neq y \text{ or } i \notin [0,\dime X].
        \end{cases}
    \end{align*}
\end{thrm}

\begin{bem} \label{remarkspanningclasses}
    Theorem \ref{bondalorlovff} uses that the involved skyscraper sheaves $\lbrace k(x) \rbrace_{x \in X}$ form a \emph{spanning class}.\ Spanning classes are in general useful for examining FM transforms, another example being given by powers of ample line bundles \cite[Sect.\ 3.2]{huybrechts}.
\end{bem}

\noindent It is also possible (and much easier) to state and prove the following orthogonality result:

\begin{lemma} \label{orthogonalitycriterion}
Let $\Phi, \Psi: \db(X) \rightarrow \db(Y)$ be FM transforms between derived categories of smooth projective varieties.\ Let $\Omega \subseteq \db(X)$ be a spanning class and suppose that
\begin{align*}
\homo_{\db(Y)}\bigl(\Phi(L'),\Psi(L)[i]\bigr) = 0 \quad \text{ for all } L',L \in \Omega \text{ and } i \in \mathbb{Z}.
\end{align*}
Then $\homo_{\db(Y)}(\Phi(\mathcal{F}^{\bullet}),\Psi(\mathcal{G}^{\bullet})) = 0$ for all $\mathcal{F}^{\bullet},\mathcal{G}^{\bullet} \in \db(X)$, i.e.\ the essential images of $\Phi$ and $\Psi$ are semi-orthogonal to each other. 
\end{lemma}

\begin{proof}
    Note that $\homo_{\db(Y)}(\Phi(-),\Psi(-)) = 0$ if and only if $\Theta := \Phi^R \circ \Psi = 0$, with $\Phi^R$ being the right adjoint of $\Phi$.\ Note also that $\Theta$ is of FM type again and has both adjoints \cite[Prop.\ 5.9 \& 5.10]{huybrechts}.\ The assumption yields that $\homo_{\db(X)}(L',\Theta(L)[i]) = 0$ for all $L',L \in \Omega$ and $i \in \mathbb{Z}$. By definition of a spanning class $\Omega$, it follows that $\Theta(L) = 0$ for all $L \in \Omega$. If $\mathcal{G}^{\bullet} \in \db(X)$ is arbitrary, then
\begin{align*}
\homo_{\db(X)}(L,\Theta^R(\mathcal{G}^{\bullet})[i]) \cong \homo_{\db(X)}(\Theta(L),\mathcal{G}^{\bullet}[i]) = 0 \text{ for all } L \in \Omega \text{ and } i \in \mathbb{Z}.
\end{align*} 
Therefore $\Theta^R = 0$ and consequently $\Theta = (\Theta^R)^L = 0$.
\end{proof}

\noindent Derived categories are triangulated categories \cite[Ch.\ 4]{gmanin}. Note that only the existence of adjoints was used to prove Lemma \ref{orthogonalitycriterion} so that its proof can be adapted to other triangulated categories containing spanning classes.\ The following notion of a semi-orthogonal sequence will be applied to $\category = \db(X)$: \par 

\begin{defi}
    Let $\category$ be a triangulated category with full triangulated subcategories $\mathbf{D}_1,\hdots,\mathbf{D}_n$. The sequence  $(\mathbf{D}_1,\hdots,\mathbf{D}_n)$ is called \emph{semi-orthogonal} provided that $\mathbf{D}_i \subseteq \mathbf{D}_j^{\perp}$ for $i < j$, i.e.\ there are no non-zero morphisms from objects of $\mathbf{D}_j$ to objects of $\mathbf{D}_i$ in $\category$. The semi-orthogonal sequence $(\mathbf{D}_1,\hdots,\mathbf{D}_n)$ is called a semi-orthogonal \emph{decomposition} (or alternatively full) provided that $\langle \mathbf{D}_1,\hdots,\mathbf{D}_n \rangle = \mathbf{C}$.\footnote{Since this work establishes only some semi-orthogonal subcategories of $\db(X^{[3]})$ that are at best \textit{part} of a semi-orthogonal decomposition, fullness is of minor importance in what follows.}
\end{defi}

\noindent See also the exposition in \cite[Sect.\ 1.1]{kuznetsovexp}.\ Often, it is additionally required that the subcategories be admissible \cite{bondalkapranov}, i.e.\ that the inclusions $\mathbf{D}_i \hookrightarrow \mathbf{C}$ admit left and right adjoints.\ This simplifies talking about orthogonal complements of subcategories (see e.g.\ Lemma \ref{complement}).\ A subcategory isomorphic to the bounded derived category of a smooth projective variety is always admissible \cite[Prop.\ 2.6 \& Thm.\ 2.14]{bondalkapranov}.\par 

\begin{lemma}[{\cite[Lemma 2.4]{homprojdual}}, after \cite{bondalrep}] \label{complement}
    Let $(\mathbf{D}_1,\hdots,\mathbf{D}_n)$ be a semi-orthogonal sequence of admissible triangulated subcategories in a triangulated category $\mathbf{C}$.\ Then there exists a semi-orthogonal decomposition $(\mathbf{D}_1,\hdots,\mathbf{D}_n,{}^{\perp}\langle \mathbf{D}_1,\hdots,\mathbf{D}_n \rangle)$ of $\mathbf{C}$.
\end{lemma}

\begin{prop} \label{appliedfibrecriterion}
Let $k = \overline{k}$.\footnote{This is only used to ensure that any closed point is $k$-valued.}\ Consider a roof $X \leftarrow V \rightarrow Y$ between smooth projective varieties.\ Here $p: V \rightarrow X$ is a smooth proper morphism and $\iota:V \rightarrow Y$ is a closed embedding.\ Denote the fibres of $p$ over $k$-points $P \in X$ by $F_P = p^{-1}(P)$.\ Let $\mathcal{E}$ and $\mathcal{G}$ in \emph{$\mathbf{\coh}$}$(V)$ be locally free sheaves that are kernels of the two Fourier--Mukai transforms $\Phi = \Phi_{\mathcal{E}}$ and $\Psi = \Psi_{\mathcal{G}}$ along $V$.\ Their restrictions are denoted by $E = \mathcal{E}_{|F_{P}}$ and $G = \mathcal{G}_{|F_{P}}$.
\begin{enumerate}
    \item[(1)] The Fourier--Mukai transform $\Phi$ is fully faithful provided that $\homo(E,E) \cong k$ on every fibre $F = F_P$ and that $H^p(F,\wedge^q(\mathcal{N}_{V/Y})_{|F} \otimes E \otimes E^{\vee}) = 0$ for $p+q > 0$.
    \item[(2)] The essential images of $\Phi$ and $\Psi$ are semi-orthogonal (~$\image \Psi \subseteq \image \Phi^{\perp}$) if on every fibre $F = F_P$, $\homo(E,G) = 0$ and $H^p(F, \wedge^q(\mathcal{N}_{V/Y})_{|F} \otimes G \otimes E^{\vee}) = 0$ for $p + q > 0$.
\end{enumerate}
\end{prop}

\begin{proof} The proof idea for (1) is essentially the same as for \cite[Prop.\ 3]{belmansflips}, namely to apply Theorem \ref{bondalorlovff} to the functor $\Phi = \Phi_{\mathcal{E}}$.\ The idea also occurs in the proof of Orlov's derived blow-up formula in \cite[Prop.\ 11.16]{huybrechts}.\ For assertion (2), Lemma \ref{orthogonalitycriterion} will be used.\ Denote the involved morphisms in the definition of $\Phi$ and $\Psi$ as follows.
\begin{center}
\begin{tikzcd}
F \arrow[r, "\iota'"', hook] \arrow[d] \arrow[rd, "\square", phantom] \arrow[rr, "j" description, bend left] & V \arrow[r, "\iota"', hook] \arrow[d, "p"] & Y \\
\lbrace P \rbrace \arrow[r, hook]                                                                            & X                                     &  
\end{tikzcd}
\end{center}
\noindent Let $\mathcal{N}' = (\mathcal{N}_{V/Y})_{|F}$ and $\mathcal{N} = \mathcal{N}_{F/Y}$ be the normal bundles. Identically to \cite[(2.9)]{belmansflips}, the relative short exact sequence of normal bundles for  $F \subseteq V \subseteq Y$ reads
\begin{align} \label{relativenormalbundles}
0 \longrightarrow \mathcal{O}_F^{\oplus d} \longrightarrow \mathcal{N} \longrightarrow \mathcal{N}' \longrightarrow 0.
\end{align}
\noindent In both cases (1) and (2), one has to examine the $\homo$-spaces
\begin{align*}
\homo_{\db(Y)}\bigl(\Phi(k(P)), \Psi(k(P))[i]\bigr) \cong \ext^i_{\mathcal{O}_Y}(j_{\ast}(\mathcal{E}_{|F}),j_{\ast}(\mathcal{G}_{|F})) \cong H^i(F,(j^!j_{\ast}G) \otimes E^{\vee}),
\end{align*}
with $\Phi = \Psi$ and $E = G$ in the setting of (1).\ There is always the following converging hypercohomology spectral sequence (see Lemma \ref{formulapushpull} below):
\begin{align} \label{hypercohomologyspectralseq}
E_2^{p,q} = H^p(F, G \otimes E^{\vee} \otimes \wedge^q\mathcal{N}) \Longrightarrow H^{p+q}(F,j^!j_{\ast}G \otimes E^{\vee}).
\end{align} 
\noindent Furthermore, there exists a filtration \cite[Exc. II.5.16]{hartshorne} induced from (\ref{relativenormalbundles}):
\begin{gather*}
0 \longrightarrow F^1 \otimes G \otimes E^{\vee} \longrightarrow \wedge^q \mathcal{N} \otimes G \otimes E^{\vee} \longrightarrow \wedge^q \mathcal{N}' \otimes G \otimes E^{\vee} \longrightarrow 0, \\
0 \longrightarrow F^2 \otimes G \otimes E^{\vee} \longrightarrow F^1 \otimes G \otimes E^{\vee} \longrightarrow \mathcal{O}_F^{\oplus d} \otimes \wedge^{q-1}\mathcal{N}' \otimes G \otimes E^{\vee} \longrightarrow 0, \\
\vdots \\
0 \longrightarrow F^q \otimes G \otimes E^{\vee} \longrightarrow F^{q-1} \otimes G \otimes E^{\vee} \longrightarrow \wedge^{q-1}\mathcal{O}_F^{\oplus d} \otimes \mathcal{N}' \otimes G \otimes E^{\vee} \longrightarrow 0, \\
0 \longrightarrow 0 \longrightarrow F^q \otimes G \otimes E^{\vee} \overset{\sim}{\longrightarrow} \wedge^q \mathcal{O}_F^{\oplus d} \otimes G \otimes E^{\vee} \longrightarrow 0.
\end{gather*}
\noindent What has to be shown for (1) is that (specializing to $\Phi  = \Psi$)\begin{align*}
\homo_{\db(Y)}(\Phi_{\mathcal{E}}(k(P)), \Phi_{\mathcal{E}}(k(P))[i]) = \begin{cases} k \cdot \id, \quad i = 0, \\ 0, \quad i \notin [0,d],~ d = \dime(X). \end{cases}
\end{align*}
\noindent This leads to proving $E_2^{p,q} = 0$ for $p+q > d$ and $E_2^{0,0} \cong k$ in (\ref{hypercohomologyspectralseq}) with $E=G$.\ Notice that by assumption, $E_2^{0,0} = \homo(E,E) \cong k$ holds.\ Let $p + q > d$ and consider $E_2^{p,q}$ with $p = 0$ and $q > d$ first.\ Then in the above filtration, $F^q \otimes E \otimes E^{\vee} = F^q \otimes G \otimes E^{\vee} = 0$ and 
\begin{align*}
H^0(F,\wedge^i\mathcal{O}_F^{\oplus d} \otimes \wedge^{q-i}\mathcal{N}' \otimes E \otimes E^{\vee}) = H^0(F,\wedge^{q-i}\mathcal{N}' \otimes E \otimes E^{\vee})^{\oplus \binom{d}{i}} = 0 
\end{align*}
by assumption for $0 \leq i < q$.\ Here the convention $\binom{d}{i} = 0$ is used if $i > d$. Inductively, 
\begin{align*}
H^0(F,\wedge^q\mathcal{N} \otimes E \otimes E^{\vee}) \cong H^0(F,F^1 \otimes E \otimes E^{\vee}) \cong \hdots \cong H^0(F,F^q \otimes E \otimes E^{\vee}) = 0.
\end{align*}
If $p > 0$, then first of all $H^p(F,F^q \otimes E \otimes E^{\vee}) = H^p(F,\wedge^q \mathcal{O}_F^{\oplus d} \otimes E \otimes E^{\vee}) = 0$, which is only a nontrivial statement if $0 \leq q \leq d$. Being zero uses the assumption that $\ext^{\ast}(E,E)$ is concentrated in degree zero. Next, $
H^p(F,\wedge^i\mathcal{O}_F^{\oplus d} \otimes \wedge^{q-i}\mathcal{N}' \otimes E \otimes E^{\vee}) = 0$ for $0 \leq i < q$ since then $p + (q-i) > q-i > 0$. The same inductive argument as for $p=0$ applied to the filtration of $\wedge^q\mathcal{N}_{F/Y} \otimes E \otimes E^{\vee}$ finishes the proof of the proposition. \par 
\noindent Statement (2) can be proven as an application of Lemma \ref{orthogonalitycriterion}:\ For any $k$-point $P \in X$, it has to be shown that $H^i(F,(j^!j_{\ast}G) \otimes E^{\vee}) = 0$ for all $i \in \mathbb{Z}$, which can be achieved by showing that $E_2^{p,q} = H^p(F,\wedge^q\mathcal{N} \otimes G \otimes E^{\vee}) = 0$ for all $p,q \geq 0$, cf.\ (\ref{hypercohomologyspectralseq}).\ This is done via the filtration from above and a similar case distinction into $p = 0$ or $p>0$.\par 
\noindent Alternatively (instead of skyscraper sheaves, see Remark \ref{remarkspanningclasses}), the spanning class of line bundles on $X$ may be used to show (2) by proving that $(\Phi^R \circ \Psi)(\mathcal{L}) = 0$ for $\mathcal{L} \in \picard(X)$.\ The right adjoint $\Phi^R$ is $Rp_{\ast}(\mathcal{E}^{\vee} \otimes \iota^!(-))$ and similarly to (\ref{hypercohomologyspectralseq}), there is a spectral sequence 
\begin{align*}
    E_2^{i,j} = R^ip_{\ast}(\wedge^j \mathcal{N}_{V/Y} \otimes p^{\ast}\mathcal{L} \otimes \mathcal{G} \otimes \mathcal{E}^{\vee}) \Longrightarrow \mathcal{H}^{i+j}\left((\Phi^R \circ \Psi)(\mathcal{L})\right).
\end{align*}
By the projection formula, $E_2^{i,j} \cong R^ip_{\ast}(\wedge^j \mathcal{N}_{V/Y} \otimes \mathcal{G} \otimes \mathcal{E}^{\vee}) \otimes \mathcal{L}$, which is zero by assumption using the theorem on cohomology and base change \cite[Thm.\ III.12.11]{hartshorne}.
\end{proof}

\subsection{Hilbert Schemes of Points} \label{preliminarieshilbert} Let $X$ be projective over a noetherian scheme $S$ and let $\mathcal{F}\in \mathbf{Coh}(X)$.\ For every desired Hilbert polynomial $P \in \mathbb{Q}[x]$, Grothendieck showed in \cite[Thm.\ 3.1]{grothendieckhilbert} that the Quot scheme $\quot^P_{\mathcal{F}/X/S}$ parametrizing quotients of $\mathcal{F}$ exists.\ In this work, the main focus lies on the case where $\mathcal{F} = \mathcal{O}_X$, $S = \spec k$ and $P \equiv n \in \mathbb{Z}_{\geq 1}$, so that $\quot^P_{\mathcal{F}/X/S} = \hilb^n_{X/k} = X^{[n]}$ is the Hilbert scheme of $n$ points and 
\begin{align*}
    X^{[n]}(T) = \lbrace Z \hookrightarrow X \times_k T \text{ ~ closed subscheme}, q = (\text{pr}_{2})_{|Z} :Z \rightarrow T \text{ ~flat,~} h_{Z_t}   \equiv n ~\forall t \in T  \rbrace
\end{align*}
for any $T \rightarrow \spec k$.\ Note that the definition implies that all fibres $Z_t = q^{-1}(t)$ are zero-dimensional of length $n$ and  consequently, $\rank q_{\ast}\mathcal{O}_Z = n$.\ Since $X^{[n]}$ represents a functor, it comes along with a universal element (usually referred to as the \emph{universal family})
\begin{align*}
    \left(\pi: \Xi_n \subseteq X \times X^{[n]} \rightarrow X^{[n]}\right), \quad \rank \pi_{\ast}\mathcal{O}_{\Xi_n} = n.
\end{align*}
The family $\Xi_n$ represents a functor itself \cite[Rem.\ 1.1.4]{goettsche1} so that 
\begin{align*}
    X^{[n]}(k) &= \lbrace Z \subseteq X \text{~closed zero-dimensional subscheme of length~} n\rbrace, \\
    \Xi_n(k) &= \lbrace (x,[Z]) : [Z] \in X^{[n]}(k), ~ x \in Z \rbrace.
\end{align*}

\begin{prop} \label{closedsubscheme}
    Let $X \rightarrow \spec k$ be a projective scheme  and let $Y \subseteq X$ be a closed subscheme, whereas $U = X \setminus Y$ is the complementary open subscheme.\ Then $Y^{[n]}$ and $U^{[n]}$ exist as closed, respectively open subschemes of $X^{[n]}$. 
\end{prop}

\begin{proof}
    The existence of $Y^{[n]}$ is clear since $Y$ is projective, and there is a closed embedding of Hilbert functors $\mathcal{H}ilb^n_{Y/k} \hookrightarrow \mathcal{H}ilb^n_{X/k}$ \cite[Lemma 5.17]{fganitsure}.\ On the other hand, the Hilbert functor $\mathcal{H}ilb^n_{U/k}$ has to be defined first as the subfunctor of $\mathcal{H}ilb^n_{X/k}$ parametrizing flat families over $T$ contained in $U \times_k T \subset X \times_k T$ \cite[Rem.\ 4.a]{grothendieckhilbert}.\ Then one may take $U^{[n]}$ to be $X^{[n]} \setminus \pi\left(\Xi_n \cap (Y \times X^{[n]})\right)$.\ In fact, this is how existence of Hilbert schemes of quasi-projective schemes $U$ can be proven in general.
\end{proof}

\noindent A special feature of the Hilbert scheme of points is its relation to the symmetric product $X^{(n)} = X^n/S_n$, the latter being a quasi-projective variety whenever $X$ is.\ Namely if the base field $k$ is algebraically closed of characteristic zero and $X$ is smooth, there exists a surjective morphism $\rho: X^{[n]}_{\red} \rightarrow X^{(n)}$ \cite[Sect.\ 7.1]{fgahilbert}, the Hilbert--Chow morphism, given on $k$-points as 
    \begin{align*}
        [W] \mapsto \sum_{P \in \supp W} \length(\mathcal{O}_{W,P}) \cdot [P], \quad [W] \in X^{[n]}(k).
    \end{align*}
Recall from the introduction that $X^{[n]}$ is smooth and connected if $\dime X \in \lbrace 1,2 \rbrace$ or if $n \in \lbrace 1,2,3\rbrace$ and $X$ is smooth and connected.\ Hence $X^{[n]} = X^{[n]}_{\red}$ in those cases.\ As in \cite{goettsche1}, $\rho$ induces a useful stratification of $X^{[n]}_{\red}$ by locally closed subschemes $X^{[n]}_{\nu}$, where 
\begin{align*}
    X^{[n]}_{\nu} = \rho^{-1}(X^{(n)}_{\nu}),\quad  X^{(n)}_{\nu} = \left\lbrace \sum \nu_i[P_i] : \text{ all $P_i \in X$ distinct } \right\rbrace, \quad \nu \dashv n ~\text{ a partition.}
\end{align*}
  In the rest of this section, the case $\nu = (n)$ is discussed.\ The closed stratum $X^{[n]}_{(n)} \subseteq X^{[n]}$ has zero-dimensional subschemes of $X$ concentrated at single points as $k$-points.\par 
\noindent Fixing a point $P \in X$, the question is how to describe all possible ideals that yield zero-dimensional subschemes of length $n$ at $P$.\ Such ideals of colength $n$ are parametrized by the \textit{Local} Punctual Hilbert scheme, of which $X^{[n]}_{(n),\red}$ is a globalization since it can also be defined as the relative Hilbert scheme\footnote{Relative Hilbert schemes are defined over more general base schemes $S$.\ Base change \cite[Prop.\ 4.4.3]{sernesi} implies that a relative $\hilb_{X/S}^n$ is a family with absolute Hilbert schemes $\hilb_{X_s/k(s)}^n$ as fibres.} of the $n$ times thickened diagonal $n\Delta \subseteq X \times_k X$ over $X$ after equipping both with their reduced scheme structure \cite[Lemma 2.1.2]{goettsche1}:

\begin{defi}
    The \emph{Local Punctual Hilbert scheme} is $\hilb^n_{d,X,P} = (\spec \mathcal{O}_{X,P}/\mathfrak{m}^n)^{[n]}_{\red}$ for a closed point $P \in X$, where $\mathfrak{m} = \mathfrak{m}_{X,P}$ and $d = \dime X$.
\end{defi}

\noindent The definition depends only on the dimension $d\geq 1$, up to \textit{non-canonical} isomorphism:\ If $Z \subseteq X$ is a closed subscheme of length $n$ concentrated at a single point $P$, then $I := \mathcal{I}_{X,P}$ has colength $n$, i.e. $\dime_k \mathcal{O}_{Z,P} = \dime_k \mathcal{O}_{X,P}/I = n$. Necessarily, $I$ contains $\mathfrak{m}^n = \mathfrak{m}_{X,P}^n$, cf.\ the proof of \cite[Lemma 1.3.2]{goettsche1}.\ The specific choice of a point $P \in X$ is not important in the definition of $\hilb^n_d = \hilb^n_{d,X,P}$ since $X$ is smooth and hence
\begin{align} \label{formalcomputation}
    \mathcal{O}_{X,P}/\mathfrak{m}^n \cong \widehat{\mathcal{O}_{X,P}}/\widehat{\mathfrak{m}}^n \cong k[[x_1,\hdots,x_d]]/\widehat{\mathfrak{m}}^n \cong k[x_1,\hdots,x_d]/\mathfrak{m}^n.
\end{align}
This justifies the (sloppy) notation $\hilb^n_{d}$, omitting the choice of local parameters.\ By Proposition \ref{closedsubscheme}, $\hilb^n_d$ is a closed subscheme of $X^{[n]}$.\ It has quite different properties compared to the \enquote{global} $X^{[n]}$ and can be realized as a closed subscheme of the Grassmannian $G(n,k[x_1,\hdots,x_d]/\mathfrak{m}^n)$.\ The case $d = 2$ has been extensively studied in \cite{iarrobino2}, cf.\ \cite[p. 10]{goettsche1}.\ There it is also explained how distinguishing colength $n$ ideals $I$ by their Hilbert function yields a stratification of $\hilb^n_d$.\ The Hilbert function is invariant under the identifications made in (\ref{formalcomputation}).

\begin{bsp} \label{examplehilb}
    If $n=1$, then $X^{[1]} \cong X$. If $n=2$, zero-dimensional subschemes $Z \subseteq X$ of length two can be the union of two reduced points or they are concentrated at a fixed point $P$.\ In the latter case, $Z$ corresponds to the choice of a tangent direction at $P \in X$, i.e.\ to a one-dimensional quotient of $\mathfrak{m}/\mathfrak{m}^2$. \par
    \noindent If $n=3$, subschemes of length $3$ are either combinations of the previous cases or they have support at a single point $P \in X$.\ There are two possible types of ideals in $\hilb^3_d$:\par
    \noindent The \emph{curvilinear} ideals define subschemes contained in the germ of a smooth curve.\ They contain second order data and are of minor importance in what follows.\ The \emph{planar} ideals are easier to handle and are locally formally of the form 
    \begin{align*}
        (u_3,\hdots,u_d) + \mathfrak{m}^2, \quad d\geq 2, \quad \degr u_i = 1,
    \end{align*}
    corresponding to the choice of a plane in the tangent space $T_PX$ of $X$ at $P$.
\end{bsp}

\noindent Example \ref{examplehilb} is an instance of the general  stratification of $\hilb_d^n$ by locally closed subschemes $Z_T$ corresponding to ideals with Hilbert function $T$ \cite{iarrobino1}.\ This can be globalized to a stratification of $X^{[n]}_{(n),\red} \subseteq X^{[n]}$ by subschemes $Z_T(X)$ \cite[Sect.\ 2.1]{goettsche1}. 

\section{The locus of planar subschemes in \texorpdfstring{$X^{[n]}$}{X[n]}}

A fully faithful functor $\db(X) \rightarrow \db(X^{[n]})$ is automatically a Fourier--Mukai transform due to results of Bondal, Orlov and van den Bergh, as explained in \cite[Thm.\ 5.14]{huybrechts}.\ From a geometric point of view, such a construction might rather be related to loci of zero-dimensional subschemes concentrated in \emph{single} points, whereas subschemes concentrated in two or more points lead to the study of products of $X$ with itself and $X^{[m]}$, $m < n$.\ As explained in the introduction, curvilinear subschemes are less important in what follows and the focus lies on the planar zero-dimensional subschemes of $X$.\ This locus has the structure of a Grassmannian bundle $\mathbb G \rightarrow X$ with Grassmannian fibre $G(l,\mathfrak{m}/\mathfrak{m}^2)$.\footnote{More precisely, the locus of planar subschemes is $Z_{(1,l,0)}(X) \cong \mathbb{G}(l,\Omega_X) = \text{Grass}(l,T_X^{\ast})$ in the notation of \cite{goettsche1}.\ The defining ideal of a planar subscheme is determined by an $l$-dimensional subspace of $T_PX$, see Section \ref{preliminarieshilbert}.\ Notice also that $Z_{(1,2,0)}(X) = B_3$ in \cite[Sect.\ 5.1]{shenvial}.}

\subsection{Functorial description of 
\texorpdfstring{$\mathbb{G}(l,\Omega_X) \hookrightarrow X^{[l+1]}$}{G -> X[l+1]}
} 

\noindent Let $X$ be a smooth projective variety over the field $k$ of dimension $d$ and let
\begin{align} \label{tautsequence}
0 \longrightarrow \mathcal{K} \longrightarrow p^{\ast}\Omega_{X/k} \longrightarrow \mathcal{Q} \longrightarrow 0
\end{align}
be the tautological exact sequence on the Grassmannian bundle $\mathbb{G} := \mathbb{G}(l,\Omega_X) \overset{p}{\longrightarrow} X$, where $\mathcal{Q}$ has rank $l \leq d$. As in \cite[Sect.\ 4.2]{krugrennemo}, a morphism $\iota: \mathbb{G} \rightarrow X^{[l+1]}$ is constructed by exhibiting a flat family $Z$ of degree $l+1$ over $\mathbb{G}$.\ Concretely, the family $Z \subseteq \mathbb G \times X$ is defined by a commutative diagram of exact sequences of sheaves on $\mathbb G \times X$ of the form 

\adjustbox{scale=0.94,center}{
\begin{tikzcd}
            & 0 \arrow[d]                                              & 0 \arrow[d]                                                    & 0 \arrow[d]                                                        &   \\
0 \arrow[r] & {(\id,p)_{\ast}\mathcal{K}} \arrow[d] \arrow[r, "="]     & {(\id,p)_{\ast}\mathcal{K}} \arrow[r] \arrow[d, "\tau"]      & 0 \arrow[r] \arrow[d]                                              & 0 \\
0 \arrow[r] & {(\id,p)_{\ast}p^{\ast}\Omega_{X/k}} \arrow[d] \arrow[r] & (p \times \id)^{\ast}\mathcal{O}_{2\Delta} \arrow[d] \arrow[r] & (p \times \id)^{\ast}\mathcal{O}_{\Delta} \arrow[r] \arrow[d, "="] & 0 \\
0 \arrow[r] & {(\id,p)_{\ast}\mathcal{Q}} \arrow[d] \arrow[r]          & \coker \tau \arrow[d] \arrow[r]                              & (p \times \id)^{\ast}\mathcal{O}_{\Delta} \arrow[r] \arrow[d]     & 0. \\
            & 0                                                        & 0                                                              & 0                                                              &  
\end{tikzcd} }

\noindent The first column is the pushforward of (\ref{tautsequence}) under $(\id,p)$, while the second row uses base change along the cartesian square
\begin{center}
    \begin{tikzcd}
\mathbb G \arrow[d, "p"'] \arrow[r, "{(\id,p)}", hook] \arrow[rd, "\square", phantom] & \mathbb G \times X \arrow[d, "p \times \id"] \\
X \arrow[r, "\Delta", hook]                                                           & X \times X                                  
\end{tikzcd}
\end{center}
and originates from the short exact sequence on $X \times X$
\begin{align} \label{differentialsequence}
    0 \longrightarrow \Delta_{\ast}\Omega_{X/k} \longrightarrow \mathcal{O}_{2\Delta} \longrightarrow \mathcal{O}_{\Delta} \longrightarrow 0.
\end{align}
Here the $n$-fold diagonal $n\Delta \subseteq X \times X$ is defined by the $n$-th power of the ideal sheaf $\mathcal{I}_{\Delta}$.\ Observe that $(p \times \id)^{\ast}\mathcal{O}_{\Delta} \cong (\id,p)_{\ast}\mathcal{O}_{\mathbb{G}} = \mathcal{O}_{\Gamma_p}$ and that the third row of the diagram is induced by the ones above.\ Now $\coker \tau$ is a quotient of $\mathcal{O}_{X \times \mathbb{G}}$, hence can be written as 
\begin{align*}
    \coker \tau = \mathcal{O}_Z
\end{align*}
for a closed subscheme $Z$ supported on $\Gamma_p$.\ The pushdown to $\mathbb{G}$ is locally free of rank $l+1$ \cite[Sect.\ 4.2]{krugrennemo}, wherefore $Z \rightarrow \mathbb{G}$ is a flat family defining a morphism $\iota: \mathbb{G} \rightarrow X^{[l+1]}$.\ By the universal property of the Hilbert scheme, there exists a cartesian diagram  

\begin{center}
\begin{tikzcd}
Z \arrow[d] \arrow[r, "q"] \arrow[rd, "\square", phantom] & \mathbb{G} \arrow[d, "\iota"] \\
\Xi_{l+1} \arrow[r, "\pi"]                                                          & {X^{[l+1]}}.                 
\end{tikzcd}
\end{center}

\begin{bem} \label{remarkthickdiagonal}
Since $X$ is smooth, $\mathcal{O}_{2\Delta}$ is flat over $X$ by (\ref{differentialsequence}).\ Inductively, the same is true for all $\mathcal{O}_{n\Delta}$ because $\mathcal{I}_{\Delta}^n/\mathcal{I}_{\Delta}^{n+1} \cong S^n(\mathcal{I}_{\Delta}/\mathcal{I}_{\Delta}^2)$, cf. \cite[Thm.\ II.8.21A(e)]{hartshorne}.
\end{bem}

\noindent The next result is probably well-known but the author could not find a suitable reference.

\begin{prop} \label{basechangeflatideal}
    Suppose there is a cartesian diagram between quasi-projective varieties 
        \begin{center}
        \begin{tikzcd}
        X' \arrow[r, "f'"] \arrow[d, "g'"'] \arrow[rd, "\square", phantom] & X \arrow[d, "g"] \\
        Y' \arrow[r, "f"]                                                                & Y .              
        \end{tikzcd}
        \end{center}
Let $Z \subseteq X$ be a closed subscheme and $Z' = X' \times_X Z$ the scheme-theoretic preimage.\ If $Z$ is flat over $Y$, $\mathcal{I}_Z$ and $\mathcal{O}_Z$ pull back to $\mathcal{I}_{Z'}$ and $\mathcal{O}_{Z'}$ via $f'$.\ If additionally $g$ is flat, then on the level of derived categories, ${Lf'}^{\ast}\mathcal{I}_Z \cong \mathcal{I}_{Z'}$ and ${Lf'}^{\ast}\mathcal{O}_Z \cong \mathcal{O}_{Z'}$.
\end{prop}

\begin{proof}
The defining short exact sequence of $Z \subseteq X$ remains exact after pullback to $X'$ since $\mathcal{O}_Z$ is flat over $Y$ \cite[Prop.\ 7.40]{goertzwedhorn}.\ By affine base change \cite[\href{https://stacks.math.columbia.edu/tag/02KG}{02KG}]{stacks-project}, $\mathcal{O}_Z$ pulls back to $\mathcal{O}_{Z'}$, which proves the first claim.\ If $g$ is a flat morphism, then the 2-out-of-3 property \cite[Lemma 2.25]{kuznetsovhyperplanes} shows that ${Lf'}^{\ast}\mathcal{O}_Z \cong \mathcal{O}_{Z'}$.
\end{proof}

\begin{kor} \label{basechangeflatideal2}
    Let $f:X \rightarrow Y$ be a morphism between varieties with $Y$ smooth.\ Then for $n \geq 1$,  $n\Delta_Y \times_Y X$ is equal to the $n$ times thickened graph defined by $\mathcal{I}_{\Gamma_f}^n \subseteq \mathcal{O}_{X\times Y}$.
\end{kor}

\begin{proof}
    The statement for $n = 1$ is purely categorical.\ More generally, consider
\begin{center}
\begin{tikzcd}
n\Delta_Y \times_Y X \arrow[d] \arrow[r, hook] \arrow[rd, "\square", phantom] & X \times Y  \arrow[d, "f \times \id"] \\
n\Delta_Y \arrow[r, hook]                                                     & Y \times Y.               
\end{tikzcd}
\end{center}
By Proposition \ref{basechangeflatideal}, $(f \times \id)^{\ast}\mathcal{I}_{\Delta} \cong \mathcal{I}_{\Gamma_f}$ and $(f \times \id)^{\ast}(\mathcal{I}_{\Delta}^n) \cong \mathcal{I}_{n\Delta_Y \times_Y X}$ (Remark \ref{remarkthickdiagonal}).\ Since $\mathcal{I}_{\Delta}^n$ is the image (i.e. kernel of the cokernel) of the multiplication morphism $\mathcal{I}_{\Delta}^{\otimes n} \rightarrow \mathcal{I}_{\Delta}$, there is an epimorphism $(f \times \id)^{\ast}(\mathcal{I}_{\Delta}^n) \twoheadrightarrow \bigl((f \times \id)^{\ast}\mathcal{I}_{\Delta}\bigr)^n \cong \mathcal{I}_{\Gamma_f}^n$.\par
\noindent This has to be an isomorphism since $X \times Y$ is integral \cite[Prop.\ 5.51]{goertzwedhorn} and since (after going to the stalks) a non-trivial module homomorphism from a non-zero ideal of the base ring to a torsion-free module has kernel equal to zero. 
\end{proof}

\noindent Back to $\iota: \mathbb{G} \rightarrow X^{[l+1]}$, induced by the flat family $Z \cong \mathbb{G} \times_{X^{[l+1]}} \Xi^{l+1} \subseteq 2\Gamma_p$ over $\mathbb{G}$.

\begin{lemma}
In the above setting, there are short exact sequences on $\mathbb{G} \times X$
\begin{align} \label{idealequationg120}
0 \longrightarrow \mathcal{I}_{2\Gamma_p} \longrightarrow \mathcal{I}_Z \longrightarrow (\id,p)_{\ast}\mathcal{K} \longrightarrow 0,
\end{align}
\begin{align} \label{idealequationg121}
0 \longrightarrow \mathcal{I}_Z \longrightarrow \mathcal{I}_{\Gamma_p} \longrightarrow (\id,p)_{\ast}\mathcal{Q} \longrightarrow 0.
\end{align}
\end{lemma}

\begin{proof}
    This follows from the Snake Lemma using the middle column, respectively the third row in the defining diagram of $\mathcal{O}_Z$.\ Note that $(p \times \id)^{\ast}\mathcal{I}_{\Delta}^2 = \mathcal{I}_{2\Gamma_p}$ and $(p \times \id)^{\ast}\mathcal{I}_{\Delta} = \mathcal{I}_{\Gamma_p}$. 
\end{proof}

\begin{prop} \label{iclosedemb}
The morphism $\iota: \mathbb{G} \rightarrow X^{[n]}$ with $n = l+1$ is a closed embedding with set-theoretical image given by planar subschemes of $X$.\ That is, $\iota$ maps a $k$-point $(P,I)$ lying inside $\lbrace P \rbrace \times G(l,\mathfrak{m}/\mathfrak{m}^2) \subseteq \mathbb{G}$ to the unique zero-dimensional subscheme $[W]$ in $X^{[n]}(k)$ with support $P \in X$ and defining ideal $I$ (modulo $\mathfrak{m}^2$).
\end{prop}

\begin{proof}
(1) Let $p:\mathbb{G} \rightarrow X$ be the projection and consider the cartesian diagram
\begin{center}
\begin{tikzcd}
W \arrow[d] \arrow[r]          & Z \arrow[d, "q"] \arrow[r]  & \Xi_{n} \arrow[d, "\pi"] \\
\spec k \arrow[r, "{f_{P,I}}"] & \mathbb{G} \arrow[r, "\iota"]  & {X^{[n]}}          
\end{tikzcd}
\end{center}
with $f_{P,I}$ denoting the inclusion of $(P,I)$ into $\mathbb{G}$ and $W = q^{-1}(\spec k)$.\ Observe that $W$ is the image point $\iota(P,I) = [W] \in X^{[n]}(k)$ and that $W_{\red} \subseteq \Gamma_p$ is just the reduced $P \in X$.\ It remains to see that the first order part of the defining ideal $\mathcal{I}_{W/X}$ is given by $I$. \par 
\noindent Restricting the sequence (\ref{idealequationg120}) to $\lbrace (P,I) \rbrace \times X \cong X$ via $f_{P,I} \times \id_X$ yields (by flatness of $(\id,p)_{\ast}\mathcal{K}$ over $\mathbb{G}$) the short exact sequence
\begin{align*}
0 \longrightarrow \bigl((p \circ f_{P,I}) \times \id\bigr)^{\ast}\mathcal{I}_{\Delta}^2 \longrightarrow \mathcal{I}_{W/X} \longrightarrow (f_{P,I} \times \id)^{\ast}(\id,p)_{\ast}\mathcal{K} \cong \mathcal{K}(P,I) \longrightarrow 0.
\end{align*}
That $\mathcal{I}_{W/X} =  (f_{P,I} \times \id)^{\ast}\mathcal{I}_{Z}$ follows from Corollary \ref{basechangeflatideal2}, which also implies that the sheaf on the left is $\mathcal{I}_{P/X}^2$.\ Going to the stalks yields $\mathcal{I}_{W,P}/\mathfrak{m}_P^2 \cong \mathcal{K}(P,I) \cong I$. \par
\noindent (2) To show that $\iota$ is a closed embedding, note that $Z$ is even a closed subscheme of $\mathbb{G} \times_X n\Delta \subseteq \mathbb{G} \times X$ so that $\iota$ factors through $X^{[n]}_{(n),\red} = (\hilb^n_{n\Delta/X})_{\red}$ (see Section \ref{preliminarieshilbert}). \par
\noindent The $X$-morphism $\mathbb{G} \rightarrow X^{[n]}_{(n),\red}$ can be checked to be a closed embedding fibrewise using \cite[Prop.\ 12.93]{goertzwedhorn}, where it coincides with the obvious inclusion of planar subschemes $Z_{(1,l,0)} = G(l,\mathfrak{m}/\mathfrak{m}^2) \hookrightarrow \hilb_d^n$ (Section \ref{preliminarieshilbert}, (\ref{formalcomputation})) according to the first half of this proof.
\end{proof} 

\subsection{Computation of the normal bundle}

\noindent For various reasons it is important to have a description of the normal bundle of $\mathbb{G} = \mathbb{G}(l,\Omega_X)$ in $X^{[n]} = X^{[l+1]}$, making only sense in the smooth situation $l \in \lbrace 1,2 \rbrace$.\ As a first step, the restriction of $\Omega_{X^{[n]}/k}$ to $\mathbb{G}$ is computed. 

\begin{lemma} \label{defirelativeext}
Let $f:X \rightarrow Y$ be projective. The $i$-th relative $\mathcal{E}xt$-sheaf on $Y$ is 
\begin{align*}
\mathcal{E}xt_f^i(\mathcal{F},\mathcal{G}) := R^i(f_{\ast}\mathcal{H}om_{\mathcal{O}_X}(\mathcal{F},\mathcal{G})) := R^i\bigl(f_{\ast}\mathcal{H}om_{\mathcal{O}_X}(\mathcal{F},-)\bigr)(\mathcal{G})
\end{align*}
for $\mathcal{F},\mathcal{G} \in \mathbf{Coh}(X)$ and can equivalently be computed as $\mathcal{H}^i\bigl(Rf_{\ast}R\mathcal{H}om_{\mathcal{O}_X}(\mathcal{F},\mathcal{G})\bigr)$.
\begin{enumerate}
\item If $Y = \spec R$ is affine, then $\mathcal{E}xt_f^i(\mathcal{F},\mathcal{G})$ is the quasi-coherent sheaf associated to the $R$-module $\ext^i_{\mathcal{O}_X}(\mathcal{F},\mathcal{G})$.
\item If $\mathcal{F}$ and $\mathcal{G}$ are coherent, then so is $\mathcal{E}xt_f^i(\mathcal{F},\mathcal{G})$.
\item If $\mathcal{F}$ is locally free, $\mathcal{E}xt_f^i(\mathcal{F},\mathcal{G}) \cong R^if_{\ast}(\mathcal{F}^{\vee} \otimes \mathcal{G})$.\ Thus $\mathcal{E}xt_f^i(\mathcal{O}_X,\mathcal{G}) =R^if_{\ast}\mathcal{G}$.
\end{enumerate}
\end{lemma}

\begin{proof}
\noindent The statements (1) - (3) can be found in \cite{langeext} and are easy to show.\ The identification $R(f_{\ast}\mathcal{H}om_{\mathcal{O}_X}(\mathcal{F},-)) \cong Rf_{\ast}R\mathcal{H}om_{\mathcal{O}_X}(\mathcal{F},-)$ uses that $\mathcal{H}om_{\mathcal{O}_X}(\mathcal{F},-)$ sends injective sheaves to flasque sheaves, cf.\ the proof of \cite[Prop.\ II.5.3]{hartshorneresidues}.
\end{proof}

\noindent There is a base change theory for relative $\mathcal{E}xt$'s analogous to that in \cite[Ch.\ III.12]{hartshorne} relating higher direct images to cohomology of the fibres, see e.g.\ \cite{langeext}. 

\begin{kor} \label{relextvanishing}
Let $f:X \rightarrow Y$ be as above and additionally smooth with equidimensional fibres of dimension $d$.\ If $\mathcal{F}$ and $\mathcal{G}$ are flat over $Y$, then
\begin{align*}
\mathcal{E}xt^i_f(\mathcal{F},\mathcal{G}) = 0 ~~ \text{ for } ~~ i > d.
\end{align*}
\end{kor}

\begin{proof}
All $\ext^i_{X_y}(\mathcal{F}_y,\mathcal{G}_y)$'s vanish for $i > d$ by Serre-duality. Thus $\mathcal{E}xt_f^i(\mathcal{F},\mathcal{G}) \otimes_{\mathcal{O}_Y} k(y) = 0$ (invoking the base change theorem \cite[Thm.\ 1.4]{langeext}) for all $y \in Y$. 
\end{proof}

\noindent The following proposition appears in \cite[Rem.\ 3.7]{lehncotangent} for Quot schemes without proof.

\begin{prop} \label{descriptionrestrictedcotangentsheaf}
Let $g:T \rightarrow X^{[n]}$ be the classifying morphism associated to a flat family $q: Z \rightarrow T$ of degree $n = l+1$, fitting into a commutative diagram of the following form.
\begin{center}
\begin{tikzcd}
Z \arrow[r, hook] \arrow[rd, "q"'] & T \times X \arrow[r, "g \times \id"] \arrow[d, "{\overline{q}}"'] \arrow[rd, "\square", phantom] & {X^{[n]} \times X} \arrow[d, "{\overline{\pi}}"] & \Xi_{n} =: \Xi \arrow[l, hook'] \arrow[ld, "\pi"] \\
    & T \arrow[r, "g"]                                                                    & {X^{[n]}}                           &                                                    
\end{tikzcd}
\end{center}
\noindent Then there is a base change isomorphism $g^{\ast}\Omega_{X^{[n]}/k} \cong \mathcal{E}xt_{\overline{q}}^d(\mathcal{O}_Z,\mathcal{I}_{Z} \otimes \omega_{\overline{q}})$.
\end{prop}

\begin{proof}
By the main theorem of \cite{lehncotangent}, $\Omega_{X^{[n]}/k} \cong \mathcal{E}xt^d_{\overline{\pi}}(\mathcal{O}_{\Xi},\mathcal{I}_{\Xi} \otimes \omega_{\overline{\pi}})$, see Proposition \ref{tangentsheafhilbert}.\ By flatness of $\overline{\pi}$, derived base change \cite[Prop.\ A.85]{fouriernahm} as well as Proposition \ref{basechangeflatideal},
\begin{align*}
Lg^{\ast}R\overline{\pi}_{\ast}R\mathcal{H}om(\mathcal{O}_{\Xi},\mathcal{I}_{\Xi} \otimes \omega_{\overline{\pi}}) &\cong R\overline{q}_{\ast}L(g \times \id)^{\ast}R\mathcal{H}om(\mathcal{O}_{\Xi},\mathcal{I}_{\Xi} \otimes \omega_{\overline{\pi}}) \\
& \cong R\overline{q}_{\ast}R\mathcal{H}om\bigl(L(g \times \id)^{\ast}\mathcal{O}_{\Xi},L(g \times \id)^{\ast}(\mathcal{I}_{\Xi} \otimes \omega_{\overline{\pi}})\bigr) \\
& \cong R\overline{q}_{\ast}R\mathcal{H}om(\mathcal{O}_{Z},\mathcal{I}_{Z} \otimes \omega_{\overline{q}}).
\end{align*}
Let $E_2^{r,s} = \mathcal{H}^r\bigl(Lg^{\ast} \mathcal{E}xt^s_{\overline{\pi}}(\mathcal{O}_{\Xi},\mathcal{I}_{\Xi} \otimes \omega_{\overline{\pi}})\bigr)$.\ There exists a spectral sequence (\cite[(3.10)]{huybrechts}) 
\begin{align*}
E_2^{r,s} \Rightarrow  \mathcal{H}^{r+s}\bigl(Lg^{\ast}R\overline{\pi}_{\ast}R\mathcal{H}om(\mathcal{O}_{\Xi},\mathcal{I}_{\Xi} \otimes \omega_{\overline{\pi}})\bigr) \cong \mathcal{E}xt^{r+s}_{\overline{q}}(\mathcal{O}_{Z},\mathcal{I}_{Z} \otimes \omega_{\overline{q}}).
\end{align*}
All fibres of $\overline{\pi}$ are isomorphic to $X$ and thus $E_2^{r,s} = 0$ for $r > 0$ or $s > d$ by Corollary \ref{relextvanishing}.\ Therefore $E_2^{0,d} = E_{\infty}^{0,d}$ and $\mathcal{E}xt^d_{\overline{q}}(\mathcal{O}_{Z},\mathcal{I}_{Z} \otimes \omega_{\overline{q}}) \cong E_2^{0,d}$, which is by definition the pullback of $\mathcal{E}xt^d_{\overline{\pi}}(\mathcal{O}_{\Xi},\mathcal{I}_{\Xi} \otimes \omega_{\overline{\pi}}) \cong \Omega_{X^{[n]}/k}$ along $g$. 
\end{proof}

\noindent Let $\overline{q}: \mathbb{G} \times X \rightarrow \mathbb{G}$ be the projection and let $\omega_{\overline{q}} = \omega_{\mathbb{G} \times X/\mathbb{G}}$.\ Proposition \ref{descriptionrestrictedcotangentsheaf} yields an isomorphism $(\Omega_{X^{[n]}/k})_{|\mathbb{G}} \cong \mathcal{E}xt_{\overline{q}}^d(\mathcal{O}_Z,\mathcal{I}_{Z} \otimes \omega_{\overline{q}})$ of locally free sheaves after considering 
\begin{align*}
    (g:T \rightarrow X^{[n]}) = (\iota: \mathbb{G} \hookrightarrow X^{[l+1]} = X^{[n]}).
\end{align*}
Recall the short exact sequence $0 \rightarrow (\id,p)_{\ast}\mathcal{Q} \rightarrow \mathcal{O}_Z \rightarrow (\id,p)_{\ast}\mathcal{O}_{\mathbb{G}} \rightarrow 0$ on $\mathbb{G} \times X$, yielding the distinguished triangle 
\begin{align*}
    R\overline{q}_{\ast}R\mathcal{H}om((\id,p)_{\ast}\mathcal{O}_{\mathbb{G}},\mathcal{I}_{Z}  \otimes \omega_{\overline{q}}) & \longrightarrow R\overline{q}_{\ast}R\mathcal{H}om(\mathcal{O}_Z,\mathcal{I}_{Z} \otimes \omega_{\overline{q}}) \\ &\longrightarrow R\overline{q}_{\ast}R\mathcal{H}om((\id,p)_{\ast}\mathcal{Q},\mathcal{I}_{Z} \otimes \omega_{\overline{q}}) \overset{+1}{\longrightarrow} \cdot
\end{align*}
that will be combined below with the result of applying $R\overline{q}_{\ast}R\mathcal{H}om((\id,p)_{\ast}\mathcal{Q},(-) \otimes \omega_{\overline{q}})$ to  (\ref{idealequationg120}) and with the result of applying $R\overline{q}_{\ast}R\mathcal{H}om((\id,p)_{\ast}\mathcal{O}_{\mathbb{G}},(-) \otimes \omega_{\overline{q}})$ to (\ref{idealequationg121}).\ Note that 
\begin{align*}
R{\overline{q}}_{\ast}R\mathcal{H}om((\id,p)_{\ast}\mathcal{E},(-) \otimes \omega_{\overline{q}}) &\cong R{\overline{q}}_{\ast}(\id,p)_{\ast}R\mathcal{H}om\bigl(\mathcal{E},(\id,p)^!((-) \otimes \omega_{\overline{q}})\bigr) \\
& \cong R\mathcal{H}om(\mathcal{E},L(\id,p)^{\ast}(-)[-d])) \\
& \cong \mathcal{E}^{\vee} \otimes L(\id,p)^{\ast}(-)[-d]
\end{align*}
for any locally free sheaf $\mathcal{E}$ by Grothendieck--Verdier duality.\ Hence the long exact sequences in cohomology yield a diagram of exact sequences on $\mathbb{G}$ of the following form.
\begin{center}
\begin{tikzcd}
 & 0                                                                                 &                                                  &                       &   \\
{\mathcal{Q}^{\vee} \otimes (\id,p)^{\ast}(\mathcal{I}_{2\Gamma_p})} \arrow[r] & {\mathcal{Q}^{\vee} \otimes (\id,p)^{\ast}\mathcal{I}_Z} \arrow[r] \arrow[u]      & \mathcal{Q}^{\vee} \otimes \mathcal{K} \arrow[r] & 0                     &   \\
 & {\mathcal{E}xt^d_{\overline{q}}(\mathcal{O}_Z,\mathcal{I}_Z \otimes \omega_{\overline{q}})}  \arrow[u] &                                                  &                       &   \\
& {(\id,p)^{\ast}\mathcal{I}_Z} \arrow[r] \arrow[u]                                 & {(\id,p)^{\ast}\mathcal{I}_{\Gamma_p}} \arrow[r] & \mathcal{Q} \arrow[r] & 0
\end{tikzcd}
\end{center}

\begin{lemma}
There are isomorphisms $(\id,p)^{\ast}\mathcal{I}_{n\Gamma_p} = (\id,p)^{\ast}(\mathcal{I}_{\Gamma_p}^n) \cong S^n p^{\ast}\Omega_{X/k}$.
\end{lemma}

\begin{proof}
Recall from Proposition \ref{basechangeflatideal2} that $\mathcal{I}_{n\Gamma_p} = \mathcal{I}_{\Gamma_p}^n \cong (p \times \id)^{\ast}\mathcal{I}_{\Delta}^n$. Note also that 
\begin{align*}
\Delta^{\ast}\mathcal{I}_{\Delta}^n \cong \Delta^{\ast}(\mathcal{I}_{\Delta}^n/\mathcal{I}_{\Delta}^{n+1})
\end{align*}
where $\Delta$ is the diagonal morphism. As $X$ is smooth, $\mathcal{I}_{\Delta}^n/\mathcal{I}_{\Delta}^{n+1} \cong S^n(\mathcal{I}_{\Delta}/\mathcal{I}_{\Delta}^2)$. So 
\begin{align*}
(\id,p)^{\ast}\mathcal{I}_{n\Gamma_p} \cong (\id,p)^{\ast}(p \times \id)^{\ast}\mathcal{I}_{\Delta}^n \cong p^{\ast}\Delta^{\ast}(\mathcal{I}_{\Delta}^n/\mathcal{I}_{\Delta}^{n+1}) \cong p^{\ast}S^n\Delta^{\ast}(\mathcal{I}_{\Delta}/\mathcal{I}_{\Delta}^2) \cong S^n p^{\ast}\Omega_{X/k},
\end{align*}
finishing the proof of the lemma.
\end{proof}

\noindent Lemma \ref{rankofpullbackI} and Lemma \ref{dimensionIMI} are similar to \cite[Lemma 4.6]{krugrennemo} and \cite[Lemma 4.5]{krugrennemo}.

\begin{lemma} \label{rankofpullbackI}
The sheaf $(\id,p)^{\ast}\mathcal{I}_Z$ is locally free of rank $d + \frac{1}{2}(l^2 - l)$.
\end{lemma}

\begin{proof}
On a reduced scheme, equidimensionality of the fibres of a coherent sheaf allows to conclude that it is locally free of the expected rank \cite[Exc.\ II.5.8]{hartshorne}.\ In fact, looking at fibres of closed points is enough:\ The locus where a coherent sheaf has a specified constant rank is locally closed and $k$-points form a very dense subset - hence the reasoning of \cite[Rem.\ 10.15]{goertzwedhorn} may be applied. \par 
\noindent Let $f_{P,I}: \spec k \rightarrow \mathbb{G}$ be a $k$-point corresponding to a pair $(P,I)$ (the notation is the same as in Proposition \ref{iclosedemb}).\ The fibre of $(\id,p)^{\ast}\mathcal{I}_Z$ at $(P,I)$ is the same as the fibre of $\mathcal{I}_W = \mathcal{I}_{W/X}$ restricted to the $k$-point $P$.\ Thus it remains to compute the dimension of $\mathcal{I}_W(P) \cong \mathcal{I}_{W,P}/\mathfrak{m}_{P}\mathcal{I}_{W,P}$ for a zero-dimensional subscheme $W \subseteq X$ with Hilbert function $(1,l,0)$ concentrated at a point $P \in X$.\  The identification $
\widehat{\mathcal{O}_{X,P}}/\widehat{\mathfrak{m}_{X,P}}^3 \cong k[x_1,\hdots,x_d]/\mathfrak{m}^3$ preserves $I/\mathfrak{m}I$ and $P = 0 \in \mathbb{A}^d$ can be assumed, where $I$ corresponds to $\mathcal{I}_{W,P}$.\ The computation of Lemma \ref{dimensionIMI} below concludes the proof.
\end{proof}

\begin{bem}
    Note that $(\id,p)^{\ast}\mathcal{I}_Z$ is the restriction of the conormal sheaf $(\mathcal{C}_{Z/\mathbb{G} \times X})_{|\Gamma_p}$.
\end{bem}

\begin{lemma} \label{dimensionIMI}
Let $I = (u_{l+1},\hdots,u_d) + \mathfrak{m}^2 \trianglelefteq k[x_1,\hdots,x_d]$ be an ideal with Hilbert function $(1,l,0)$ and linear forms $u_i$.\ Then $I/I\mathfrak{m} \cong I/\mathfrak{m}^2 \oplus \mathfrak{m}^2/I\mathfrak{m}$ and $\dime_k I/I\mathfrak{m} = d + \frac{1}{2}(l^2 - l)$.
\end{lemma}

\begin{proof}
Given the homogeneous ideal $I \trianglelefteq k[x_1,\hdots,x_d]$, the linearly independent linear forms $u_{l+1}, \hdots, u_d$ can be extended to a $k$-basis $u_1,\hdots,u_d$ of $k[x_1,\hdots,x_d]_1 \cong_k \mathfrak{m}/\mathfrak{m}^2$. It is easy to reduce to the case $u_i = x_i$.\ There is a short exact sequence of $k[x_1,\hdots,x_d]$-modules 
\begin{align*}
0 \rightarrow \mathfrak{m}^2/I\mathfrak{m} \rightarrow I/I\mathfrak{m} \rightarrow I/\mathfrak{m}^2 \rightarrow 0,
\end{align*}
of course split as a sequence of $k$-vector spaces (in fact, it is even split as a sequence of $k[x_1,\hdots,x_d]$-modules).\ The dimension of $I/\mathfrak{m}^2$ is $d-l$ by definition of the prescribed Hilbert function.\ The surviving monomials in $\mathfrak{m}^2/I\mathfrak{m}$ with $I\mathfrak{m} = (x_{l+1},\hdots,x_d)\mathfrak{m} + \mathfrak{m}^3$ are all those $\binom{l+1}{2} = \frac{1}{2}(l^2 + l)$ monomials of degree $2$ that only contain $x_1,\hdots,x_l$. 
\end{proof}

\noindent To sum up, there is a diagram of exact sequences of \emph{vector bundles} on $\mathbb{G}$ of the form
\begin{center}
\begin{tikzcd}
  & 0                                                                                                                      &                                                  &                       &   \\
\mathcal{Q}^{\vee} \otimes S^2p^{\ast}\Omega_{X/k} \arrow[r, "\id \otimes \alpha"] & {\mathcal{Q}^{\vee} \otimes (\id,p)^{\ast}\mathcal{I}_Z} \arrow[u] \arrow[r, "\id \otimes \beta"] & \mathcal{Q}^{\vee} \otimes \mathcal{K} \arrow[r] & 0                     &   \\
& {(\Omega_{X^{[l+1]}/k})_{|\mathbb{G}}} \arrow[u]                                                                       &                                                  &                       &   \\
& { (\id,p)^{\ast}\mathcal{I}_Z} \arrow[r] \arrow[u]                                                                     & p^{\ast}\Omega_{X/k} \arrow[r]                   & \mathcal{Q} \arrow[r] & 0.
\end{tikzcd}
\end{center}

\begin{bem} 
In the case $l=1$ of $X^{[2]}$, the middle column is already exact with a zero on the bottom. In the case $l=2$ of $X^{[3]}$, there is a short exact sequence
\begin{align*}
    0 \rightarrow \mathcal{K} \rightarrow (\Omega_{X^{[3]}/k})_{|\mathbb{G}} \rightarrow \mathcal{Q}^{\vee} \otimes (\id,p)^{\ast}\mathcal{I}_Z \rightarrow 0.
\end{align*}
\noindent This is because counting ranks of locally free sheaves for $l = 1$, the middle column is exact since $d + 1 \cdot d = 2d = \dime X^{[2]}$. For $l=2$, a morphism $\mathcal{K} \rightarrow (\Omega_{X^{[3]}/k})_{|\mathbb{G}}$ has to be constructed first using (\ref{idealequationg120}) in a way such that $S^2p^{\ast}\Omega_{X/k} \rightarrow (\id,p)^{\ast}\mathcal{I}_Z \rightarrow (\Omega_{X^{[3]}/k})_{|\mathbb{G}}$ vanishes on the fibres.\ Then $(d-2) + 2\cdot (d+1) = 3d = \dime X^{[3]}$.\ The proof is omitted here since this remark is not needed later on.
\end{bem}

\noindent The next task is to understand the kernel of the map $\id \otimes \beta$.\ Recall that the long exact sequence of interest (forgetting $\mathcal{Q}^{\vee} \otimes (-)$ for a moment) is obtained from (\ref{idealequationg120}) as
\begin{align*}
\mathcal{K} \otimes p^{\ast}\Omega_{X/k} \cong L^1(\id,p)^{\ast}(\id,p)_{\ast}\mathcal{K} \overset{m}{\longrightarrow} S^2p^{\ast}\Omega_{X/k} \overset{\alpha}{\longrightarrow} (\id,p)^{\ast}\mathcal{I}_Z \overset{\beta}{\longrightarrow} \mathcal{K} \longrightarrow 0
\end{align*}
by pullback along $(\id,p)$, where the isomorphism on the left hand side is a consequence of the following lemma together with the identification $\mathcal{N}_{\Gamma_p / \mathbb G \times X}^{\vee} \cong p^{\ast}\mathcal{N}_{\Delta/X \times X}^{\vee} \cong p^{\ast}{\Omega_{X/k}}$.

\begin{lemma}[{\cite[Prop.\ A.6]{caldararuks}, \cite[Prop.\ 11.8]{huybrechts}}] \label{formulapushpull}
Let $j: Z \hookrightarrow Y$ be a closed embedding of smooth varieties and let $\mathcal{E}$ be a locally free sheaf on $Z$.\ Then for all $i \geq 0$, there are natural isomorphisms $L^ij^{\ast}j_{\ast}\mathcal{E} \cong \mathcal{E} \otimes \wedge^i\mathcal{N}_{Z/Y}^{\vee}$ and $\mathcal{H}^i(j^!j_{\ast}\mathcal{E}) \cong \mathcal{E} \otimes \wedge^i\mathcal{N}_{Z/Y}$.
\end{lemma}

\begin{prop} \label{multiplicationmap}
The morphism $m$ above is the natural multiplication map and there is a short exact sequence of locally free sheaves of ranks $\binom{l+1}{2}, d + \frac{1}{2}(l^2 - l)$, resp. $d - l$:
\begin{align*}
0 \longrightarrow S^2\mathcal{Q} \longrightarrow (\id,p)^{\ast}\mathcal{I}_Z \overset{\beta}{\longrightarrow} \mathcal{K} \longrightarrow 0.
\end{align*}
\end{prop}

\begin{proof}
    The second assertion follows from the first, for if $m$ is the multiplication map with image $\mathcal{K}\cdot p^{\ast}\Omega_{X/k}$, its cokernel ($\cong \ker \beta$) is $S^2\mathcal{Q}$ since there exists a short exact sequence 
    \begin{align*}
        0 \longrightarrow \mathcal{K}\cdot p^{\ast}\Omega_{X/k} \longrightarrow S^2p^{\ast}\Omega_{X/k} \longrightarrow S^2\mathcal{Q} \longrightarrow 0
    \end{align*}
    resulting from (\ref{tautsequence}), cf. \cite[\href{https://stacks.math.columbia.edu/tag/01CJ}{01CJ}]{stacks-project}.\ The short exact sequences (\ref{idealequationg120}) and (\ref{idealequationg121}) corresponding to $\mathcal{I}_{2\Gamma_p} \subseteq \mathcal{I}_Z \subseteq \mathcal{I}_{\Gamma_p}$ can be combined to a big commutative diagram using the defining sequence of the sheaf of differentials. Considering the long exact sequences for the derived pullback $L(\id,p)^{\ast}$ applied to this diagram, it follows with the aid of Lemma \ref{formulapushpull} that there is an even bigger commutative diagram of the following form.
\begin{center}
\begin{tikzcd}
0 \arrow[d, dashed]                                                                             &                                                                             & \mathcal{Q} \otimes p^{\ast}\Omega_{X/k} \arrow[d] \arrow[r, dashed]         & 0 \arrow[d]                               &    \\
\mathcal{K} \otimes p^{\ast}\Omega_{X/k} \arrow[r, "m"] \arrow[d, dashed]  & S^2p^{\ast}\Omega_{X/k} \arrow[r,"\alpha"] \arrow[d, "="]  & {(\id,p)^{\ast}\mathcal{I}_Z} \arrow[d] \arrow[r,"\beta"]  & \mathcal{K} \arrow[r] \arrow[d]           & 0  \\

p^{\ast}\Omega_{X/k} \otimes p^{\ast}\Omega_{X/k} \arrow[r,"m'"] \arrow[rruu, dashed, bend left]     & S^2p^{\ast}\Omega_{X/k} \arrow[r, "0"]                                      & p^{\ast}\Omega_{X/k} \arrow[d] \arrow[r, "="]                                & p^{\ast}\Omega_{X/k} \arrow[r] \arrow[d] & 0 \\
&                                                                             & \mathcal{Q} \arrow[d]  \arrow[r, "="]                                                        &       \mathcal{Q} \arrow[d]                                    &    \\
&                                                                             & 0                                                                            &     0                                      &   
\end{tikzcd}
\end{center}
The morphisms between the tautological sheaves on $\mathbb{G}$ in the right column a priori constitute \enquote{some} natural non-split cokernel sequence.\par 
\noindent Using that on each fibre of $p$, $\ext^1(Q,K)$ is one-dimensional (cf.\ the methods used in Appendix \ref{AppendixA}), all morphisms between the tautological sheaves come from the tautological sequence (\ref{tautsequence}) and are thus the expected ones. \par 
\noindent The claim reduces to $m'$ being the canonical projection.\ Since $L(\id,p)^{\ast}(p \times \id)^{\ast} \cong p^{\ast}L\Delta^{\ast}$, it is enough to show that the morphism $\partial$ in the long exact $L\Delta^{\ast}$-sequence 
\begin{align*}
\hdots \longrightarrow L^1\Delta^{\ast}\mathcal{I}_{\Delta} \cong \wedge^2\Omega_{X/k} \overset{L^1\Delta^{\ast}\pi}{\longrightarrow} \Omega_{X/k}^{\otimes 2} \overset{\partial}{\longrightarrow} S^2\Omega_{X/k} \overset{0}{\longrightarrow} \Omega_{X/k} \overset{\sim}{\longrightarrow} \Omega_{X/k} \longrightarrow 0
\end{align*}
associated to $0 \rightarrow \mathcal{I}_{\Delta}^2 \rightarrow \mathcal{I}_{\Delta} \overset{\pi}{\rightarrow} \Delta_{\ast}\Omega_{X/k} \rightarrow 0$ on $X \times X$ is the canonical projection.\ Equivalently, $L^1\Delta^{\ast}\pi$ is the monomorphism corresponding to $a \wedge b \mapsto a \otimes b - b \otimes a$ with cokernel $S^2\Omega_{X/k}$.\footnote{This works in any characteristic of the base field.}\ This can be checked on the fibres, where the claim reduces to understanding 
\begin{align*}
\tor^1_{\mathcal{O}_{X,P}}(\mathfrak{m},\mathcal{O}_{X,P}/\mathfrak{m}) \longrightarrow  \tor^1_{\mathcal{O}_{X,P}}(\mathfrak{m}/\mathfrak{m}^2,\mathcal{O}_{X,P}/\mathfrak{m}).
\end{align*}
More elegantly (and without the need of reducing from $L(\id,p)^{\ast}$ to $L\Delta^{\ast}$), one applies the global Lemma \ref{technicalpullbacksymmetric} to $i_{\ast}\mathcal{N}_{Y/X}^{\vee} = (\id,p)_{\ast}(\mathcal{N}_{\Gamma_p/\mathbb{G} \times X}^{\vee})$.
\end{proof}

\noindent All in all, the following diagram of exact sequences is obtained for $l \in \lbrace 1,2 \rbrace$ after invoking the relative cotangent sequence for $\mathbb{G} \overset{p}{\rightarrow} X \rightarrow \spec k$, the identity $\Omega_{p} \cong \mathcal{Q}^{\vee} \otimes \mathcal{K}$ (see Proposition \ref{tangentsheafgrass}) and the conormal sequence for $\mathbb{G} \subseteq X^{[l+1]} = X^{[n]}$ with $\mathcal{N} = \mathcal{N}_{\mathbb{G}/X^{[n]}}$.

\begin{center}
\begin{tikzcd}
& 0   &        &        &   \\
& \mathcal{Q} \arrow[u]     & 0     & 0              &   \\
0 \arrow[r] & p^{\ast}\Omega_{X/k} \arrow[r, "(dp)^{\ast}"] \arrow[u] \arrow[rd, "(a)", phantom] & \Omega_{\mathbb{G}/k} \arrow[r] \arrow[u] \arrow[rd, "(b)", phantom] & \Omega_p \arrow[u] \arrow[r]                   & 0 \\
& { (\id,p)^{\ast}\mathcal{I}_Z} \arrow[u] \arrow[r]   & {(\Omega_{X^{[l+1]}/k})_{|\mathbb{G}}} \arrow[r] \arrow[u, "(d\iota)^{\ast}"] & {\mathcal{Q}^{\vee} \otimes (\id,p)^{\ast}\mathcal{I}_Z} \arrow[r] \arrow[u] & 0 \\   &                    & \mathcal{N}^{\vee} \arrow[u]                                        & \mathcal{Q}^{\vee} \otimes S^2\mathcal{Q} \arrow[u]                &   \\   &              & 0 \arrow[u]              & 0 \arrow[u]      &  
\end{tikzcd}
\end{center}

\begin{thrm}[Thm.\ \ref{maintheoremnormalbundleshort}] \label{maintheoremnormalbundle}
The above diagram commutes.\ So there is an exact sequence 
\begin{align*}
S^2\mathcal{Q} \longrightarrow \mathcal{N}^{\vee} \longrightarrow \mathcal{Q}^{\vee} \otimes S^2\mathcal{Q} \longrightarrow \mathcal{Q} \longrightarrow 0.
\end{align*}
Hence if $l = 1$, the normal bundle of $\mathbb{G} = \mathbb{G}(l,\Omega_{X/k})$ inside $X^{[2]}$ is isomorphic to $S^2\mathcal{Q}^{\vee}$, which under the identification $\mathbb{G} = \mathbb{P}(\Omega_{X/k}) = X^{[2]}_{(2),\red}$ is the line bundle $\mathcal{O}_p(-2)$ in accordance with the already known result \cite[Thm.\ 4.1(ii)]{krugploogsosna}.\ If $l = 2$, then the normal bundle of $\mathbb{G}$ inside $X^{[3]}$ has rank $4$ and is part of the short exact sequence
\begin{align} \label{definingsequenceN}
0 \longrightarrow \mathcal{Q}^{\vee} \longrightarrow \mathcal{Q} \otimes S^2\mathcal{Q}^{\vee} \longrightarrow \mathcal{N} \longrightarrow 0.
\end{align}
\end{thrm}

\begin{proof}
The only things to prove are the commutativity of (b) and of (a) since after identifying the kernel of $(\id,p)^{\ast}\mathcal{I}_Z \rightarrow p^{\ast}\Omega_{X/k}$ with $S^2\mathcal{Q}$, the Snake Lemma can be applied:  \par \vspace{3pt}
\noindent Namely if $l = 1$, $\mathcal{Q}^{\vee} \otimes S^2\mathcal{Q} \rightarrow \mathcal{Q}$ is an epimorphism between locally free sheaves of the same rank, hence an isomorphism.\ The same argument applies to $S^2\mathcal{Q} \rightarrow \mathcal{N}^{\vee} \rightarrow 0$ and dualizing yields the claim.\ If $l = 2$, then $\mathcal{N}^{\vee}$ is locally free of rank $4 = \codim(\mathbb{G},X^{[3]})$.\ Thus counting ranks and dualizing afterwards yields (\ref{definingsequenceN}).\ Consider the square (b) now. \par \vspace{3pt}
\noindent It is sufficient to show commutativity of a square after restricting to fibres of $p$ because commutativity can be checked on the fibres at closed points.\ The restriction of $\Omega_{\mathbb{G}/k} \rightarrow \Omega_p$ to a fibre $G := G(l,\mathfrak{m}/\mathfrak{m}^2) \cong p^{-1}(P)$ is the same as the codifferential of the inclusion $G \hookrightarrow \mathbb{G}$.\ The chain of inclusions $G \hookrightarrow \mathbb{G} \hookrightarrow X^{[n]}$ coincides with $G \hookrightarrow \hilb^{n}_d \hookrightarrow X^{[n]}$, where as usual $n= l+1$.\ The claim for (b) is that dualizing 
\begin{align} \label{sequenceofinterest}
    (\Omega_{X^{[n]}/k})_{|\mathbb{G}} \cong \mathcal{E}xt^d_{\overline{q}}(\mathcal{O}_Z,\mathcal{I}_Z \otimes \omega_{\overline{q}}) \longrightarrow \mathcal{Q}^{\vee} \otimes (\id,p)^{\ast}\mathcal{I}_Z \longrightarrow \mathcal{Q}^{\vee} \otimes \mathcal{K} \cong \Omega_p
\end{align}
and restricting to $[(P,I)] \in \mathbb{G}(k)$ yields a sequence of vector spaces 
\begin{align} \label{restrictedhomseq}
    \homo_k(I/\mathfrak{m}^2,\mathfrak{m}/I) \rightarrow \homo_{\mathcal{O}_X}(I,\mathfrak{m}/I) \rightarrow \homo_{\mathcal{O}_X}(I,\mathcal{O}_X/I)
\end{align}
agreeing precisely with the composition of the differentials of $p^{-1}(P) \hookrightarrow \hilb^{n}_d \hookrightarrow X^{[n]}$.\footnote{This is equal to the composition $p^{-1}(P) \hookrightarrow (\spec \mathcal{O}_{X,P}/\mathfrak{m}^n)^{[n]} \hookrightarrow X^{[n]}$.}\par \vspace{3pt}
\noindent \emph{Proof of the claim.}\ The differentials are described first.\ There is no harm in restricting to an open affine $\spec A \subseteq X$ containing $P = [\mathfrak{p}]$ where $I$ defines a zero-dimensional subscheme $W$ of length $n$. Let $\overline{A} = A_{\mathfrak{p}}/(\mathfrak{p}A_{\mathfrak{p}})^n = \mathcal{O}/\mathfrak{m}^n$ and let $\overline{I} \trianglelefteq \overline{A}$ be the corresponding ideal.\ The three tangent space identifications (see Prop.\ \ref{tangentsheafgrass}(4) and Prop.\ \ref{tangentsheafhilbert}(3))
\begin{gather*}
    T_{[I/\mathfrak{m}^2]}G \cong \homo_k(I/\mathfrak{m}^2,\mathfrak{m}/I),\quad  T_{[W]}X^{[n]} \cong \homo_{A}(I,A/I) \quad \text{and}\\  T_{[W]}\bigl(\spec \mathcal{O}_{X,P}/\mathfrak{m}^n\bigr)^{[n]} \cong \homo_{\overline{A}}(\overline{I},A/I)
\end{gather*}
 can be applied to the following commutative diagram of vector spaces.
\begin{center}
\begin{tikzcd}
 & {T_{[I/\mathfrak{m}^2]}G} \arrow[ld, hook] \arrow[d, hook] \arrow[rd, hook] &                                                 \\
{T_{[W]}X^{[n]}} & T_{[W]}\bigl(\spec \mathcal{O}_{X,P}/\mathfrak{m}^n\bigr)^{[n]} \arrow[l, hook'] \arrow[r, hook]                                     & {T_{[\overline{I}]}G(n,\overline{A})}
\end{tikzcd}
\end{center}

\noindent Examining the proof of the tangent space description for Hilbert schemes\footnote{See the summary in Appendix \ref{subsectionhilberttangent} or directly \cite[Lemma 5.8]{geomalgcurv}.}, it is easy to see that the left horizontal arrow sends $f: \overline{I} \rightarrow A/I$ to itself as an $A$-module homomorphism, precomposed with $I \rightarrow \overline{I}$.\ The right horizontal arrow just sends $f$ to itself as a $k$-linear map. The right diagonal map is the differential of an inclusion of Grassmannians $G(l,\mathfrak{m}/\mathfrak{m}^2) \subseteq G(n,A/\mathfrak{m}^3)$, readily checked to send $g: I/\mathfrak{m}^2 \rightarrow \mathfrak{m}/I$ to the composition with $\overline{I} \rightarrow I/\mathfrak{m}^2$ and $\mathfrak{m}/I \hookrightarrow A/I$.\ Hence the middle vertical arrow and eventually the whole composition $T_{[I/\mathfrak{m}^2]}G \rightarrow T_{[W]}\hilb^{n}_d \rightarrow T_{[W]}X^{[n]}$ just maps any tangent vector $g$ to itself modulo some projections.\ This is now compared to (\ref{sequenceofinterest}):\par \vspace{3pt}
\noindent Observe that $\mathcal{E}xt^d_{\overline{q}}(\mathcal{O}_Z,\mathcal{I}_Z \otimes \omega_{\overline{q}}) \rightarrow \mathcal{Q}^{\vee} \otimes (\id,p)^{\ast}\mathcal{I}_Z \rightarrow \mathcal{Q}^{\vee} \otimes \mathcal{K}$ is just
\begin{align*}
\mathcal{E}xt^d_{\overline{q}}(\mathcal{O}_Z,\mathcal{I}_Z \otimes \omega_{\overline{q}}) \rightarrow \mathcal{E}xt^d_{\overline{q}}((\id,p)_{\ast}\mathcal{Q},\mathcal{I}_Z \otimes \omega_{\overline{q}}) \rightarrow \mathcal{E}xt^d_{\overline{q}}((\id,p)_{\ast}\mathcal{Q},(\id,p)_{\ast}\mathcal{K} \otimes \omega_{\overline{q}}),
\end{align*}
which pulled back to the $k$-point $\spec k \rightarrow \mathbb{G}$ corresponding to $(P,I) = [W]$ becomes
\begin{align*}
\ext^d_{\mathcal{O}_X}(\mathcal{O}_W,\mathcal{I}_W \otimes \omega_X) \longrightarrow \ext^d_{\mathcal{O}_X}(\mathcal{Q}(P),\mathcal{I}_W \otimes \omega_X) \longrightarrow \ext^d_{\mathcal{O}_X}(\mathcal{Q}(P),\mathcal{K}(P) \otimes \omega_X)
\end{align*}
with the natural morphisms in between.\ Dualizing (\ref{sequenceofinterest}), applying Serre duality on $X$ and replacing $X$ by an affine open afterwards therefore yields the desired sequence (\ref{restrictedhomseq}), so that (b) commutes.\ Notice that implicitly, naturality of base change for relative $\mathcal{E}xt$ sheaves (Lemma \ref{extbasechangeisnatural}) as well as the fact that the identifications 
\begin{align*}
    (\Omega_{X^{[n]}/k})_{|\mathbb{G}} \cong \iota^{\ast}\mathcal{E}xt^d_{\overline{\pi}}(\mathcal{O}_{\Xi},\mathcal{I}_{\Xi} \otimes \omega_{\overline{\overline{\pi}}}) \cong  \mathcal{E}xt^d_{\overline{q}}(\mathcal{O}_Z,\mathcal{I}_Z \otimes \omega_{\overline{q}}) \quad \text{and} \quad \Omega_p \cong \mathcal{Q}^{\vee} \otimes \mathcal{K}
\end{align*}
induce identity maps on tangent spaces up to scaling were used (see Prop.\ \ref{tangentsheafhilbert} and the naturality statement in \cite[Prop.\ F.212]{goertzwedhorn2}).\ These subtle technical statements are postponed until Appendix \ref{AppendixB}.\par \vspace{3pt}
\noindent Commutativity of the square (a) is also shown by going to tangent spaces:\ It has been implicitly used several times that for any morphism $f:V \rightarrow W$ of $k$-schemes such that $f(v) = w$, the tangent space at a point can be described functorially in equivalent ways, i.e.\ that there is a commutative diagram (setting $k[\overline{\varepsilon}] := k[\varepsilon]/(\varepsilon^2)$)
\begin{center}
\begin{tikzcd}
\Omega_{V/k}(v)^{\vee} \arrow[d, "df_{|v}"] \arrow[r, "\sim"] & {(\mathfrak{m}_{V,v}/\mathfrak{m}_{V,v}^2)^{\vee}} \arrow[d, "(f^{\sharp})^{\vee}"] & {\mor_k\bigl((\spec k[\overline{\varepsilon}],(\varepsilon)),(V,v)\bigr)} \arrow[d, "f\circ"] \arrow[l, "\sim"'] \\
\Omega_{W/k}(w)^{\vee} \arrow[r, "\sim"]                      & {(\mathfrak{m}_{W,w}/\mathfrak{m}_{W,w}^2)^{\vee}}                                  & {\mor_k\bigl((\spec k[\overline{\varepsilon}],(\varepsilon)),(W,w)\bigr).} \arrow[l, "\sim"']                    
\end{tikzcd}
\end{center}
\noindent The middle description of the tangent space is a very special case of the tangent space to a Hilbert scheme, namely to $X^{[1]} \cong X$: $\homo_{\mathcal{O}_{X,P}/\mathfrak{m}}(\mathfrak{m}/\mathfrak{m}^2,\mathcal{O}_{X,P}/\mathfrak{m}) \cong \homo_{\mathcal{O}_{X,P}}(\mathfrak{m},\mathcal{O}_{X,P}/\mathfrak{m})$.\par 
\noindent The pullback of (a) to a $k$-point $f_{P,I}: \spec k \rightarrow \mathbb{G}$ with image $[W] \in \mathbb{G}(k)$\footnote{$[W] = [(P,I)]$ corresponds as always to $P \in X$ with added scheme structure encoded in an ideal $I$.} is the dual of 
\begin{center}
\begin{tikzcd}
{\homo_{\mathcal{O}_{X,P}}(\mathfrak{m},\mathcal{O}_{X,P}/\mathfrak{m})} \arrow[d] & {\mor_k\bigl((\spec k[\overline{\varepsilon}],(\varepsilon)),(\mathbb{G},[W])\bigr)} \arrow[l, "dp"'] \arrow[d, "d\iota"] \\
{\homo_{\mathcal{O}_{X,P}}(I,\mathcal{O}_{X,P}/\mathfrak{m})}                      & {\homo_{\mathcal{O}_{X,P}}(I,\mathcal{O}_{X,P}/I).} \arrow[l]                                                
\end{tikzcd}
\end{center}
\noindent The two arrows pointing to the left bottom corner are the obvious ones, cf.\ Lemma \ref{extbasechangeisnatural}.\ Notice that $(\iota,p) = (\iota \times  \id) \circ (\id,p) : \mathbb{G} \hookrightarrow \mathbb{G} \times X \hookrightarrow X^{[n]} \times X$ factors through the universal family $\Xi_n$.\ The involved differentials $d\iota$ and $dp$ can be understood as follows.\par \vspace{3pt}
\noindent Let $t: \spec k[\overline{\varepsilon}] \rightarrow \mathbb{G}$ be a morphism (a tangent vector) restricting to $f_{P,I}$.\ Then $(\iota \circ t,p \circ t)$ is a tuple of tangent vectors to $[W] \in \mathbb{G}(k) \subseteq X^{[n]}(k)$, respectively to $p([W]) = P \in X^{[1]}$.\ By the above remark on $\Xi_n$, $(d\iota(t),dp(t)) = (\iota \circ t,p \circ t) \in T_{[(P,W)]}\Xi_n$, which is just a rephrasement of the containment $P \in W \subseteq X$.\par \vspace{3pt}
\noindent The question concerning (a) is whether the homomorphisms $d\iota(t): I \rightarrow \mathcal{O}_{X,P}/I$ and $dp(t): \mathfrak{m} \rightarrow \mathcal{O}_{X,P}/\mathfrak{m}$ fulfill $dp(t)_{|I} = d\iota(t) \modd \mathfrak{m}$.\ This holds true by the description of tangent spaces of the universal family $\Xi_n$ \cite[Lemma 8.8]{geomalgcurv}.\ Commutativity of (a) follows immediately, apart from the same subtleties mentioned at the end of (b).\par \vspace{3pt}
\noindent According to Proposition \ref{multiplicationmap} and its proof, there is a commutative diagram 
\begin{center}
\begin{tikzcd}
\mathcal{K} \arrow[rd, hook]                                          &                                           &             \\
{(\id,p)^{\ast}\mathcal{I}_Z} \arrow[u, "\beta", two heads] \arrow[r] & p^{\ast}\Omega_{X/k} \arrow[r, two heads] & \mathcal{Q},
\end{tikzcd}
\end{center}
\noindent hence the kernel of $(\id,p)^{\ast}\mathcal{I}_Z \rightarrow p^{\ast}\Omega_{X/k}$ coincides with $\ker \beta = S^2\mathcal{Q}$.\ The snake lemma applied to the whole diagram with squares (a) and (b) yields the exact sequence of locally free sheaves $S^2\mathcal{Q} \longrightarrow \mathcal{N}^{\vee} \longrightarrow \mathcal{Q}^{\vee} \otimes S^2\mathcal{Q} \longrightarrow \mathcal{Q} \longrightarrow 0$.
\end{proof}

\begin{kor} \label{moreexplicitsplitting}
    Let $d\geq 3$\footnote{These are the only cases of interest, cf.\ Conjecture (\ref{conjecture}) for $d \leq 2$.} and $\chara k = 0$.\ Then $\mathcal{N} \cong S^3 Q^{\vee} \otimes \det Q$.
\end{kor}

\begin{proof}
    Consider the decomposition (here $\chara k = 0$ becomes important, see Remark \ref{rules}) of  $\mathcal{Q} \otimes S^2\mathcal{Q}^{\vee}$ into $\Sigma^{2,-1}\mathcal{Q}^{\vee}$ and $\mathcal{Q}^{\vee}$ (Example \ref{pieriexample}) as well as the projection morphism
    \begin{align*}
        s: \mathcal{Q} \otimes S^2\mathcal{Q}^{\vee} \cong \Sigma^{2,-1}\mathcal{Q}^{\vee} \oplus \mathcal{Q}^{\vee} \longrightarrow \mathcal{Q}^{\vee}.
    \end{align*}
    See Appendix \ref{AppendixA} for the definition of Schur functors applied to locally free sheaves and how to perform computations with them.\ Composing the inclusion in (\ref{definingsequenceN}) with $s$ yields an endomorphism of $\mathcal{Q}^{\vee}$, but $\enndo_{\mathcal{O}_{\mathbb{G}}}(\mathcal{Q}^{\vee})$ is approximated by the Grothendieck spectral sequence with second page 
    \begin{align*}
        E_2^{i,j} = H^i\bigl(X,R^jp_{\ast}(\mathcal{Q}^{\vee} \otimes \mathcal{Q})\bigr).
    \end{align*}
    Here $p:\mathbb{G} \rightarrow X$ is the projection and $R^{j}p_{\ast}(\mathcal{Q}^{\vee} \otimes \mathcal{Q}) = 0$ for $j> 0$ by cohomology and base change since the objects $Q^{\vee}$ on the Grassmannian fibres are exceptional if $d \geq 3$ (Lemma \ref{exampleextcomputation}).\ Because $p_{\ast}(\mathcal{Q}^{\vee} \otimes \mathcal{Q}) \cong \mathcal{O}_X$, $\enndo_{\mathcal{O}_{\mathbb{G}}}(\mathcal{Q}^{\vee})$ is one-dimensional, so (\ref{definingsequenceN}) splits.\ Since $\mathbb{G}$ is a projective variety, $\Sigma^{2,-1}\mathcal{Q}^{\vee} \oplus \mathcal{Q}^{\vee} \cong \mathcal{N} \oplus \mathcal{Q}^{\vee}$ implies $\Sigma^{2,-1}\mathcal{Q}^{\vee} \cong \mathcal{N}$ according to \cite[Thm.\ 3]{atiyah}.\ It remains to observe that $\Sigma^{2,-1}\mathcal{Q}^{\vee} \cong S^3 \mathcal{Q}^{\vee} \otimes \det \mathcal{Q}$.
\end{proof}

\begin{bem} 
    Already \emph{after} formulating a proof of Theorem \ref{maintheoremnormalbundle}, the author noticed that it is easy to deduce from \cite[Rem.\ 2.5.10 \& Prop.\ 2.5.11]{goettsche1} that over $k = \mathbb{C}$ and for $d = \dime X \geq 2$, the projectivization $\mathbb{P}(\mathcal{N}^{\vee})$ is isomorphic to $\mathbb{P}(S^3Q)$.\ This shows that the normal bundle $\mathcal{N}$ of the planar locus is a twist of $S^3Q^{\vee}$ by a line bundle $\mathcal{L}$ (over $\mathbb{C}$).\ It is still not clear to the author how to obtain $\mathcal{L}$ from this faster without more effort.\ An attempt to use $\mathcal{N} \cong S^3 Q^{\vee} \otimes \mathcal{L}$ would be to compute its determinant and to apply the adjunction formula, but this leads at most to understanding $\mathcal{L}^{\otimes 4}$ better.\par 
    \noindent Since the interest lies in more general algebraically closed fields $k = \overline{k}$ anyway, it makes indeed sense to prove Theorem \ref{maintheoremnormalbundle} first and then deduce from (\ref{definingsequenceN}) that $\mathcal{N} \cong S^3 Q^{\vee} \otimes \det Q$ like it is done in Corollary \ref{moreexplicitsplitting}, i.e.\ that $\mathcal{L} = \det Q$.\ Then $\mathbb{P}(\mathcal{N}^{\vee}) \cong \mathbb{P}(S^3Q)$ follows (not only over $\mathbb{C}$), generalizing the result in \cite{goettsche1}.\ 
\end{bem}

\begin{bem}
    Apparently, the normal bundle sequence (\ref{definingsequenceN}) does not depend on the dimension $d = \dime X$.\ Intuitively, this is related to the fact that any length $3$ subscheme of $X$ is contained in a plane (assuming the existence of local, formal coordinates).\par \noindent  It seems likely that (\ref{definingsequenceN}) allows for a deeper interpretation from a deformation-theoretic point of view in the following sense: \par 
    \noindent The sheaf $\mathcal{Q} \otimes S^2\mathcal{Q}^{\vee} = \mathcal{H}om(S^2\mathcal{Q},\mathcal{Q})$ locally describes deformations of a planar subscheme to other length $3$ subschemes inside a plane containing it.\ This deformation space contains $\mathcal{Q}^{\vee} \cong \mathcal{H}om(\mathcal{Q},\mathcal{O}_{\mathbb{G}})$, corresponding intuitively to deformations of planar subschemes inside a plane that only change the underlying point of $X$.\par 
    \noindent It is an interesting question whether the \emph{global} sequence (\ref{definingsequenceN}) above can also be obtained the other way around, starting with the \emph{local} deformation-theoretic description.
\end{bem}

\section{The functors \texorpdfstring{$\Phi_{\alpha}$}{Phialpha}} \label{sectionphialpha}

The normal bundle computed in Theorem \ref{maintheoremnormalbundle} and Corollary \ref{moreexplicitsplitting} is the key to examine certain relative Fourier--Mukai transforms along the roof $X \overset{p}{\longleftarrow} \mathbb{G} \overset{\iota}{\longrightarrow} X^{[3]}$.

\begin{defi}
The functors $\db(X) \rightarrow \db(X^{[3]})$, $\mathcal{F}^{\bullet} \mapsto \iota_{\ast}(\Sigma^{\alpha}\mathcal{Q}^{\vee} \otimes p^{\ast}\mathcal{F}^{\bullet})$, are denoted by $\Phi_{\alpha}$.\ Here $\alpha = (\alpha_1,\alpha_2)$ is a partition of an integer $\leq 2\cdot (d-2)$, visualizable as a Young diagram with at most two rows and at most $d-2$ columns, i.e. $0 \leq \alpha_2 \leq \alpha_1 \leq d-2$. 
\end{defi}

\begin{bem}
Working with $\mathcal{Q}^{\vee}$ instead of $\mathcal{K}$ is way easier since $\mathcal{Q}^{\vee}$ only has rank $2$, hence Young diagrams with less rows occur.\ The upcoming calculations with Schur functors and Young diagrams all rely on Appendix \ref{AppendixA}.\ It is important to assume $\chara k = 0$ from now on to avoid representation-theoretic difficulties.
\end{bem}

\subsection{An auxiliary lemma} \par \vspace{0.2cm}

\noindent Let $G = G(2,d)$ be the Grassmannian of $2$-dimensional quotients (occuring as fibres of $p:\mathbb{G} \rightarrow X$) and let $Q$ be the tautological quotient bundle on $G$.\ The notation $\alpha \prec \beta$ is used for the total order on Young diagrams explained in Remark \ref{correctordering} of the appendix, in particular $|\beta| < |\alpha| \Rightarrow \alpha \prec \beta$ and $\alpha \prec \beta \Rightarrow |\beta| \leq |\alpha|$.

\begin{lemma} \label{auxiliary}
    Let $\alpha \prec \beta$ or $\alpha = \beta$.\ Then the cohomology $H^{\ast}(G, \wedge^q \mathcal{N}' \otimes \Sigma^{-\beta}Q^{\vee} \otimes\Sigma^{\alpha}Q^{\vee})$ vanishes in every degree for all integers $q > 0$ provided that $\alpha_1 - \beta_2 \leq d-5$.\ Here $\mathcal{N}'$ is the restriction of the normal bundle of $\iota: \mathbb{G} \hookrightarrow X^{[3]}$ to a fibre $G$ of $p:\mathbb{G} \rightarrow X$.
\end{lemma}

\noindent Let $\lambda_{\alpha} = \alpha_1 - \alpha_2$ and $\lambda_{\beta} = \beta_1 - \beta_2$.\ Lemma \ref{auxiliary} will make use of the identity
\begin{align} \label{extcomputation}
H^{\ast}\bigl(G, \wedge^q \mathcal{N}' \otimes \Sigma^{-\beta}Q^{\vee} \otimes\Sigma^{\alpha}Q^{\vee}\bigr) = \bigoplus_{\gamma = 0}^{\min\lbrace \lambda_{\alpha},\lambda_{\beta} \rbrace} H^i\bigl(G,\wedge^q \mathcal{N}' \otimes \Sigma^{\alpha_1 - \beta_2 - \gamma, \alpha_2 - \beta_1 + \gamma}Q^{\vee}\bigr),
\end{align}
see the proof of Lemma \ref{exampleextcomputation}.\ The abbreviation $\Sigma^{a,b}Q^{\vee} = \Sigma^{\alpha_1 - \beta_2 - \gamma,\alpha_2 - \beta_1 + \gamma}Q^{\vee}$ is used and Lemma \ref{exampleextcomputation} also shows that $\alpha \prec \beta \Rightarrow 1 \leq a \leq d-2$.\ By assumption, even $a \leq d-5$ holds.\ Observe that $\alpha \prec \beta$ and $\alpha = \beta$ both imply that $d \geq 5$, so Corollary \ref{moreexplicitsplitting} applies.

\begin{proof}[Proof of Lemma \ref{auxiliary}]
    Recall that on each fibre of $p$, there is a short exact sequence\footnote{For the calculations where $q \in \lbrace 2,4 \rbrace$, the short exact sequence is just as useful as Corollary \ref{moreexplicitsplitting}.}
\begin{align} \label{restricteddefiningsequenceN}
    0 \longrightarrow Q^{\vee} \longrightarrow Q \otimes S^2 Q^{\vee} \longrightarrow \mathcal{N}' \longrightarrow 0 \quad \quad \Longrightarrow \quad \mathcal{N}' \cong S^3Q^{\vee} \otimes \det Q
\end{align}
obtained as the restriction of (\ref{definingsequenceN}).\ Let $q = 1$.\ Then $\mathcal{N}' \otimes \Sigma^{a,b}Q^{\vee}$ is isomorphic to 
\begin{align*}
& ~S^3Q^{\vee} \otimes \Sigma^{(a-b,0)}Q^{\vee} \otimes \det Q^{\vee \otimes (b-1)} \\
\cong & ~\Sigma^{(a + 2 ,b-1)}Q^{\vee} \oplus \Sigma^{(a+1 ,b)}Q^{\vee} \oplus \Sigma^{(a,b+1)}Q^{\vee} \oplus \Sigma^{(a-1,b+2)}Q^{\vee}, \quad a \leq d-5,
\end{align*}
according to Pieri's formula \ref{pieriexample}.\ Here $\Sigma^{(a+1 ,b)}Q^{\vee}$ only shows up if $a-b \geq 1$, $\Sigma^{(a,b+1)}Q^{\vee}$ only shows up if $a-b \geq 2$ and $\Sigma^{(a-1,b+2)}Q^{\vee}$ only shows up if $a-b \geq 3$.\ This makes sense for $\alpha \prec \beta$ and also for $\alpha = \beta$, where in the latter case the summands in the decomposition of $\mathcal{N}' \otimes \Sigma^{k,-k}Q^{\vee}$ only show up if $k > 0$ (respectively $k> 1$).\ It has to be excluded that $a - 1 = 0$ in the case where $\alpha_1 - \beta_2 - \alpha_2 + \beta_1 - 2\gamma =  a-b  \geq 3$: \par
\noindent Note that $\alpha_2 - \beta_1 \geq \beta_2 - \alpha_1$ because $\alpha \prec \beta \Rightarrow |\alpha| \geq |\beta|$ or $\alpha = \beta$.\ Therefore
\begin{align*}
a = \alpha_1 - \beta_2 - \gamma &\geq  3 + \alpha_2 - \beta_1 + \gamma = 3 + b \\
&\geq 3 + \beta_2 - \alpha_1 + \gamma = 3- a,
\end{align*}
hence $2a \geq 3$, implying $a > 1$.\ Corollary \ref{vanishinga1} yields $H^{\ast}(G,\mathcal{N}' \otimes \Sigma^{a,b}Q^{\vee}) = 0$ if $\alpha \prec \beta$ and also if $\alpha = \beta$ (where $a = k = -b$).\par
\noindent The next easy case is $q=4$. Note that $\wedge^{4}\mathcal{N}' \cong (\det Q^{\vee})^{\otimes 2}$ by Lemma \ref{computationwedge4}.\ Consequently, $\wedge^{4}\mathcal{N}' \otimes \Sigma^{a,b}Q^{\vee} \cong \Sigma^{a+2,b+2}Q^{\vee}$.\ Since $1 \leq a \leq d-5$ if $\alpha \prec \beta$, obviously $1 \leq a+2 \leq d-2$.\ Therefore $H^{\ast}(G,\wedge^{4}\mathcal{N}' \otimes \Sigma^{a,b}Q^{\vee}) = 0$ by Corollary \ref{vanishinga1}.\ This also works if $\alpha = \beta$ and $\wedge^{4}\mathcal{N}' \otimes \Sigma^{k,-k}Q^{\vee} \cong \Sigma^{k+2,-k+2}Q^{\vee}$ for $0 \leq a = k \leq d-5$.\par 
\noindent The case $q = 2$ is where the assumption $\alpha_1 - \beta_2 \leq d-5$ will really be exploited.\ There is a filtration of second exterior powers for (\ref{restricteddefiningsequenceN}) \cite[Exc. II.5.16]{hartshorne} that reads as follows after tensoring with $\Sigma^{a,b}Q^{\vee}$ (or $\Sigma^{k,-k}Q^{\vee}$ if $\alpha = \beta$).\ 
\begin{gather*}
0 \longrightarrow F^1 \otimes \Sigma^{a,b}Q^{\vee} \longrightarrow \wedge^2(S^2Q^{\vee} \otimes Q) \otimes \Sigma^{a,b}Q^{\vee} \longrightarrow \wedge^2\mathcal{N}' \otimes \Sigma^{a,b}Q^{\vee} \longrightarrow 0 \\
0 \longrightarrow \det Q^{\vee} \otimes \Sigma^{a,b}Q^{\vee} \longrightarrow F^1 \otimes \Sigma^{a,b}Q^{\vee} \longrightarrow Q^{\vee} \otimes \mathcal{N}' \otimes \Sigma^{a,b}Q^{\vee} \longrightarrow 0
\end{gather*}
\noindent The term $\wedge^2(S^2Q^{\vee} \otimes Q)$ decomposes as
\begin{align} \label{wedgedecomposition}
    \wedge^2(S^2Q^{\vee} \otimes Q) \cong (\Sigma^{3,-1}Q^{\vee})^{\oplus 2} \oplus (\Sigma^{1,1}Q^{\vee})^{\oplus 2} \oplus S^2Q^{\vee}
\end{align}
according to Lemma \ref{computationwedge2}.\ Step by step, it is now shown that all cohomology of the objects in the filtration vanishes.\ Consider the second short exact sequence first.
\begin{itemize}
    \item Since $\det Q^{\vee} \otimes \Sigma^{a,b}Q^{\vee} \cong \Sigma^{a+1,b+1}Q^{\vee}$ and $1 \leq a+1 \leq d-2$, this locally free sheaf has no cohomology at all on $G$ (Corollary \ref{vanishinga1}).
    \item Next comes the cohomology of $Q^{\vee} \otimes \mathcal{N}' \otimes \Sigma^{a,b}Q^{\vee}$, which is computed using the decomposition of $\mathcal{N}' \otimes \Sigma^{a,b}Q^{\vee}$ derived in the previous step $q=1$, i.e.
    \begin{align*}
        Q^{\vee} \otimes \mathcal{N}' \otimes \Sigma^{a,b}Q^{\vee} \cong Q^{\vee} \otimes \bigl(\Sigma^{(a + 2 ,b-1)}Q^{\vee} \oplus \Sigma^{(a+1 ,b)}Q^{\vee} \oplus \Sigma^{(a,b+1)}Q^{\vee} \oplus \Sigma^{(a-1,b+2)}Q^{\vee}\bigr).
    \end{align*}
    Tensoring with $Q^{\vee}$ yields $Q^{\vee} \otimes \Sigma^{m,n}Q^{\vee} \cong \Sigma^{m+1,n}Q^{\vee} \oplus \Sigma^{m,n+1}Q^{\vee}$ in general, so that the only distinct summands in the decomposition of $Q^{\vee} \otimes \mathcal{N}' \otimes \Sigma^{a,b}Q^{\vee}$ are
    \begin{align} \label{listofsummands}
        \Sigma^{a+3,b-1}Q^{\vee}, \quad \Sigma^{a+2,b}Q^{\vee}, \quad \Sigma^{a+1,b+1}Q^{\vee}, \quad \Sigma^{a,b+2}Q^{\vee} \quad \text{and} \quad  \Sigma^{a-1,b+3}Q^{\vee}.
    \end{align}
    Some summands may not show up if $a-b$ is too small, in particular $\Sigma^{a-1,b+3}Q^{\vee}$ only occurs for $a-b \geq 4$.\ In the special case $\alpha = \beta$ where $\Sigma^{a,b}Q^{\vee} = \Sigma^{k,-k}Q^{\vee}$, this means that the last two (one) are omitted if $k = 0$ ($k=1$).\ As in the proof of $q = 1$, $a-b \geq 3$ leads to $a > 1$.\ Now Corollary \ref{vanishinga1} can be applied once more to finish the argument and to conclude that all in all, $H^{\ast}(G,Q^{\vee} \otimes \mathcal{N}' \otimes \Sigma^{a,b}Q^{\vee}) = 0$.
\end{itemize}
\noindent It follows that in the filtration above, $F^1 \otimes \Sigma^{a,b}Q^{\vee}$ has no cohomology on $G$.
\begin{itemize}    
    \item  What remains to show is the vanishing of $H^{\ast}(G,\wedge^2(S^2Q^{\vee} \otimes Q) \otimes \Sigma^{a,b}Q^{\vee})$.\ According to (\ref{wedgedecomposition}), this leads to showing 
    \begin{align*}
        H^{\ast}(G,\Sigma^{3,-1}Q^{\vee} \otimes \Sigma^{a,b}Q^{\vee}) = H^{\ast}(G,\det Q^{\vee} \otimes \Sigma^{a,b}Q^{\vee}) = H^{\ast}(G,\Sigma^{2,0}Q^{\vee} \otimes \Sigma^{a,b}Q^{\vee}) = 0.
    \end{align*}
    The middle cohomology has been dealt with before.\ Observe that $\Sigma^{2,0}Q^{\vee} \otimes \Sigma^{a,b}Q^{\vee}$ is isomorphic to $\Sigma^{a+2,b}Q^{\vee} \oplus \Sigma^{a+1,b+1}Q^{\vee} \oplus \Sigma^{a,b+2}Q^{\vee}$, the last two summands only showing up if $a-b \geq 1$, respectively $a-b \geq 2$.\ If $\alpha = \beta$ and $\Sigma^{a,b}Q^{\vee} = \Sigma^{k,-k}Q^{\vee}$, this means that only the first summand shows up if $k=0$.\ These objects have no cohomology on $G$ by Corollary \ref{vanishinga1}.\ Finally, it is an easy consequence of Pieri's formula that the summands in the decomposition of $\Sigma^{3,-1}Q^{\vee} \otimes \Sigma^{a,b}Q^{\vee}$ are the same as in (\ref{listofsummands}).\ Using $\alpha_1 - \beta_2 \leq d-5$, their cohomologies have already been shown to vanish with the aid of Corollary \ref{vanishinga1}.
\end{itemize}
\noindent Finally, let $q = 3$.\ Notice that $\wedge^3\mathcal{N}' \cong \det \mathcal{N}' \otimes (\mathcal{N}')^{ \vee}$, cf.\ \cite[Exc. II.5.16(b)]{hartshorne}.\ Using $\det \mathcal{N}' \cong (\det Q^{\vee})^{\otimes 2}$, $(\mathcal{N}')^{ \vee} \cong S^3Q \otimes \det Q^{\vee}$ as well as the rules \ref{pieriexample} and \ref{dualformula}, it follows that $\wedge^3\mathcal{N}' \cong S^3Q^{\vee}$.\ What remains to show is $H^{\ast}\bigl(G,S^3Q^{\vee} \otimes S^{a-b}Q^{\vee} \otimes   (\det Q^{\vee})^{\otimes b}\bigr) = 0$.\ This can be done by applying Pieri's formula.\ Firstly, this means observing
\begin{align*}
    S^3Q^{\vee} \otimes S^{a-b}Q^{\vee} \otimes (\det Q^{\vee})^{\otimes b} \cong \Sigma^{a + 3,b}Q^{\vee} \oplus \Sigma^{a +2,b+1}Q^{\vee} \oplus \Sigma^{a + 1,b+2}Q^{\vee} \oplus \Sigma^{a,b+3}Q^{\vee},
\end{align*}
with the second (third, fourth) summand occuring only if $a-b \geq 1$ ($\geq 2, \geq 3$) and $\alpha \prec  \beta$.\ For $\alpha = \beta$ and $\Sigma^{a,b}Q^{\vee} = \Sigma^{k,-k}Q^{\vee}$, only the first summand shows up if $k=0$ and the last one is missing if $k = 1$.\ Corollary \ref{vanishinga1} applies.\ All cases $q = 1,2,3,4$ are finished.
\end{proof}

\subsection{Fully faithfulness}

\noindent 

\begin{thrm}[Thm.\ \ref{fullyfaithfulthm}] \label{fullyfaithfultheorem1}
Let $d \geq 5$ and let $\alpha$ be a partition inscribed into a rectangle with at most $2$ rows and $d-2$ columns.\ Suppose $\lambda_{\alpha} = \alpha_1 - \alpha_2 \leq d-5$.\ Then the assumptions of Proposition \ref{appliedfibrecriterion}(1) are fulfilled and $\Phi_{\alpha}: \db(X) \rightarrow \db(X^{[3]})$ is fully faithful.
\end{thrm}

\begin{proof}
Let $\mathcal{N}'$ be the restriction of $\mathcal{N}_{\mathbb{G}/X^{[3]}}$ to a fibre of $p:\mathbb{G} \rightarrow X$.\ Before applying Proposition \ref{appliedfibrecriterion}, note that the kernel $\Sigma^{\alpha}\mathcal{Q}^{\vee}$ of the Fourier--Mukai transform $\Phi_{\alpha}$ restricts to $\Sigma^{\alpha}Q^{\vee}$ on the fibres of $p$ and that $\Sigma^{\alpha}Q \otimes \Sigma^{\alpha}Q^{\vee}$ is naturally isomorphic to $\bigoplus_{k=0}^{\lambda_{\alpha}} \Sigma^{k,- k}Q^{\vee}$ (Lemma \ref{exampleextcomputation}).\ The goal is to show that on the fibres $G = p^{-1}(P)$, 
\begin{align*}
    \homo(\Sigma^{\alpha}Q^{\vee},\Sigma^{\alpha}Q^{\vee}) = k \quad \text{and} \quad H^p(G, \wedge^q \mathcal{N}' \otimes \Sigma^{\alpha}Q \otimes \Sigma^{\alpha}Q^{\vee}) = 0 \quad \text{ for } p + q > 0.
\end{align*}
This leads to showing that the vector space $H^0(G,\bigoplus \Sigma^{k,-k}Q^{\vee})$ is one-dimensional and that $H^p(G, \wedge^q \mathcal{N}' \otimes \Sigma^{k,-k}Q^{\vee}) = 0$ for $p+q > 0$ and $0 \leq k \leq \lambda_{\alpha} \leq d-5$.\par
\noindent The one-dimensionality as well as the case $q = 0$ follow from exceptionality of the bundles $\Sigma^{\alpha}Q^{\vee}$ proven in Lemma \ref{exampleextcomputation}.\ The rest is a direct consequence of Lemma \ref{auxiliary}.
\end{proof}

\begin{bem}
The author suspects Theorem \ref{fullyfaithfultheorem1} to remain true in the case of a smooth \emph{quasi}-projective variety $X$,\footnote{Potentially using a smooth compactification and base change arguments for FM transforms.} though the above proof needs to be modified if $X$ is not projective since it relied implicitly on the Bondal--Orlov-criterion \ref{bondalorlovff}: \par 
\noindent Theorem \ref{bondalorlovff} uses that skyscraper sheaves form a spanning class in $\db(X)$ (Remark \ref{remarkspanningclasses}), which remains true in the quasi-projective regular case.\ But this is only useful for testing fully-faithfulness if the functor in question has \emph{both adjoints} \cite[Tag 0G25]{stacks-project}, guaranteed for example when working with projective varieties.\ 
\end{bem}

\subsection{Orthogonality}

\noindent It is time to relate different functors $\Phi_{\alpha}, \Phi_{\beta}: \db(X) \rightarrow \db(X^{[3]})$ to each other.\ The usual notation $\alpha \prec \beta$ from Remark \ref{correctordering} is used for the total order on the Young diagrams of interest.

\begin{thrm}[Thm.\ \ref{shortversioncollection}] \label{semiorthogonalsequencethrm}
For two partitions $\alpha \prec \beta$ of an integer $\leq 2(d-2)$ such that $\alpha_1 - \beta_2 \leq d-5$ holds, $\homo_{\db(X^{[3]})}\left(\Phi_{\beta}(\db(X)),\Phi_{\alpha}(\db(X))\right) = 0$.\ In particular, considering all fully faithful functors $\Phi_{\alpha}$ from Theorem \ref{fullyfaithfultheorem1} with the additional restriction that $\alpha_2 \geq 3$ yields a collection of $\binom{d-3}{2}$ semi-orthogonal subcategories of $\db(X^{[3]})$.
\end{thrm}

\noindent Recall in advance that on the absolute Grassmannian $G = G(2,d)$ with tautological quotient bundle $Q$, $\ext^{\ast}(\Sigma^{\beta}Q^{\vee},\Sigma^{\alpha}Q^{\vee}) = 0$ for partitions $\alpha \prec \beta$.\ This is Example \ref{exampleextcomputation}, a special case of Kapranov's famous result \ref{kapranov85main}. 

\begin{proof}[Proof of Theorem \ref{semiorthogonalsequencethrm}]
The proof is very similar to the one of Theorem \ref{fullyfaithfultheorem1}.\ According to Proposition \ref{appliedfibrecriterion}(2), the first thing to check is that $\homo(\Sigma^{\beta}Q^{\vee},\Sigma^{\alpha}Q^{\vee})= 0$ on each fibre $G$ of $p:\mathbb{G} \rightarrow X$.\ This is part of Lemma \ref{exampleextcomputation}, i.e.\ exceptionality of the sequence $\lbrace \Sigma^{\alpha}Q^{\vee}\rbrace_{\alpha}$.\ Let again $\mathcal{N}'$ be the restriction of $ \mathcal{N}_{\mathbb{G}/X^{[3]}}$ to $G$.\ The next thing to check is
\begin{align*}
H^p\bigl(G, \wedge^q \mathcal{N}' \otimes \Sigma^{-\beta}Q^{\vee} \otimes\Sigma^{\alpha}Q^{\vee}\bigr) =: \bigoplus H^p\bigl(G, \wedge^q \mathcal{N}' \otimes \Sigma^{a,b}Q^{\vee}\bigr) = 0 \quad \text{ for } p + q > 0,
\end{align*}
where $\Sigma^{a,b}Q^{\vee} = \Sigma^{\alpha_1 - \beta_2 - \gamma,\alpha_2 - \beta_1 + \gamma}Q^{\vee}$, cf.\ (\ref{extcomputation}).\ This is proven in Lemma \ref{auxiliary} for $q > 0$.\ If $q = 0$, the rest of the statement follows from $\ext^{\ast}(\Sigma^{\beta}Q^{\vee},\Sigma^{\alpha}Q^{\vee}) = 0$.
\end{proof}

\section{Application to generalized Kummer Varieties}

\noindent Recall the definition of relative Fourier--Mukai transforms over a base scheme $S$ given in Section \ref{subsectionderivedcats}.\ With Theorem \ref{fullyfaithfultheorem1} in mind, the question arises whether restricting (or more generally base-changing) a fully faithful FM transform yields again a fully faithful FM transform between interesting smooth projective varieties.\ More precisely:\par 
\noindent Starting with a roof $X \overset{p}{\longleftarrow} V \overset{q}{\longrightarrow} Y$ between projective varieties and a  relative FM transform $\Phi: \db(X) \rightarrow \db(Y)$ with kernel $\mathcal{K}^{\bullet} \in \db(V)$, suppose that there exists a target variety $S$ with sufficiently well-behaved morphisms $X \rightarrow S$ and $Y \rightarrow S$ so that $(V \rightarrow X \rightarrow S) = (V \rightarrow Y \rightarrow S)$, inducing a proper morphism $\iota: V \rightarrow X \times_S Y$.\par 
\noindent If it can be shown that $\Phi = \Phi_{R\iota_{\ast}\mathcal{K}^{\bullet}}$ being fully faithful implies the same for a base change $\Phi_T$, it will follow that there is an embedding $\db(X \times_S T) \hookrightarrow \db(Y \times_S T)$.\ Optimally, this procedure also preserves semi-orthogonality of the images of several FM transforms.

\begin{defi}[{\cite[Sect.\ 2.8]{homprojdual}}] \label{defibasechangeFMtrafo}
    Let $\Phi = Rq_{\ast}\bigl(\mathcal{K}^{\bullet} \otimes^L Lp^{\ast}(-)\bigr): \db(X) \rightarrow \db(Y)$ with kernel  $\mathcal{K}^{\bullet} \in \db(X \times_S Y)$ be a relative FM transform over a base $S$.\ Fix a base change morphism $\phi: T \rightarrow S$.\ The \emph{base change of $\Phi$ along $\phi$} is defined via the following cartesian diagram as $\Phi_T = R\overline{q}_{\ast}(Lj^{\ast}\mathcal{K}^{\bullet} \otimes^L \overline{p}^{\ast}(-)): \derived(X_T) \rightarrow \derived(Y_T)$.\par 
    
\begin{center} 
\begin{tikzcd}
X_T \times_T Y_T \arrow[rr, "\overline{q}"' ] \arrow[dd, "\overline{p}" ] \arrow[rd, "j" description] &                                               & Y_T \arrow[dd] \arrow[rrd]    &  &              \\
                                                                                                  & X \times_S Y \arrow[dd, "p" description] \arrow[rrr, crossing over] &                               &  & Y \arrow[dd] \\
X_T \arrow[rr] \arrow[rd]                                                                         &                                               & T \arrow[rrd, "\phi"] &  &              \\
                                                                                                  & X \arrow[from=uu, crossing over] \arrow[rrr]                                 &                               &  & S         
\end{tikzcd}
\end{center}
\end{defi}

\begin{bem} \label{remcartesiancube}
    The notations $X_T = X \times_S T$ and $Y_T = Y \times_S T$ are used.\ If $T = \spec k$ corresponds to a $k$-point $s \in S(k)$, this construction leads to the \textit{restriction} of $\Phi$ to $s$, denoted by $\Phi_s: \derived(X_s) \rightarrow \derived(Y_s)$ and also defined in \cite[Sect.\ 1.3]{symmquotstacks}.\ In Definition \ref{defibasechangeFMtrafo}, it is not claimed that the base change of $\Phi$ is well-defined on the level of bounded derived categories (unless additional restrictions on the involved morphisms are imposed).\
\end{bem}

\noindent Let $\chara k = 0$.\ Recall that in the setting of Section \ref{sectionphialpha}, $\Phi_{\alpha}$ is the FM transform along
\begin{center}
    \begin{tikzcd}
{X^{[3]}_{(3)}} \arrow[d] \arrow[r, hook] \arrow[rd, "\square", phantom] & {X^{[3]}} \arrow[d, "\rho"] \\
X \arrow[r, "\overline\Delta_{123}", hook]                               & X^{(3)}                    
\end{tikzcd}
\end{center}
with kernel given by the pushforward of $\Sigma^{\alpha}\mathcal{Q}^{\vee}$ under the closed embedding of  $\mathbb{G}(2,\Omega_X)$ into $ X^{[3]}_{(3),\red} \subseteq X^{[3]}_{(3)}$, cf.\ Proposition \ref{iclosedemb}.\ The morphisms $\overline\Delta_{123}$ and $\rho$ are certainly not smooth in general, e.g.\ the Hilbert--Chow morphism usually has non-smooth fibres isomorphic to $\hilb^n_d$ on the one hand as well as zero-dimensional fibres.\par 
\noindent Nonetheless, In the case of an abelian variety $X = A$, there is an $S_3$-invariant summation morphism $\Sigma: A^3 \rightarrow A$, descending to $\overline{\Sigma}: A^{(3)} \rightarrow A$.\ Together with the multiplication $\mu_3 = \overline{\Sigma} \circ \overline{\Delta_{123}}:A \rightarrow A$ by $3$, this yields a commutative diagram 
\begin{center}
    \begin{tikzcd}
{A^{[3]}_{(3)}} \arrow[d] \arrow[r] & {A^{[3]}} \arrow[d, "\overline{\Sigma} \circ \rho"] \\
A \arrow[r, "{\mu_3}"]                & A.                     
\end{tikzcd}
\end{center}
Now $\mathbb{G}(2,\Omega_A) \cong A \times G(2,d) \hookrightarrow A^{[3]}_{(3)} \rightarrow A \times_{\mu_3,A} A^{[3]}$ remains\footnote{Using e.g.\ the characterization \cite[Cor.\ 12.92]{goertzwedhorn} of closed embeddings.} a closed embedding since 
\begin{align*}
    \mathbb{G}(2,\Omega_A) \rightarrow A \times_{\mu_3,A} A^{[3]} \hookrightarrow A \times A^{[3]}
\end{align*}
is a closed embedding as observed before.\ Since $\chara k = 0$ in this section, $\mu_3$ is smooth \cite[\href{https://stacks.math.columbia.edu/tag/0BFH}{Tag 0BFH}]{stacks-project}.\ The fibre of $\mu_3:A \rightarrow A$ over $0 \in A$ is a collection of $3^{2d}$ points since $k = \overline{k}$ as well \cite[\href{https://stacks.math.columbia.edu/tag/03RP}{Tag 03RP}]{stacks-project}.\ The fibre of $\overline{\Sigma} \circ \rho: A^{[3]} \rightarrow A$ is part of Definition \ref{defikummer}:

\begin{defi} \label{defikummer}
    Let $A$ be an abelian variety (not necessarily a surface) over an algebraically closed field of characteristic zero.\ The \emph{generalized Kummer variety} $K_n(A)$ is defined by the following fibred product. 
    \begin{center}
\begin{tikzcd}
K_n(A) \arrow[d] \arrow[r] \arrow[rd, "\square", phantom] & {A^{[n]}_{\red}} \arrow[d, "\overline{\Sigma} \circ \rho"] \\
\lbrace 0 \rbrace  \arrow[r]                                 & A.                                                        
\end{tikzcd}
    \end{center}
\end{defi}

\begin{bem}
    Defined this way, $K_n(A)$ is smooth for $n \leq 3$ since then $A^{[n]}_{\red} = A^{[n]}$ and $\overline{\Sigma} \circ \rho$ is generically smooth.\ Equipping $A^{[n]}$ with a group action coming from addition on $A$ that makes $\overline{\Sigma} \circ \rho$ equivariant, smoothness of $\overline{\Sigma} \circ \rho$ everywhere follows.
\end{bem}

\noindent Now that only smooth morphisms are involved, it becomes more plausible that the steps alluded at the beginning of this section are actually feasible.

\begin{prop} \label{basechangefullyfaithful}
    Suppose that there is a cartesian diagram of projective $k$-schemes
    \begin{center}
    \begin{tikzcd}
X \times_S Y \arrow[d, "p"'] \arrow[r, "q"] & Y \arrow[d, "\pi_Y"] \\
X \arrow[r, "\pi_X"]                        & S .                 
    \end{tikzcd}
    \end{center}
    Suppose furthermore that $\pi_X$ and $\pi_Y$ are smooth morphisms and that $S$ is smooth.\
    \begin{enumerate}
        \item[(1)]  If $\Phi: \db(X) \rightarrow \db(Y)$ is a fully faithful FM transform with bounded kernel $\mathcal{K}^{\bullet}$ as above, then the restrictions 
\begin{align*}
    \Phi_s: \db(X_s) \rightarrow \db(Y_s)
\end{align*}
    for all $s \in S(k)$ are fully faithful, cf.\ Definition \ref{defibasechangeFMtrafo} and Remark \ref{remcartesiancube}.\ 
    \item[(2)] Suppose that there is a semi-orthogonal sequence arising from fully faithful functors $\Phi_i: \db(X_i) \rightarrow \db(Y)$ with bounded kernels $\mathcal{K}_i^{\bullet}$ for $i = 1, \hdots, n$, and all pairs $(X_i,Y)$ and projections obey the assumptions from above.\ Then this induces semi-orthogonal sequences $\left(\Phi_i(\db(X_{i,s}))\right)_i$ of $ \db(Y_s)$ for all $s \in S(k)$. 
    \end{enumerate}
    The statements remain valid for the restrictions to open subsets $\Phi_U: \db(X_U) \rightarrow \db(Y_U)$.
\end{prop}

\noindent A way to prove Proposition \ref{basechangefullyfaithful} is to consider adjoints of FM transforms.\ This has been done in \cite[Prop. 2.44]{kuznetsovhyperplanes} and many of the results are summarized in \cite{homprojdual}.

\begin{proof}[{Proof of Proposition \ref{basechangefullyfaithful}}]
    (1) According to \cite[Prop.\ 2.39]{homprojdual}, whether a base change $\Phi_T$ along a morphism $\phi:T \rightarrow S$ is fully faithful boils down to the question whether $\phi$ is \emph{faithful for the pair} $(X,Y)$, which can be checked as follows \cite[2.35 - 2.36]{homprojdual}.\par 
    \noindent The morphism $\phi$ has to be faithful with respect to $\pi_X$, $\pi_Y$ and $\pi_X \circ p = \pi_Y \circ q$.\ This holds if $\phi$ is flat, proving already the case $T = U \subseteq S$.\ For the case of the closed embedding $T = \lbrace s \rbrace \hookrightarrow S$, it is enough to observe \cite[Lemma 2.36]{homprojdual} that $T$ is smooth, $X,Y,S$ are smooth, $\pi_X$ and $\pi_Y$ are smooth by assumption, $\pi_X \circ p = \pi_Y \circ q$ is smooth as a composition of smooth morphisms, and thus $X \times_S Y$ is also smooth.\ In particular, all involved fibre products are of the expected dimension.\ The claim follows. \par 
\noindent (2) There is a result on base change of \emph{full} semi-orthogonal sequences \cite[Thm\ 2.40]{homprojdual}, i.e.\ semi-orthogonal decompositions, that does not apply here (see however Lemma \ref{complement}).\ The necessary technical statement whose assumptions are fulfilled with the same arguments as in (1) can be found in \cite[Prop.\ 2.44(i)]{kuznetsovhyperplanes}.
\end{proof}

\begin{prop}[Thm.\ \ref{shortversionkummer}] \label{longversionkummer}
Let $A$ be an abelian variety of dimension $d$ over an algebraically closed field of characteristic zero.\ Every fully faithful Fourier--Mukai transform $\Phi: \db(A) \rightarrow \db(A^{[3]})$ from Theorem \ref{fullyfaithfultheorem1} induces $3^{2d}$ exceptional objects in $\db(K_3(A))$.\  In particular, Theorem \ref{semiorthogonalsequencethrm} yields an exceptional sequence of length $\binom{d-3}{2}\cdot 3^{2d}$ in $\db(K_3(A))$.
\end{prop}

\begin{proof}
This follows directly from Proposition \ref{basechangefullyfaithful} and from the considerations at the beginning of this section.
\end{proof}

\begin{bem}
    Notice that the same reasoning can be applied to the semi-orthogonal decomposition \cite[Thm.\ 4.1(ii)]{krugploogsosna} of $A^{[2]}$ to obtain exceptional objects inside $\db(K_2(A))$.\ A more direct approach to obtain a semi-orthogonal decomposition of $\db(K_2(A))$ is formulated in \cite[Cor.\ 6.1]{krugploogsosna}.
\end{bem}

\appendix
\section{Cohomology on Grassmannians}  \label{AppendixA}

Let $G = G(k,V)$ be the Grassmannian of $k$-dimensional quotients of a vector space $V$.\ In this appendix, the necessary material needed to understand the construction of Kapranov's strong, full exceptional sequence \cite{kapranov1} in the derived category $\db(G)$ is briefly explained and summarized without proofs.\ This includes parts of the representation theory of the general linear group $\GL_d(k)$, calculations with Schur functors $\Sigma^{\alpha}$ and the application of the Borel--Weil--Bott theorem for $\GL_d(k)$ to Grassmannians.

\begin{thrm}[{\cite{akinbuchsbaumweyman}, \cite[Thm.\ 2.2.9-2.2.10]{weyman}}] 
    Let $\chara k = 0$ and $V = k^d$. Then the group $\GL_d(k)$ is linearly reductive, and any irreducible rational representation is of the form $\Sigma^{\alpha}V \otimes (\det V)^{\otimes m}$.\ Here $m\in \mathbb{Z}$ and $\alpha$ is a partition consisting of at most $d-1$ non-negative integers.\ This becomes false if $k$ is an infinite field of \emph{positive} characteristic.
\end{thrm}

\noindent There are different ways to define the Schur functors $\Sigma^{\alpha}$, and $V = k^d$ above can be replaced by a finite free module over a commutative ring or even by more general modules.\ The following overview includes the most common definitions.\ Assume first that $\alpha$ contains only non-negative integers $\alpha_i \geq 0$.\ Write $|\alpha| = \sum \alpha_i$ and $\alpha'$ for the transpose partition.

\begin{itemize}
    \item In \cite[Ch. 4 \& 6]{fultonharrisrep}, the representation theory of the symmetric group $S_m$ over $k=\mathbb{C}$ is explained first before defining Schur functors as images of symmetrizing operators.\ More precisely, $S_m$ acts on $V^{\otimes m}$ from the right and the Young symmetrizer $c_{\alpha} \in k[S_m]$ defines the $\GL_d(k)$-representation $\Sigma^{\alpha}V := V^{\otimes m}c_{\alpha}$.\ That all irreducible complex representations of $\GL_d(k)$ arise via Schur functors and powers of the determinant representation is proven in \cite[Ch.\ 15]{fultonharrisrep}. 
    \item In \cite[Ch.\ 8]{fultonyoung}, Schur functors are characterized via universal properties and defined as quotients of $\wedge^{\mu_1}V \otimes \hdots \otimes \wedge^{\mu_r}V$, where $\mu = \alpha'$ is the transpose partition.\ Here $V$ can be any module over a commutative ring $k$.\ The representation theory of $\GL_d(k)$ for $k=\mathbb{C}$ is summarized with reference to \cite{fultonharrisrep}.
    \item In \cite{weyman}, free modules over any commutative ring are treated, but the notational conventions are different.\ The relationship is $\Sigma^{\alpha}V = L_{\mu}V$, $\mu = \alpha'$.\par 
    \noindent The module $L_{\mu}V$ is defined similarly to \cite{fultonyoung} and later shown to be free again, as well as isomorphic to the image of a \emph{Schur map} $\phi_{\mu}$.\ The definition of this map is somewhat combinatorial, for details see \cite[p.\ 37]{weyman}.\ The relationship to \cite{fultonharrisrep} is explained in \cite[Sect.\ 2.2]{weyman}.
    \item An even more general exposition can be found in \cite{akinbuchsbaumweyman}.\ Schur functors are first defined in the free case as images of Schur maps and afterwards described as quotient modules (more precisely:\ as cokernels)  \cite[Sect.\ II.2]{akinbuchsbaumweyman}.\ The same is done in \cite[Sect.\ V.2]{akinbuchsbaumweyman} for finitely generated modules, and the cokernel description allows to conclude that \emph{Schur functors commute with base change}.
\end{itemize} 

\noindent In characteristic zero, the definition from \cite{fultonharrisrep} seems to be the most amenable one.\ Note that $\Sigma^{(m)}V = S^mV$ is the symmetric power whereas the transpose partition yields $\Sigma^{(1,\hdots,1)}V \cong \wedge^mV$.\ It has been mentioned several times that the formation of $\Sigma^{\alpha}$ is functorial, preserves free modules and commutes with base change.\ Over a field of characteristic zero, a dimension formula can be found in \cite[Thm. 6.3]{fultonharrisrep}.\ In particular, $\Sigma^{\alpha}V = 0$ if $\alpha$ has more than $d = \dime V$ entries.\footnote{This is true in arbitrary characteristic, cf.\ the definition in \cite[Sect.\ 2.1]{weyman}.}\ Let $k$ be a commutative ring of characteristic zero.

\begin{thrm} \label{rules}
    Consider two free modules $V,W$ over $k$.
    \begin{enumerate}
        \item[(1)] \cite[Cor.\ 2.3.3]{weyman} There are isomorphisms $S^m(V \otimes W) \cong \bigoplus_{|\alpha| = m} \Sigma^{\alpha}V \otimes \Sigma^{\alpha}W$ and $\wedge^m(V \otimes W) \cong \bigoplus_{|\alpha| = m} \Sigma^{\alpha}V \otimes \Sigma^{\alpha'}W$, natural in $V$ and $W$.
        \item[(2)] \cite[Cor.\ 2.3.4]{weyman} The tensor product of Schur functors $\Sigma^{\alpha}V \otimes \Sigma^{\beta}V$ decomposes into $\bigoplus_{|\nu| = |\alpha| + |\beta|} (\Sigma^{\nu}V)^{\oplus N_{\alpha\beta\nu}}$, where $N_{\alpha\beta\nu}$ are the Littlewood--Richardson coefficients (see also \cite[App.\ A.1]{fultonharrisrep}).\ A special case where $\beta = (m)$ or $\beta = (1,\hdots,1)$ are Pieri's formulas \cite[Cor. 2.3.5]{weyman}, explained in Example \ref{pieriexample}.
        \item[(3)] \cite[Prop.\ 2.3.8 \& 2.3.9]{weyman} The \emph{inner plethysm} problem is solvable in the following easy cases: $S^m(S^2V)$, $S^m(\wedge^2V)$, $\wedge^m(S^2V)$, $\wedge^m(\wedge^2V)$.\ In more detail:
        \begin{gather*}
            S^m(S^2V) \cong \bigoplus_{\substack{|\alpha| = 2m, \\ \alpha_i \text{ even}}} \Sigma^{\alpha}V, \quad \quad S^m(\wedge^2V) \cong \bigoplus_{\substack{|\alpha| = 2m, \\ \alpha_i' \text{ even}}} \Sigma^{\alpha}V, \\
            \wedge^m(S^2V) \cong \bigoplus_{\substack{|\alpha| = 2m,\\ \alpha = (\mathbf{u}|\mathbf{v}),\\ v_i = u_i + 1}}\Sigma^{\alpha'}V, \quad \quad \wedge^m(\wedge^2V) \cong \bigoplus_{\substack{|\alpha| = 2m,\\ \alpha = (\mathbf{u}|\mathbf{v}),\\ u_i = v_i + 1}}\Sigma^{\alpha'}V,
        \end{gather*}
        using the hook notation $\alpha = (\mathbf{u}|\mathbf{v})$ of \cite{weyman} in the last two identities.
    \end{enumerate}
\end{thrm}

\begin{bem} \label{filtrationchar}
    In positive characteristic, items $(1)$ and $(2)$ in Theorem \ref{rules} only remain true with filtrations instead of direct sums.
\end{bem} 

\begin{bsp}[Pieri's formulas] \label{pieriexample}
Consider $\beta = (m)$ and $\Sigma^{\alpha}V \otimes \Sigma^{\beta}V = \Sigma^{\alpha}V \otimes S^mV$ first. Pieri's formulas state that this tensor product decomposes into $\GL_d(k)$-representations $\Sigma^{\nu}V$, where $\nu$ is obtained from $\alpha$ by adding $m$ boxes with no two in the same column. Similarly for $\beta' = (1,\hdots,1)$, where $\Sigma^{\alpha}V \otimes \Sigma^{\beta'}V = \Sigma^{\alpha}V \otimes \wedge^mV$ decomposes into $\Sigma^{\nu}V$ with $\nu$ obtained by adding $m$ boxes with no two in the same row.\par 
\noindent In particular, tensoring with $\det V$ adds one box to $\alpha$ in each row.\ For example, let $d = \dime V = 3$ and $\alpha = (2,1)$.\ Then $\Sigma^{\alpha}V \otimes S^2V \cong \Sigma^{4,1,0}V \oplus \Sigma^{3,2,0}V \oplus \Sigma^{3,1,1}V \oplus \Sigma^{2,2,1}V$: \par \vspace{0.4cm}
\begin{center}
{\scriptsize \ydiagram[*(CadetBlue)]
{2+2}
*[*(white)]{4,1} \hspace{0.8cm }
\ydiagram[*(CadetBlue)]
{2+1,1+1}
*[*(white)]{3,2} \hspace{0.8cm }
\ydiagram[*(CadetBlue)]
{2+1,0,1}
*[*(white)]{3,1,1} \hspace{0.8cm }
\ydiagram[*(CadetBlue)]
{0,1+1,1}
*[*(white)]{2,2,1}}
\end{center} \vspace{0.4cm}
\noindent and similarly $\Sigma^{\alpha}V \otimes \wedge^2V \cong \Sigma^{3,2,0}V \oplus \Sigma^{3,1,1}V \oplus \Sigma^{2,2,1}V$ (notice $\Sigma^{2,1,1,1}V = 0$), visualized by adding the coloured boxes. 
\end{bsp} 

\noindent To obtain all rational representations of $\GL_d(k)$, (duals of) the determinant representation $\det V$ have to be considered.\ In view of Pieri's formula, this allows to define $\Sigma^{\alpha}V$ for any non-increasing sequence $(\alpha_i)_i$ of integers:\ $\Sigma^{\alpha}V := \Sigma^{\alpha + m}V \otimes \det V^{\otimes -m}$, where $(\alpha + m)_i$ is given by $\alpha_i + m$.\ With this definition, the $\Sigma^{\alpha}V$ are precisely the representations of $\GL_d(k)$ with highest weight $\alpha \in \mathbb{Z}^d$, cf.\ \cite[Ch.\ 8.2]{fultonyoung}, \cite[0.2]{kapranov1} and \cite[p.\ 54]{weyman}.

\begin{lemma} \label{dualformula}
    Let $\chara k = 0$ and $\alpha = (\alpha_1,\hdots,\alpha_r)$.\ Then $(\Sigma^{\alpha}V)^{\vee}  \cong \Sigma^{\alpha}(V^{\vee}) \cong \Sigma^{-\alpha}V$, where $-\alpha$ denotes the sequence of integers $(-\alpha_r,\hdots,-\alpha_1)$.
\end{lemma}

\begin{proof}
    This is mentioned in \cite{kapranov1}.\ For a proof, see \cite[Prop. II.4.1 \& II.4.2]{akinbuchsbaumweyman}.
\end{proof}

\noindent As put in \cite[p. 66]{weyman}, \enquote{all the formulas proven [above] are functorial, so they extend to vector bundles}.\ Though this is certainly true, it is worth mentioning that the concept of Schur functors can be applied to all quasi-coherent sheaves or can (and maybe should) be categorified.\ For a quasi-coherent sheaf $\mathcal{F}$ on a ringed space $X$, $\Sigma^{\alpha}\mathcal{F}$ can be defined as the sheaf associated to the presheaf $U \mapsto \Sigma^{\alpha}\Gamma(U,\mathcal{F})$, using the definition of Schur functors over arbitrary commutative rings.\ If $X$ is a scheme, $\Sigma^{\alpha}\mathcal{F}$ is easily seen to be quasi-coherent again (just like for symmetric or wedge powers), essentially because of the compatibility of $\Sigma^{\alpha}$ with base change.\ Theorem \ref{rules} also extends to sheaves and restricts to the full subcategory of finite locally free sheaves.\ Returning to the Grassmannian $X = G = G(k,V)$ of $k$-dimensional quotients and the tautological short exact sequence (\ref{tautsequence}), recall first of all the following classical result of Kapranov:

\begin{thrm} \label{kapranov85main}
    The derived category $\db(G)$ possesses a full \cite[Prop.\ 1.4]{kapranov1} and strong exceptional collection \cite[Prop.\ 2.2]{kapranov1} consisting of objects $\Sigma^{\alpha}K$ (dually: $\Sigma^{\alpha'}Q^{\vee}$), where $\alpha$ is a Young diagram with $\leq k$ columns and $\leq d-k$ rows. 
\end{thrm}

\begin{bem}[The correct ordering] \label{correctordering}
    As explained in \cite{grasscharp}, many different total orderings $\alpha \prec \beta$ on partitions yield an exceptional sequence.\ The important property when considering $\Sigma^{\alpha}Q^{\vee}$ is that $|\beta| < |\alpha|$ implies $\alpha \prec \beta$, i.e.\ that partitions with more entries occur \emph{first}.\ On the subsets of partitions of same length ($|\alpha| = |\beta|$), any total order (e.g.\ the lexicographic one) can be chosen.\ With this convention, $\alpha \prec \beta$ implies $|\beta| \leq |\alpha|$.
\end{bem}

\noindent Fullness of the collection relies on a resolution of the diagonal $\Delta \subseteq G \times G$, whereas strong exceptionality involves computing $\ext^{\ast}_{\mathcal{O}_G}(\Sigma^{\beta}K,\Sigma^{\alpha}K) \cong H^{\ast}(G,\Sigma^{\alpha}K \otimes \Sigma^{-\beta}K)$.\par 
\noindent In general, the locally free sheaf $\Sigma^{\alpha}K \otimes \Sigma^{-\beta}K$ (similarly for $\Sigma^{\alpha}Q^{\vee} \otimes \Sigma^{-\beta}Q^{\vee}$) decomposes via the Littlewood--Richardson rule into sheaves $\Sigma^{\gamma}K$ with $-k \leq \gamma_i \leq k$.\ Instead of $\Sigma^{\gamma}K$, consider more generally the bundles $\Sigma^{\gamma}K \otimes \Sigma^{\delta}Q^{\vee}$.\ Kapranov first relates these bundles to invertible sheaves $\mathcal{L}_{\alpha}$ with $\alpha = (-\gamma,\delta)$ on the full flag variety $F$ lying over $G = G(k,V)$ \cite[2.5]{kapranov1}.\ The cohomology of $\mathcal{L}_{\alpha}$ can then be computed using the following theorem.

\begin{thrm}[Borel--Weil--Bott for $\GL_d$ {\cite[Cor. 4.1.7]{weyman}}] \label{BWB}
    Consider the cohomology $H^{\ast}(G,\Sigma^{\gamma}K \otimes \Sigma^{\delta}Q^{\vee}) \cong H^{\ast}(F,\mathcal{L}_{\alpha})$.\ Let $\rho = (d,d-1,\hdots,2,1) \in \mathbb{Z}^d$ be the Weyl vector and consider $\alpha + \rho \in \mathbb{Z}^d$.\ If this sequence has two identical entries, the cohomology $H^{\ast}(F,\mathcal{L}_{\alpha})$ vanishes completely.\ If not, reorder the sequence using a unique permutation $\sigma$ of length $l = l(\sigma)$ to obtain a non-increasing sequence $\beta :=  \sigma(\alpha + \rho) - \rho$.\ Then $H^{\ast}(F,\mathcal{L}_{\alpha})$ is concentrated in degree $l$ where it equals $\Sigma^{\beta}V^{\vee}$.
\end{thrm}

\noindent The statement concerning $H^{\ast}(F,\mathcal{L}_{\alpha})$ holds for all reductive algebraic groups $H$ over an algebraically closed field of characteristic zero after replacing $F$ by $H/B$ ($B$ a Borel subgroup containing a maximal torus $T$), $\alpha$ by a character of $T$ and the definition of $\mathcal{L}_{\alpha}$ by a more general one etc.\ See \cite[Ch. 4.3]{weyman} or the summary in \cite[Rem. 3.2]{kuznetsovisograssmannians}.

\begin{kor} \label{vanishinga1}
    Let $k=2$ and let $Q$ be the tautological quotient bundle of rank $2$ on $G = G(2,V)$.\ Suppose that $1 \leq a \leq d-2$.\ Then $H^{\ast}(G,\Sigma^{a,b}Q^{\vee}) = 0$ in all degrees.
\end{kor}

\begin{proof}
    According to Theorem \ref{BWB}, this follows from $\alpha = (0,\hdots,0,a,b) \in \mathbb{Z}^d$ and
    \begin{align*}
        \alpha + \rho = (d,\hdots,3;a+2,b+1),
    \end{align*}
    which contains a double entry if $a+2 \in \lbrace 3, \hdots, d \rbrace$.
\end{proof}

\noindent As an application, exceptionality in Theorem \ref{kapranov85main} can be reproven easily for $k=2$:

\begin{lemma}[Kapranov's collection, $\rank Q = 2$] \label{exampleextcomputation}
Assume that $k=2$.\ Then the objects $\lbrace \Sigma^{\alpha}Q^{\vee }\rbrace_{\alpha}$ form a strong exceptional collection in the bounded derived category of the Grassmannian $G(2,d)$, where $\alpha$ is a Young diagram with $\leq d-2$ columns and $\leq 2$ rows.
\end{lemma}

\begin{proof}
The first step for proving exceptionality of the collection is to show
\begin{align*}
\ext^i_{\mathcal{O}_G}(\Sigma^{\beta}Q^{\vee},\Sigma^{\alpha}Q^{\vee}) = \bigoplus_{\gamma = 0}^{\min\lbrace \lambda_{\alpha},\lambda_{\beta} \rbrace} H^i\bigl(G,\Sigma^{\alpha_1 - \beta_2 - \gamma, \alpha_2 - \beta_1 + \gamma}Q^{\vee}\bigr)
\end{align*}
with $\lambda_{\alpha} = \alpha_1 - \alpha_2$ and $\lambda_{\beta} = \beta_1 - \beta_2$.\ This also yields (\ref{extcomputation}).\ The vector bundle $\Sigma^{\alpha}Q^{\vee} \otimes \Sigma^{-\beta}Q^{\vee}$ can be rewritten as $\Sigma^{(\alpha_1,\alpha_2)}Q^{\vee} \otimes \Sigma^{(\lambda_{\beta},0)}Q^{\vee} \otimes (\det Q^{\vee})^{-\beta_1}.$\ The product of the first two factors decomposes into a direct sum according to Pieri's formula since $\Sigma^{(\lambda_{\beta},0)}Q^{\vee}$ is the $\lambda_{\beta}$'th symmetric power.\ Afterwards $\beta_1$ blocks are subtracted in each row from the obtained Young diagrams in the decomposition due to  $(\det Q^{\vee})^{-\beta_1}$.\ A case distinction is made. \par
\noindent Suppose first that $\lambda_{\alpha} \leq \lambda_{\beta}$.\ Then the Young diagrams corresponding to direct summands in the decomposition of $\Sigma^{(\alpha_1,\alpha_2)}Q^{\vee} \otimes \Sigma^{(\lambda_{\beta},0)}Q^{\vee}$ result from adding $\lambda_{\beta}$ boxes to $\alpha$, with no two in the same column. More precisely, the summands are
\begin{align*}
\Sigma^{(\alpha_1+ \lambda_{\beta},\alpha_2)}Q^{\vee} \oplus \Sigma^{(\alpha_1+ (\lambda_{\beta} - 1),\alpha_2 +1)}Q^{\vee} \oplus \hdots \oplus \Sigma^{(\alpha_1+ (\lambda_{\beta} - \lambda_{\alpha}),\alpha_2 + \lambda_{\alpha})}Q^{\vee},
\end{align*}
and not more since the newly added boxes would then overlap column-wise.\ If on the other hand $\lambda_{\beta} \leq \lambda_{\alpha}$, one obtains a decomposition of the first two factors into
\begin{align*}
\Sigma^{(\alpha_1+ \lambda_{\beta},\alpha_2)}Q^{\vee} \oplus \Sigma^{(\alpha_1+ (\lambda_{\beta} - 1),\alpha_2 +1)}Q^{\vee} \oplus \hdots \oplus \Sigma^{(\alpha_1,\alpha_2 + \lambda_{\beta})}Q^{\vee}.
\end{align*}
Subtracting $\beta_1$ everywhere yields the decomposition of $\ext^i_{\mathcal{O}_G}(\Sigma^{\beta}Q^{\vee},\Sigma^{\alpha}Q^{\vee})$. \par 
\noindent If $\alpha = \beta$, then $\ext^i_{\mathcal{O}_G}(\Sigma^{\alpha}Q^{\vee},\Sigma^{\alpha}Q^{\vee}) = \bigoplus_{k=0}^{\lambda_{\alpha}} H^i(G,\Sigma^{k,- k}Q^{\vee}$).\ Only for $k=0$, the cohomology $H^i(G,\Sigma^{k,- k}Q^{\vee}) = H^i(G,\mathcal{O}_G)$ is non-zero in degree $i = 0$.\ All other summands vanish for all $i$ since then $1 \leq k \leq \lambda_{\alpha} \leq d-2$.\ Thus all $\Sigma^{\alpha}Q^{\vee}$ are exceptional. \par 
\noindent If $\alpha \neq \beta$, all sheaves $\Sigma^{a,b}Q^{\vee} = \Sigma^{\alpha_1 - \beta_2 - \gamma, \alpha_2 - \beta_1 + \gamma}Q^{\vee}$ in the above sum fulfill $d-2 \geq a \geq 1$, proving strong exceptionality of Kapranov's sequence in this special case: \par
\noindent Assume $|\beta| < |\alpha|$.\ If $\lambda_{\alpha} \leq \lambda_{\beta}$, then $a \in \lbrace \alpha_1 - \beta_2, \hdots, \alpha_1 - \beta_2 - \lambda_{\alpha} = \alpha_2 - \beta_2 \rbrace$, and
\begin{align*}
d-2 \geq \alpha_1 \geq \alpha_1 - \beta_2 \geq \hdots \geq \alpha_2 - \beta_2 \geq 1,
\end{align*}
where $\alpha_2 - \beta_2 \geq 1$ holds because $\alpha_2 \leq \beta_2$ would lead to $|\beta| < |\alpha| \leq \alpha_1 + \beta_2$, i.e. $\beta_1 < \alpha_1$. Summing up yields $\alpha_2 + \beta_1 < \alpha_1 + \beta_2$ or equivalently $\lambda_{\beta} < \lambda_{\alpha}$, a contradiction.\ The case $\lambda_{\beta} \leq \lambda_{\alpha}$ works identically. \par
\noindent Next, assume that $|\beta| = |\alpha|$.\ There are again two cases to handle and the question in the first case is if indeed $\alpha_2 - \beta_2 \geq 1$.\ Assume first that $\lambda_{\alpha} \leq \lambda_{\beta}$.\ If $\lambda_{\alpha} = \lambda_{\beta}$, then together with $|\alpha| = |\beta|$ this implies $\alpha = \beta$, a contradiction.\ Hence $\lambda_{\alpha} < \lambda_{\beta}$ automatically. Assume $\alpha_2 \leq \beta_2$.\ Then together with $\lambda_{\alpha} < \lambda_{\beta}$, it follows that 
\begin{align*}
\alpha_1 + \beta_2 < \alpha_2 + \beta_1 \leq \beta_2 + \beta_1 \Longrightarrow \alpha_1 < \beta_1.
\end{align*}
But then $|\alpha| = \alpha_1 + \alpha_2 < \beta_1 + \beta_2 = |\beta|$, a contradiction to $|\alpha| = |\beta|$.\ The argument is identical if $\lambda_{\beta} \leq \lambda_{\alpha}$.\ By Corollary \ref{vanishinga1}, none of the $\Sigma^{a,b}Q^{\vee}$'s has any cohomology.
\end{proof}

\noindent The following lemmata were also used in the proofs of Theorem \ref{fullyfaithfultheorem1} and \ref{semiorthogonalsequencethrm}.

\begin{lemma} \label{computationwedge4}
    Recall the short exact sequence $0 \rightarrow Q^{\vee} \rightarrow Q \otimes S^2Q^{\vee} \rightarrow \mathcal{N}' \rightarrow 0$ (\ref{restricteddefiningsequenceN}).\ There is an isomorphism $\det \mathcal{N}' = \wedge^4 \mathcal{N}' \cong (\det Q^{\vee})^{\otimes 2} \cong \Sigma^{2,2}Q^{\vee}$ of locally free sheaves.
\end{lemma}

\begin{proof}
    Note first that $\det Q^{\vee} \otimes \det \mathcal{N}' \cong \det( Q \otimes S^2Q^{\vee})$ \cite[Exc. II.5.16]{hartshorne}.\ Furthermore, there is the general identity for finite locally free sheaves $\mathcal{E}$ and $\mathcal{F}$
    \begin{align*}
        \det(\mathcal{E} \otimes \mathcal{F}) \cong \det(\mathcal{E})^{\otimes \rank \mathcal{F}} \otimes\det(\mathcal{F})^{\otimes \rank \mathcal{E}},
    \end{align*}
    being a special case of Theorem \ref{rules} (1).\ This implies that 
    \begin{align*}
        \det \mathcal{N}' \cong \det Q \otimes (\det Q)^{\otimes 3} \otimes \det(S^2 Q^{\vee})^{\otimes 2} \cong (\det Q)^{\otimes 4} \otimes \det(S^2 Q^{\vee})^{\otimes 2}.
    \end{align*}
    \noindent \textit{Claim.}\ There is an isomorphism $\det S^2Q^{\vee} = \wedge^3S^2Q^{\vee} \cong (\det Q^{\vee})^{\otimes 3}$.\ This implies the statement of the lemma since then $ \det \mathcal{N}' \cong (\det Q)^{\otimes 4} \otimes (\det Q^{\vee})^{\otimes 3\cdot 2} \cong (\det Q^{\vee})^{\otimes (6-4)}$. \par 
    \noindent \textit{Proof of the Claim.} Computing  a decomposition into irreducibles of $\wedge^3S^2Q^{\vee}$ is an instance of Theorem \ref{rules} (3).\ The partitions $\alpha$ occuring in the decomposition are given in hook notation \cite[pp.\ 8--9]{weyman} by  sequences $(u_1 > \hdots > u_r|v_1 > \hdots > v_r)$ with $u_i > 0$ and $v_i = u_i + 1$.\ The $u_i$ (resp.\ $v_i$) stand for the diagonal arm (resp. leg) lengths of the Young diagram associated to $\alpha$.\ In this fashion, $u_1 = 1$ leads to $\alpha = (1|2) = (1,1)$, so that $|\alpha| = 2$.\ Next, let $u_1 = 2$.\ Only partitions with $|\alpha| = 6$ are needed, and indeed the only such partition is $(2,1|3,2) = (2,2,2)$.\footnote{Since $|\alpha| = 4$ for $(2|3) = (2,1,1)$ and $|\alpha| > 6$ for $u_1 \geq 3$.} This means $\wedge^3S^2Q^{\vee} \cong \Sigma^{(2,2,2)'}Q^{\vee} \cong \Sigma^{3,3}Q^{\vee}$.
\end{proof}

\begin{lemma} \label{computationwedge2}
    There is an isomorphism of locally free sheaves between  $\wedge^2(S^2Q^{\vee} \otimes Q)$ and the direct sum $(\Sigma^{3,-1}Q^{\vee})^{\oplus 2} \oplus (\Sigma^{1,1}Q^{\vee})^{\oplus 2} \oplus S^2Q^{\vee}$.
\end{lemma}

\begin{proof}
    Theorem \ref{rules} will be applied several times, first of all (1) to obtain
    \begin{align*}
        \wedge^2(S^2Q^{\vee} \otimes Q) \cong \bigl(S^2S^2Q^{\vee} \otimes \wedge^2 Q\bigr) \oplus \bigl(\wedge^2S^2Q^{\vee} \otimes S^2 Q\bigr).
    \end{align*}
    By Theorem \ref{rules} (3), $S^2S^2Q^{\vee} \cong S^4Q^{\vee} \oplus \Sigma^{2,2}Q^{\vee}$ and also $\wedge^2S^2Q^{\vee} \cong \Sigma^{(2,1,1)'}Q^{\vee} \cong \Sigma^{3,1}Q^{\vee}$ like in the proof of Lemma \ref{computationwedge4}.\ Therefore 
    \begin{align*}
        S^2S^2Q^{\vee} \otimes \wedge^2 Q \cong S^2S^2Q^{\vee} \otimes (\det Q^{\vee})^{\otimes -1} \cong \Sigma^{3,-1}Q^{\vee} \oplus \Sigma^{1,1}Q^{\vee}
    \end{align*}
    and $\wedge^2S^2Q^{\vee} \otimes S^2 Q \cong \wedge^2S^2Q^{\vee} \otimes S^2 Q^{\vee} \otimes (\det Q^{\vee})^{\otimes -2}$ decomposes into
    \begin{align*}
        \Sigma^{3,-1}Q^{\vee} \oplus \Sigma^{2,0}Q^{\vee} \oplus \Sigma^{1,1}Q^{\vee} = \Sigma^{3,-1}Q^{\vee} \oplus S^2Q^{\vee} \oplus \Sigma^{1,1}Q^{\vee}.
    \end{align*}
    Collecting all summands together yields the claim of the lemma.
\end{proof}

\noindent Lemma \ref{tangentbundlesimple} is needed in Appendix \ref{AppendixB}. It states that the tangent bundle of $G$ is \textit{simple}.

\begin{lemma} \label{tangentbundlesimple}
    The cotangent bundle of the Grassmannian fulfills $\enndo_{\mathcal{O}_G}(\Omega_{G/k}) \cong k \cdot \id$.
\end{lemma}

\begin{proof}
    Note that by Proposition \ref{tangentsheafgrass} (2), $\enndo_{\mathcal{O}_G}(\Omega_{G/k}) \cong H^0(G, K \otimes K^{\vee} \otimes Q \otimes Q^{\vee})$. Using Lemma \ref{dualformula} and Pieri's formula, there is an isomorphism 
\begin{align*}
 K \otimes K^{\vee} \otimes Q \otimes Q^{\vee} \cong \bigl( \mathcal{O}_G \oplus \Sigma^{(1,0,\hdots,0,-1)}K \bigr) \otimes \bigl( \mathcal{O}_G \oplus \Sigma^{(1,0,\hdots,0,-1)}Q\bigr), \quad \rank K, \rank Q \geq 2.
\end{align*} 
Thus $\enndo_{\mathcal{O}_G}(\Omega_{G/k})$ splits up into the direct sum of $H^0(G,\mathcal{O}_G)$, $H^0(G,\Sigma^{(1,0,\hdots,0,-1)}K)$, $H^0(G,\Sigma^{(1,0,\hdots,0,-1)}Q^{\vee})$ and $H^0(G,\Sigma^{(1,0,\hdots,0,-1)}K \otimes \Sigma^{(1,0,\hdots,0,-1)}Q^{\vee})$ if $\rank K, \rank Q \geq 2$. \par
\noindent Of course, $H^{\ast}(G,\mathcal{O}_G) \cong k[0]$.\ The bundles $\Sigma^{(1,0,\hdots,0,-1)}K$ and $\Sigma^{(1,0,\hdots,0,-1)}Q^{\vee}$ have no cohomology at all\footnote{The bundles may not even show up if $K$, respectively $Q$ is a line bundle.} since $(1,0,\hdots,0,-1)$ contains the positive entry $1$.\ It remains to compute $H^0(G,\Sigma^{(1,0,\hdots,0,-1)}K \otimes \Sigma^{(1,0,\hdots,0,-1)}Q^{\vee})$ using Theorem \ref{BWB} if $\rank K, \rank Q \geq 2$.\ Note that
\begin{align*}
\alpha &= (1,0,\hdots,0,-1;1,0,\hdots,0,-1) \in \mathbb{Z}^{d-k} \times \mathbb{Z}^k, \\
\alpha + \rho &= (d+1,d-1,\hdots,k+2,k;k+1,k-1,\hdots,2,0).
\end{align*}
No entries of $\alpha + \rho$ coincide and the $\sigma \in S_d$ from Theorem \ref{BWB} is the transposition permuting the entries $k$ and $k+1$ of $\alpha + \rho$. Hence $H^{\ast}(G,\Sigma^{(1,0,\hdots,0,-1)}K \otimes \Sigma^{(1,0,\hdots,0,-1)}Q^{\vee})$ lives only in degree $l(\sigma) = 1$, where it equals $\Sigma^{(1,0,\hdots,0,-1)}V^{\vee}$.\ To sum up,  
\begin{align*}
\ext_{\mathcal{O}_G}^{\ast}(\Omega_{G/k},\Omega_{G/k}) \cong k[0] \oplus \Sigma^{(1,0,\hdots,0,-1)}V^{\vee}[-1] \quad \text{if } \rank K, \rank Q \geq 2,
\end{align*}
and $\enndo_{\mathcal{O}_G}(\Omega_{G/k}) = \ext_{\mathcal{O}_G}^{0}(\Omega_{G/k},\Omega_{G/k})$ is one-dimensional in any case.\ Note that neither the cotangent nor the tangent bundle on the Grassmannian $G$ is rigid if $\rank K, \rank Q \geq 2$.
\end{proof}

\section{Tangent Spaces, Differentials and Relative Ext} \label{AppendixB}

The subsections in this Appendix can be understood independently from each other.

\subsection{Derived pullbacks of short exact sequences associated to normal bundles}

\begin{lemma}  \label{technicalpullbacksymmetric}
    Let $Y$ be a reduced, local complete intersection inside a smooth variety $X$. If $\mathcal{J}$ denotes the ideal sheaf of the closed immersion $i: Y \hookrightarrow X$, the derived pullback $Li^{\ast}$ applied to $0 \rightarrow \mathcal{J}^2 \rightarrow \mathcal{J} \rightarrow i_{\ast}(\mathcal{N}_{Y/X}^{\vee}) \rightarrow 0$ induces a long exact sequence
\begin{align*}
\wedge^2\mathcal{N}_{Y/X}^{\vee} \overset{\kappa}{\longrightarrow} \mathcal{N}_{Y/X}^{\vee \otimes 2} \longrightarrow S^2\mathcal{N}_{Y/X}^{\vee} \longrightarrow \mathcal{N}_{Y/X}^{\vee} \overset{\sim}{\longrightarrow} \mathcal{N}_{Y/X}^{\vee} \longrightarrow 0
\end{align*}
such that $\kappa$ is the natural homomorphism $\kappa: a \wedge b \mapsto a \otimes b - b \otimes a$.
\end{lemma}

\begin{proof}
    Locally around each point $P \in X$, $Y$ is the vanishing locus $V(s)$ of a (regular, see \cite[Thm. 8.21A(c)]{hartshorne}) section $s \in H^0(\mathcal{E})$ for some locally free sheaf $\mathcal{E}$ on an open neighbourhood $U$ of $P$. This implies existence of a Koszul resolution on $U$. Shrinking $U$ further, it can be assumed that $U$ is affine and also that $\mathcal{E} \cong \mathcal{O}_U^{\oplus c}$, where $c = \codim(Y,X)$.\ Let $i'$ be the restriction of $i$ to $U \cap Y \hookrightarrow U$. Since restriction to an open subscheme is exact, there are functorial isomorphisms 
    \begin{align*}
        (L^1i^{\ast}\mathcal{J})_{|U \cap Y} \cong L^1 {i'}^{\ast}\mathcal{J}_{|U} \quad \text{ and } \quad (L^1i^{\ast}i_{\ast}\mathcal{N}_{Y/X}^{\vee})_{|U \cap Y} \cong L^1 {i'}^{\ast}({i'}_{\ast}\mathcal{N}_{U \cap Y/U}^{\vee}),
    \end{align*}
    so the whole question becomes affine-local.\ By definition, $L^1i'^{\ast}\bigl(\mathcal{J}_{|U} \rightarrow {i'}_{\ast}\mathcal{N}_{U \cap Y/U}^{\vee}\bigr)$ can be calculated from free resolutions on $U$, even explicitly in the current situation: Resolve $\mathcal{J}_{|U}$ by the truncated Koszul complex and ${i'}_{\ast}\mathcal{N}_{U \cap Y/U}^{\vee} \cong \mathcal{E}^\vee \otimes \mathcal{O}_{U \cap Y}$ by the standard Koszul complex tensored with $\mathcal{E}^{\vee} \cong \mathcal{O}_U^{\oplus c}$ after choosing a basis $\lbrace e_i \rbrace_i$ for $\mathcal{O}_U^{\oplus c}$. This gives
    \begin{center}
    \begin{tikzcd}
    \hdots \arrow[r] & \wedge^3\mathcal{E}^\vee  \arrow[r]                          & \wedge^2\mathcal{E}^\vee  \arrow[r] \arrow[d, dashed] & \mathcal{E}^\vee  \arrow[r] \arrow[d, "="] & \mathcal{J}_{|U} \arrow[r] \arrow[d]                 & 0 \\
    \hdots \arrow[r] & \mathcal{E}^\vee \otimes \wedge^2\mathcal{E}^\vee  \arrow[r] & \mathcal{E}^\vee  \otimes \mathcal{E}^\vee  \arrow[r] & \mathcal{E}^\vee  \arrow[r]                & {i'}_{\ast}\mathcal{N}_{U \cap Y/U}^{\vee} \arrow[r] & 0,
    \end{tikzcd}
    \end{center}
    where commutativity of the right square is tautologically fulfilled. Only the dashed arrow has to be constructed.\footnote{This suffices because either the next arrows can be constructed in an identical fashion or one uses that on the affine scheme $U$, projective resolutions exist (only the degree $(-1)$-part is of interest here).}\ To this end, just send $e_i \wedge e_j \mapsto e_i \otimes e_j - e_j \otimes e_i$, which works since the horizontal arrows are contraction with the components of the regular section $s$.
\end{proof}

\subsection{Different tangent space descriptions for Hilbert schemes} \label{subsectionhilberttangent}

\noindent Since general Hilbert and Quot schemes are constructed inside possibly very large Grassmannians, it makes sense to talk about tangent spaces and tangent sheaves of Grassmannians first.

\begin{prop} \label{tangentsheafgrass}
    Let $\mathbb{G} = \mathbb{G}(k,\mathcal{E}) \overset{p}{\longrightarrow} X$ be the relative Grassmannian of $k$-dimensional quotients over a base scheme $X$, where $\mathcal{E} \in$ \emph{$\coh(X)$}. 
    \begin{enumerate}
        \item[(1)] The sheaf of relative differentials associated to the projection $p$ is $\Omega_p \cong \mathcal{H}om(\mathcal{Q},\mathcal{K})$, where $\mathcal{K}$ and $\mathcal{Q}$ are the tautological bundles on $\mathcal{G}$.
    \end{enumerate}
    In the case where $X = \spec k$ for a field $k$ and $\mathcal{E} \cong \widetilde{V}$, $V \cong k^d$, this implies:
    \begin{enumerate}
        \item[(2)] The sheaf of differentials $\Omega_{G/k}$ on the absolute Grassmannian $G = G(k,d)$ is locally free and isomorphic to $\mathcal{H}om(Q,K) \cong Q^{\vee} \otimes K$.
        \item[(3)] Restricting the sheaf of differentials and the isomorphism of (a) to a $k$-point $[W]$ in $G(k)$ yields the identification of tangent spaces $T_{[W]}G \cong \homo_k(W,V/W)$.
        \item[(4)] Without using (3), there is a natural one-to-one-correspondence between pointed morphisms from the dual numbers $\spec k[\overline{\varepsilon}] = \spec k[\varepsilon]/(\varepsilon^2)$ to $G$ (mapping $(\varepsilon)$ to the $k$-point $[W]$) and the vector space $\homo_k(W,V/W)$. 
    \end{enumerate}
    Furthermore, the two a priori different identifications of $T_{[W]}G$ with $\homo_k(W,V/W)$ coincide, i.e.\ there is a commutative diagram of the following form relating (3) and (4).
\begin{center}
\begin{tikzcd}
{\Omega_{G/k}([W])^{\vee}} \arrow[r, "(3)"] \arrow[d, "\sim"'] & {\mathcal{H}om(K,Q)([W])} \arrow[r,  "\sim"]                                                    & {\homo_k(W,V/W)} \\
{T_{[W]}G} \arrow[r, "\sim"]                                   & \mor_k\bigl((\spec k[\overline{\varepsilon}], (\varepsilon)),(G,[W])\bigr) \arrow[ru, "(4)"'] &                 
\end{tikzcd}
\end{center}
\end{prop}

\begin{proof}
    Claim (1) is proven in \cite[Prop.\ 4.6.1]{sernesi} for $\mathcal{E}$ locally free and in \cite{lehncotangent} for a general coherent sheaf $\mathcal{E}$.\ Then (2) and (3) are clear.\ Claim (4) is a standard fact, the correspondence being explained e.g. in \cite[Ex. 8.22]{goertzwedhorn}.\ The main statement (needed in the proof of Theorem \ref{maintheoremnormalbundle}) is that the two different approaches of computing the tangent space to an absolute Grassmannian yield the same, which can be proven as follows. \par
    \noindent Firstly, there is yet another way to describe the tangent bundle of $G$, namely by patching trivial vector bundles, i.e.\ by globalizing the functorial tangent space description (4) to an isomorphism $\psi: \mathcal{T}_{G/k} \cong \mathcal{H}om(K,Q)$ of vector bundles:\ Abstractly, this can be accomplished as in \cite[Thm.\ 17.46]{goertzwedhorn2}, relying implicitly on the open cover of the Grassmannian $G$ by affine spaces \cite[Lemma 8.13]{goertzwedhorn} over which the tautological bundles trivialize.\  The reader preferring a more geometric argument can also follow the construction of $\psi$ in \cite[Thm.\ 3.5]{3264} (globalizing the functorial description of $T_{[W]}G$ as wanted).\footnote{The proof in \cite{3264} is given for classical varieties, but can easily be rephrased scheme-theoretically.}\par
    \noindent Then one uses  $\enndo_{\mathcal{O}_G}(\mathcal{T}_{G/k})= k \cdot \id$ (Lemma \ref{tangentbundlesimple}) to conclude that $\psi$ and the dual of the isomorphism $\Omega_p \cong \mathcal{H}om(Q,K)$ from (2) are scalar multiples of each other (in particular, this holds true on all fibres).
\end{proof}

\noindent In fact, the analogue to Proposition \ref{tangentsheafgrass} is true for Hilbert schemes, where the relative case (1) is neglected here (but remains true according to \cite[Thm. 3.1]{lehncotangent}).

\begin{prop} \label{tangentsheafhilbert}
    Let $X$ be a smooth projective variety of dimension $d$ over a field $k$ and let $X^{[n]} = \hilb^n_{X/k}$ be the Hilbert scheme of $n$ points.
    \begin{enumerate}
        \item[(1)] \cite{lehncotangent} The sheaf of differentials $\Omega_{X^{[n]}/k}$ on the Hilbert scheme $X^{[n]}$ is isomorphic to the locally free sheaf $\mathcal{E}xt^d_{\overline{\pi}}(\mathcal{O}_{\Xi},\mathcal{I}_{\Xi} \otimes \omega_{\overline{\pi}})$, where $\overline{\pi}: X \times X^{[n]} \rightarrow X^{[n]}$ is the projection and the relative $\mathcal{E}xt$-sheaf is defined in Lemma \ref{defirelativeext}.
        \item[(2)] Restricting the isomorphism of (1) to a $k$-point $[W] \in X^{[n]}(k)$ yields an isomorphism $T_{[W]}X^{[n]} \cong \homo_{\mathcal{O}_X}(\mathcal{I}_W,\mathcal{O}_W) \cong \homo_{\mathcal{O}_X}(\mathcal{I}_W/\mathcal{I}_W^2,\mathcal{O}_X/\mathcal{I}_W)$.
        \item[(3)] Without using (2) and also for $X$ quasi-projective, there is a natural one-to-one-correspondence between pointed morphisms from the dual numbers $\spec k[\overline{\varepsilon}]$ to $X^{[n]}$ (assigning $[W]$ to $(\varepsilon)$) and the vector space $\homo_{\mathcal{O}_X}(\mathcal{I}_W,\mathcal{O}_W)$.
    \end{enumerate}
    Furthermore, the two a priori different identifications in $(2)$ and $(3)$ coincide, which is important in the proof of Theorem \ref{maintheoremnormalbundle}.
\end{prop}

\begin{proof}
    The statement (1) is \cite[Thm.\ 3.1]{lehncotangent}, proven by realizing the Hilbert scheme inside a Grassmannian.\ This will be explained below.\ Claim (2) follows from (1) using the base change argument for relative $\mathcal{E}xt$ from Proposition \ref{descriptionrestrictedcotangentsheaf} or Lemma \ref{extbasechangeisnatural} together with Serre duality.\ The classical description (3) is explained in \cite[Ch.\ IX.5]{geomalgcurv}:\par 
    \noindent First of all, it is possible to reduce to the affine case $X = \spec A$, where a zero-dimensional subscheme of $X$ corresponds to an ideal $I \trianglelefteq A$ such that $\dime_k A/I = n$.\ A tangent vector $\spec k[\overline{\varepsilon}] \rightarrow X^{[n]}$ at the ideal $I$ corresponds to an ideal $J \trianglelefteq A[\overline{\varepsilon}] = A \otimes_k k[\overline{\varepsilon}]$ restricting to $I$ mod $\varepsilon$ such that $A[\overline{\varepsilon}]/J$ is flat of degree $n$ over $k[\overline{\varepsilon}]$.\ These embedded first-order deformations \cite[Lemma 5.8]{geomalgcurv} can be identified with $\homo_A(I,A/I)$ using the following maps, inverse to each other.\ The first map is
\begin{align*}
\bigl(J \trianglelefteq A[\overline{\varepsilon}]\bigr) \mapsto (\varphi_J: I \longrightarrow A/I), \quad \varphi_J(i) = \overline{h},\quad  i- \varepsilon h \text{ is a lift of } i \in I. 
\end{align*}
\noindent Writing $I = (i_1,\hdots,i_r)$, the second map is defined by
\begin{align*}
    (\varphi: I \longrightarrow A/I) \mapsto \bigl(J_{\varphi} \trianglelefteq A[\overline{\varepsilon}]\bigr), \quad J_{\varphi} = (j_1,\hdots,j_r), \quad j_k := i_k - \varepsilon a_k, \quad \varphi(i_k) = \overline{a_k}.
\end{align*}
\noindent It remains to explain why the dual of the restriction of the isomorphism in (1) and the functorial isomorphism of (3) coincide in the same way as in Proposition \ref{tangentsheafgrass}, i.e.\ why 
\begin{center}
\begin{tikzcd}
{T_{[W]}X^{[n]} = \Omega_{X^{[n]}/k}([W])^{\vee}} \arrow[r, "{(2)}"] & {\ext_{\mathcal{O}_X}^d(\mathcal{O}_W,\mathcal{I}_W \otimes \omega_X)^{\vee}} \arrow[r, "\sim"] & {\homo_{\mathcal{O}_X}(\mathcal{I}_W,\mathcal{O}_W)}
\end{tikzcd}
\end{center}
is the same as the isomorphism of (3) deduced via dual numbers.\ To this end, Lehn's proof of (1) in the case of Hilbert schemes is summarized first:
\begin{enumerate}
\item[(i)] For a large enough integer $m > 0$, the Hilbert scheme $X^{[n]} = \hilb^n_{X/k}$ is a closed subscheme of the Grassmannian $G_m = G\bigl(n,H^0(X,\mathcal{O}_{X}(m)\bigr)$.\ The closed embedding $\mathcal{H}ilb^n_{X/k} \hookrightarrow h_{G_m}$ of functors is given on $T$-points as follows\footnote{Of course, the important thing is that the map to $q_{\ast}(\mathcal{O}_W(m))$ stays surjective.\ Trying to ensure this for $m \gg 0$ leads to the concept of $m$-regularity.} for all $m \gg 0$.
\begin{align*}
\left[ q: W \hookrightarrow X \times T \overset{f}{\rightarrow} T \right] \mapsto \left[ f_{\ast}(\mathcal{O}_{X \times T}(m)) = H^0(X,\mathcal{O}(m)) \otimes \mathcal{O}_T \twoheadrightarrow q_{\ast}(\mathcal{O}_W(m))\right]
\end{align*}
\item[(ii)] Inside a finite product $\prod_{m_i \gg 0} G_{m_i}$ of such Grassmannians, $X^{[n]}$ can also be realized as a connected component of the vanishing locus of a morphism between locally free sheaves \cite[Lemma 3.2]{lehncotangent}.\ The cotangent sheaves of Grassmannians and those of vanishing loci have convenient global descriptions, and $\Omega_{X^{[n]}/k}$ is a cokernel of a morphism $\Lambda$ between these.
\item[(iii)] A further examination of $\Lambda$ shows that some components of it are already epimorphisms, wherefore the Snake lemma yields a description of $\Omega_{X^{[n]}/k}$ as the cokernel of an easier (to describe) morphism $\overline{\Lambda}$ \cite[Cor.\ 3.6]{lehncotangent} of the following form.
\begin{align*}
    \overline{\Lambda}: \bigoplus_{i=1}^l \mathcal{H}om(\mathcal{O}_{m_i},\overline{\pi}_{\ast}\mathcal{I}'(m_i)) \longrightarrow \mathcal{H}om(\mathcal{O}_{m_0},\mathcal{I}_{m_0})
\end{align*}
\noindent Here $\mathcal{O}_{m_i} = \overline{\pi}_{\ast}\mathcal{O}_{\Xi}(m_i)$, $\mathcal{I}_{m_i} = \overline{\pi}_{\ast}\mathcal{I}_{\Xi}(m_i)$ and $\mathcal{I}' = \kernel \bigl(\overline{\pi}^{\ast}(\mathcal{I}_{m_0})(-m_0) \twoheadrightarrow \mathcal{I}_{\Xi}\bigr)$, where $\overline{\pi}:X \times X^{[n]} \rightarrow X^{[n]}$ is the projection and $\Xi_n$ denotes the universal family.\
\item[(iv)] Applying (iii), there exists a commutative square with exact rows
\begin{center}
\begin{tikzcd}
{\mathcal{E}xt_{\overline{\pi}}^d(\mathcal{O}_{\Xi},(\overline{\pi}^{\ast}\mathcal{I}_{m_0})(-m_0) \otimes \omega_{\overline{\pi}})} \arrow[d, "\text{G.V.}"] \arrow[r] & {\mathcal{E}xt_{\overline{\pi}}^d(\mathcal{O}_{\Xi},\mathcal{I}_{\Xi}\otimes \omega_{\overline{\pi}})} \arrow[d, "\sim", dashed] \arrow[r] & 0  \\
{\mathcal{H}om_{{X^{[n]}}}(\mathcal{O}_{m_0},\mathcal{I}_{m_0} )} \arrow[r]                                                                         & \coker \overline{\Lambda} \cong \Omega_{X^{[n]}/k} \arrow[r]                                                                                                        & 0.
\end{tikzcd}
\end{center}
\noindent The right vertical map is induced from the left one and yields the desired global isomorphism (1).\ The abbreviation G.V. stands for Grothendieck--Verdier duality.
\end{enumerate}
Going to the fibres at the $k$-point $[W] \in X^{[n]}(k)$ and dualizing, this diagram becomes
\begin{center}
\begin{tikzcd}
{\ext_{\mathcal{O}_X}^d(\mathcal{O}_W(m),I_m \otimes_k\omega_X)^{\vee}} & {\ext_{\mathcal{O}_X}^d(\mathcal{O}_W,\mathcal{I}_W \otimes \omega_X)^{\vee}} \arrow[l] & 0 \arrow[l] \\
{\homo_k(I_m,H^0(\mathcal{O}_W)_m)} \arrow[u, "\sim"']                                  & {T_{[W]}X^{[n]}} \arrow[u, "{(2)}"'] \arrow[l]                                 & 0, \arrow[l]
\end{tikzcd}
\end{center}
where $m = m_0 \gg 0$ and $I_m = \mathcal{I}_m([W]) = H^0(X,\mathcal{I}_W(m))$ is the degree $m$-part of the ideal cutting out $W \subseteq X$. Even more precisely, this diagram can be enlarged to \par \vspace{0.25cm}
\adjustbox{scale=0.85,center}{
\begin{tikzcd}
{\ext_{\mathcal{O}_X}^d(\mathcal{O}_W(m),I_m \otimes_k  \omega_X)^{\vee}} &                                                                                  & {\ext_{\mathcal{O}_X}^d(\mathcal{O}_W,\mathcal{I}_W \otimes \omega_X)^{\vee}} \arrow[ll, "j"] \arrow[ld, "\text{S.D.}"'] & 0 \arrow[l] \\
                                                                          & {\homo_{\mathcal{O}_X}(\mathcal{I}_W, \mathcal{O}_W)} \arrow[ld, "\alpha", hook'] &                                                                                                              &             \\
{\homo_k(I_m,H^0(\mathcal{O}_W)_m)} \arrow[uu, "\gamma"]                        &                                                                                  & {T_{[W]}X^{[n]}} \arrow[uu, "(2)"'] \arrow[ll] \arrow[lu, "(3)"] \arrow[ld, "\beta", hook']                 & 0, \arrow[l] \\
                                                                          & {T_{[H^0(\mathcal{O}_W)_m]}G_m} \arrow[lu, "(\ref{tangentsheafgrass})"]                          &                                                                                                              &            
\end{tikzcd} }
\noindent where a priori only the rectangle and the triangle on the bottom involving $\beta$ commute. The monomorphism $\alpha$ is given by applying $H^0(X,(-) \otimes \mathcal{O}_X(m))$, which is the same as
\begin{align*}
\homo_{\mathcal{O}_X}(\mathcal{I}_W, \mathcal{O}_W) \hookrightarrow \homo_{\mathcal{O}_X}(I_m \otimes_k \mathcal{O}_X,\mathcal{O}_W(m)) \cong \homo_k(I_m,H^0(\mathcal{O}_W)_m).
\end{align*}
The inclusion $\beta$ is the differential of the embedding $X^{[n]} \hookrightarrow G_m$.\ Remember that the goal was to explain why $(3) = (2) \circ (\text{S.D.})$ holds.\ This can be checked after post-composing (3) and $(\text{S.D.}) \circ (2)$ with the monomorphism $\gamma \circ \alpha$.\par 
\noindent But $\alpha \circ (3)$ is the same\footnote{That is, via the functorial tangent space descriptions, the differential of the inclusion $X^{[n]} \hookrightarrow G_m$ sends a homomorphism $\mathcal{I}_W \rightarrow \mathcal{O}_W$ to its degree-$m$-part, cf.\ \cite[p.\ 13]{moduliofcurves}.} as $\beta$ composed with one of the equivalent (!)\ tangent space identifications of the Grassmannian from Proposition\ \ref{tangentsheafgrass}, thus equal to the horizontal map $T_{[W]}X^{[n]} \rightarrow \homo_k(I_m,H^0(\mathcal{O}_W)_m)$. This implies $\gamma \circ \alpha \circ (3) = j \circ (2)$. \par 
\noindent It remains to see why $\gamma \circ \alpha \circ (\text{S.D.}) \circ (2) = j \circ (2)$ as well, which holds true if $\gamma^{-1}$ is the composition of Serre-duality with adjunction, i.e. if $\gamma^{-1}$ coincides with
\begin{align*}
\ext_{\mathcal{O}_X}^d(\mathcal{O}_W(m),I_m \otimes_k  \omega_X)^{\vee} \cong \homo_{\mathcal{O}_X}(I_m \otimes_k \mathcal{O}_X, \mathcal{O}_W(m)) \cong \homo_k(I_m,H^0(\mathcal{O}_W)_m).
\end{align*}
In other words, $\gamma$ (the dual of the pullback of G.V.\ duality to $[W]$ in $X^{[n]}(k)$ plus a natural isomorphism of vector spaces) should be expressible via the chain of isomorphisms above.\par 
\noindent Recall that Grothendieck--Verdier duality in the absolute case is Serre duality, see for example \cite[Thm.\ 25.127]{goertzwedhorn2}.\ Thus, one only has to argue why pulling back relative Grothendieck--Verdier duality to a $k$-point yields absolute Grothendieck--Verdier duality.\ This follows from Lemma \ref{basechangegvduality} applied to $\mathcal{F} = \mathcal{O}_{\Xi}(m)$, $\mathcal{K} =  \mathcal{I}_{m}$ and $p = \overline{\pi}$.
\end{proof}

\subsection{Naturality of base change for relative Ext, Grothendieck--Verdier duality}

\begin{lemma} \label{extbasechangeisnatural}
Let $f:X \rightarrow Y$ be a smooth and projective morphism with $d$-dimensional fibres and consider the restriction of $\mathcal{E}xt_f^d(-,-)$ to coherent sheaves that are flat over $Y$.\ Then the relative $\mathcal{E}xt$'s commute with base change, \emph{naturally} in both arguments.\ In other words, for every point $\iota: \spec k(y) \rightarrow Y$ there is an isomorphism of bifunctors 
\begin{align*}
    \mathcal{E}xt_f^d(-,-)\otimes_{\mathcal{O}_Y} k(y) \cong \ext_{\mathcal{O}_{X_y}}^d((-)_y,(-)_y).
\end{align*}
\end{lemma}

\begin{proof}
Identically to the proof of Proposition \ref{descriptionrestrictedcotangentsheaf}, the following steps can be performed in order to obtain base change isomorphisms:
\begin{align*}
\iota^{\ast}\mathcal{H}^d(Rf_{\ast}R\mathcal{H}om(-,-)) &\cong \mathcal{H}^d(L\iota^{\ast}Rf_{\ast}R\mathcal{H}om(-,-)) \cong \mathcal{H}^d(R\tilde{f}_{\ast}R\mathcal{H}om(L\tilde{\iota}^{\ast}(-),L\tilde{\iota}^{\ast}(-))) \\
&\cong \mathcal{H}^d(R\tilde{f}_{\ast}R\mathcal{H}om(\tilde{\iota}^{\ast}(-),\tilde{\iota}^{\ast}(-))) \cong \ext_{\mathcal{O}_{X_y}}^d((-)_y,(-)_y).
\end{align*}
Here $\widetilde{(-)}$ denotes the base change of a morphism and the second-to-last identity uses the flatness assumptions.\ The point of this lemma is to understand why the first isomorphism is natural in both arguments since the spectral sequence argument used in Proposition \ref{descriptionrestrictedcotangentsheaf} makes it a priori difficult to trace back all involved morphisms.\par 
\noindent Abstractly, one can refer to the naturality statement in \cite[Prop.\ F.212]{goertzwedhorn2}.\ Alternatively, it is not difficult to see that naturality boils down to proving the following elementary statement:\ For a morphism $h: \mathcal{A}^{\bullet} \rightarrow \mathcal{B}^{\bullet}$ in the (non-derived) category of bounded complexes over $\coh(Y)$, where $\mathcal{A}^{\bullet}$ and $\mathcal{B}^{\bullet}$ are locally free and have no cohomology in degrees $\geq d + 1$, there are identifications $\iota^{\ast}\mathcal{H}^d(\mathcal{A}^{\bullet}) \cong \mathcal{H}^d(\iota^{\ast}\mathcal{A}^{\bullet})$ and $\iota^{\ast}\mathcal{H}^d(\mathcal{B}^{\bullet}) \cong \mathcal{H}^d(\iota^{\ast}\mathcal{B}^{\bullet})$ so that $\iota^{\ast}\mathcal{H}^d(\mathcal{A}^{\bullet}) \rightarrow \iota^{\ast}\mathcal{H}^d(\mathcal{B}^{\bullet})$ and $\mathcal{H}^d(\iota^{\ast}\mathcal{A}^{\bullet}) \rightarrow \mathcal{H}^d(\iota^{\ast}\mathcal{B}^{\bullet})$ agree.
\end{proof}

\begin{lemma} \label{basechangegvduality}
Let $p:X \rightarrow Y$ be a smooth, proper morphism of varieties over $k$ of relative dimension $d$ and consider the following cartesian diagram induced by a $k$-point of $Y$.
    \begin{center}
    \begin{tikzcd}
    X' \arrow[d, "p'"'] \arrow[r, "\iota'"] \arrow[rd, "\square", phantom] & X \arrow[d, "p"] \\
    \spec k \arrow[r, "\iota"]                                       & Y               
    \end{tikzcd}
    \end{center}
Assume that $\mathcal{F} \in$ \emph{$\coh(X)$} is flat over $Y$, that $\mathcal{K}$ and $Rp_{\ast}\mathcal{F}$ are locally free sheaves on $Y$ concentrated in degree zero and that the natural map $\iota^{\ast}p_{\ast}\mathcal{F} \rightarrow p'_{\ast}\iota'^{\ast}\mathcal{F}$ is an isomorphism.\ Observe that by relative Grothendieck--Verdier duality for $p$ ($G_p$ for short), 
\begin{align*} 
    \mathcal{H}^d(G_p): \mathcal{E}xt^d_p(\mathcal{F},p^{\ast}\mathcal{K} \otimes \omega_p) \cong \mathcal{H}^d(R\mathcal{H}om(Rp_{\ast}\mathcal{F}, \mathcal{K})[-d]) = \mathcal{H}om_{\mathcal{O}_Y}(p_{\ast}\mathcal{F},\mathcal{K}).
\end{align*}
Then the pullback of $\mathcal{H}^d(G_p)$ along $\iota$ combined with base change $\alpha: L\iota^{\ast}Rp_{\ast} \cong Rp'_{\ast}L\iota'^{\ast}$ agrees with Grothendieck--Verdier duality for $p'$.\ In other words, $\iota^{\ast}\mathcal{H}^d(G_p)$ agrees with the homomorphism  $\mathcal{H}^d(G_p'): \ext^d_{\mathcal{O}_{X'}}(F,p'^{\ast} K \otimes \omega_X) \cong \homo_k(p'_{\ast}F,K)$, where $F$ and $K$ are the pullbacks of $\mathcal{F}$ and $\mathcal{K}$ to $X'$, respectively to $\spec k$.
\end{lemma}

\begin{proof}
    Note that $p^{\ast}\mathcal{K}$ is flat over $Y$ since $p$ is flat.\ The isomorphism $\iota^{\ast}\mathcal{H}^d(G_p)$ is given as 
    \begin{align*}
        \mathcal{H}^d(L\iota^{\ast}Rp_{\ast}R\mathcal{H}om(\mathcal{F},p^{\ast}\mathcal{K} \otimes \omega_p)) \overset{\sim}{\longrightarrow} \mathcal{H}^d(L\iota^{\ast}R\mathcal{H}om(p_{\ast}\mathcal{F}, \mathcal{K})[-d]) \cong  \mathcal{H}om_k(\iota^{\ast}p_{\ast}\mathcal{F},K),
    \end{align*}
    where $\iota^{\ast}p_{\ast}\mathcal{F} \cong p'_{\ast}F$ by assumption.\ This uses Lemma \ref{extbasechangeisnatural}.\ Precomposing with $\alpha$ and pulling $L\iota'^{\ast}$ into $R\mathcal{H}om$ gives the following isomorphism denoted by $\mathcal{H}^d(\overline{G_p})$:
    \begin{align*}
        \mathcal{H}^d(Rp'_{\ast}R\mathcal{H}om(L\iota'^{\ast}\mathcal{F},p'^{\ast}K \otimes \omega_{X'})) = \ext^d_{\mathcal{O}_{X'}}(F,p'^{\ast} K \otimes \omega_{X'}) \overset{\sim}{\longrightarrow} \homo_k(p'_{\ast}F,K).
    \end{align*}
    Here $L\iota'^{\ast}\mathcal{F} = \iota'^{\ast}\mathcal{F} = F$ by flatness of $\mathcal{F}$ and $p$. The claim is that $G_{p'} \cong \overline{G_p}$.\par 
    \noindent The existence of the functorial isomorphisms $G_p$ and $G_{p'}$ is equivalent to having adjunctions $Rp_{\ast} \dashv p^{!}$ and $Rp'_{\ast} \dashv p'^{!}$.\footnote{See also \cite[Sect.\ 25.4]{goertzwedhorn2}.\ As $p$ is smooth and proper, $p^{!} \cong p^{\times} \cong (-) \otimes \omega_p[d]$.}\ Therefore the claim reduces to proving that the (co-)units $\varepsilon: Rp_{\ast}p^{!} \Rightarrow \id$ and $\id \Rightarrow p^{!}Rp_{\ast}$ associated to G.V.-duality for $p$ restrict to the (co-)units for $p'$ via base change along $L\iota^{\ast}$ and $L\iota'^{\ast}p^{!} \cong p'^{!}L\iota^{\ast}$.\ This is precisely \cite[\href{https://stacks.math.columbia.edu/tag/0AWG}{Tag 0AWG}]{stacks-project}, and the argument is completely formal.
\end{proof}

\end{document}